\documentclass[12pt,a4paper]{article}
\usepackage[]{fontenc}
\usepackage[latin1]{inputenc}
 \usepackage[colorlinks]{hyperref}
\usepackage{a4,amsmath,amssymb,amsfonts,amsthm}
\usepackage{amsrefs}
\usepackage[all]{xy}

\def\mZ{\mathbb{Z}}
\def\mN{\mathbb{N}}           % parce que vraiment on les utilise trop souvent
\def\mQ{\mathbb{Q}}

\def\mG{\mathbb{G}}
\def\mP{\mathbb{P}}

\def\Spec{{\rm Spec}}
\def\Sec{{\rm Sec}}

\def\Hom{{\rm Hom}}

 \def\Tr{{\rm Tr}}
 
 \def\Gal{{\rm Gal}}
 
 \def\det{{{\rm det}}}

  \def\Spec{{\rm Spec}}
\def\Id{{\rm Id}}

\def\Hilb{{\rm Hilb}}

\def\Pic{{\rm Pic}}

\def\et{{\rm et}}
\def\sep{{\rm sep}}

\def\Gr{{\rm Gr}}
\def\ker{{\rm ker}}
\def\Mor{{\rm Mor}}
\def\rk{{\rm rk}}

\gdef\Im{{\rm Im}}

\def\beginProof{\par{\bf Proof. }}
 \def\endProof{${\qed}$\par\smallskip}

 \def\mQ{{\Bbb Q}}
 
 \def\mZ{{\Bbb Z}}
 
 \def\mN{{\Bbb N}}
 \def\mF{{\Bbb F}}
 
  \def\mG{{\Bbb G}}
 
 \def\CE{{\mathcal E}}
 \def\CJ{{\mathcal J}}
 \def\CH{{\mathcal H}}
 \def\CO{{\mathcal O}}
 \def\CA{{\mathcal A}}
 \def\CD{{\mathcal D}}
 \def\CB{{\mathcal B}}
  \def\CW{{\mathcal W}}
   \def\CU{{\mathcal U}}
 \def\CC{{\mathcal C}}

 \def\CM{{\mathcal M}}

\def\CN{{\mathcal N}}
 \def\wt#1{\widetilde{#1}}
 
 \def\refeq#1{(\ref{#1})}

\def\P1{{{\bf P}^1}}

\def\GL{{\rm GL}}

\def\red{{\rm red}}

\def\min{{\rm min}}

\def\max{{\rm max}}
\def\min{{\rm min}}

\def\hn{{\rm hn}}

\def\brk{{\overline{\rm rk}}}
\def\unr{{\rm unr}}
\def\NB{{\bf N.B.}}
\def\ab{{\rm ab}}
\def\CTC{{\rm CT}}
\def\sh{{\rm sh}}

\def\Lie{{\rm Lie}}
\def\HN{{\rm HN\,}}
\def\binf{{\rm binf}}
\def\twp{{^{(p)}}}
\def\twpn{{^{(p^n)}}}
\def\IVD{{\rm IVD}}
\def\coLie{{\rm coLie}}
\def\order{{\rm order}}
\def\sm{{\rm sm}}

\def\Tor{{\rm Tor}}

\def\perf{{\rm perf}}

\def\c1{{\rm c}_1}

 \newtheorem{theor}{Theorem}[section]
 \newtheorem{prop}[theor]{Proposition}

 \newtheorem{lemma}[theor]{Lemma}

 \newtheorem{conj}[theor]{Conjecture}

\newtheorem{rem}[theor]{Remark}

\begin{document}

 \author{Damian R\"OSSLER\footnote{Mathematical Institute, 
University of Oxford, 
Andrew Wiles Building, 
Radcliffe Observatory Quarter, 
Woodstock Road, 
Oxford OX2 6GG, 
United Kingdom}}
 \title{ On the group of purely inseparable points of an abelian variety defined over a function field of positive characteristic II}
\maketitle
\begin{abstract}
% \begin{center}{\bf\Large PRELIMINARY VERSION; PLEASE DO NOT DISTRIBUTE}\end{center}

Let $A$ be an abelian variety over the function field  $K$ of a curve over a finite field. 
We describe several mild geometric conditions ensuring that the group $A(K^\perf)$ is finitely generated and 
that the $p$-primary torsion subgroup of $A(K^\sep)$ is finite. 
This gives partial answers to questions of Scanlon, Ghioca and Moosa, and Poonen and Voloch. We also describe a simple theory (used to prove our results) relating the Harder-Narasimhan filtration of vector bundles to the structure 
of finite flat group schemes of height one over projective curves over perfect fields. Finally, we  
use our results to give a complete proof of a conjecture of Esnault and Langer on 
Verschiebung divisibility of points in abelian varieties over function fields.
\end{abstract}

\tableofcontents

\section{Introduction}

Let $k$ be a finite field characteristic $p>0$ and let 
$S$ be a smooth, projective and geometrically connected curve over $k$. Let $K:=\kappa(S)$ be its function field. Let $A$ be an abelian variety of dimension $g$ over $K$. Choose an algebraic closure $\bar K$ of $K$. Let $K^\perf\subseteq\bar K$ be the maximal purely inseparable extension of $K$, let $K^\sep\subseteq\bar K$ be the maximal separable extension of $K$ and let $K^\unr\subseteq K^\sep$ be the maximal separable 
extension of $K$, which is unramified above every place of $K$. Finally, we let 
$\CA$ be a smooth commutative group scheme over $S$ such that $\CA_K=A$. 
We shall write $\omega_{\CA}:=\epsilon_{\CA/S}^*(\Omega_{\CA/S})$ for the restriction 
of the cotangent sheaf of $\CA$ over $S$ via the zero section 
$\epsilon_{\CA/S}:S\to\CA$ of $\CA$. We shall say that $\omega_{\CA}$ is the 
Hodge bundle of $\CA$. 

If $G$ is an abelian group, we shall write
$$
\Tor_p(G):=\{x\in G\,|\,\exists n\geq 0:p^n\cdot x=0\}
$$
and
$$
\Tor^p(G):=\{x\in G\,|\,\exists n\geq 0:n\cdot x=0\,\wedge\,(n,p)=1\}.
$$

The aim of this text is to prove the following two theorems and to give a proof 
of a conjecture of Esnault and Langer (see further below).

\begin{theor}
{\rm (a)} Suppose that $A$ is geometrically simple. If $A(K^\perf)$ is finitely generated and of \mbox{rank $>0$} then $\Tor_p(A(K^\sep))$ is a finite group.

{\rm (b)} Suppose that $A$ is an ordinary (not necessarily simple) abelian variety. If $\Tor_p(A(K^\sep))$ is a finite group then $A(K^\perf)$ is finitely generated.
\label{THA}
\end{theor}

\begin{theor}
Suppose that $\CA$ is a semiabelian scheme and that $A$ is a geometrically simple abelian 
variety over $K$.  
If $\Tor_p(A(K^\sep))$ is infinite, then 
\begin{itemize}
\item[{\rm (a)}] $\CA$ is an abelian scheme;
\item[{\rm (b)}]  there is $r_A\geq 0$ such that $p^{r_A}\cdot \Tor_p(A(K^\sep))\subseteq\Tor_p(A(K^\unr)).$
\end{itemize}
Furthermore, there is
\begin{itemize}
\item[{\rm (c)}]  an abelian scheme $\CB$ over $S$;
\item[{\rm (d)}]  an $S$-isogeny $\CA\to \CB$, whose degree is a power of $p$ and such that 
the corresponding isogeny $\CA_K\to\CB_K$ is \'etale;
\item[{\rm (e)}] an \'etale $S$-isogeny $\CB\to \CB$ whose degree is $>1$ and is a power of $p$,
\end{itemize}
and 
\begin{itemize}
\item[{\rm (f)}]  {\rm (Voloch)} if $A$ is ordinary then the Kodaira-Spencer rank of $A$ is not maximal;
\item[{\rm (g)}]  if $\dim(A)\leqslant 2$ then $\Tr_{\bar K|\bar k}(A_{\bar K})\not=0$;
\item[{\rm (h)}] for all closed points $s\in S$, the $p$-rank of $\CA_s$ is $>0$.
\end{itemize}
\label{THB}
\end{theor}

Here $\Tr_{\bar K|\bar k}(A_{\bar K})$ is the $\bar K|\bar k$-trace of $A_{\bar K}$. 
This is an abelian variety over $\bar k$. See subsection \ref{trss}. 

Theorems \ref{THA} and \ref{THB} (b) have applications in the context of the work of 
Poonen and Voloch on the Brauer-Manin obstruction over function fields. In particular Theorems \ref{THA} and \ref{THB} (b) show that the conclusion of \cite[Th. B]{Poonen-Voloch-BM} holds whenever the underlying abelian variety 
is geometrically simple, has semistable reduction and violates any of the conditions in Theorem \ref{THB}, in particular if it has a point of bad reduction. Theorems \ref{THA} and \ref{THB} (b) also feed into the "full" Mordell-Lang conjecture. See \cite[after Claim 4.4]{Scanlon-A-positive} and  \cite[Intro.]{AV-Toward} for this conjecture.  In particular, in conjunction with the main result of \cite{Ghioca-Moosa-Division} Theorems \ref{THA} and \ref{THB} (b) show that the "full" Mordell-Lang conjecture holds if the underlying abelian variety 
is ordinary, geometrically simple, has semistable reduction and violates any of the conditions in Theorem \ref{THB}, in particular if it has a point of bad reduction.

Let now $L$ be a field, which is finitely generated as a field over an algebraically closed field $l_0$ of characteristic $p$. Let 
$C$ be an abelian variety over $L$. 
 
\begin{conj}[Esnault-Langer]
Suppose that for all $\ell\geqslant 0$ we are given a point 
$x_\ell\in C^{(p^\ell)}(L)$ and suppose that for all $\ell\geqslant 1$, we have $V_{C^{(p^\ell)}/L}(x_{\ell})=x_{\ell-1}$. 
Then the image of $x_0$ in $C(L)/\Tr_{L|l_0}(C)(l_0)$ is a torsion point, which is of order prime to $p$.
\label{ELconj}
\end{conj}

See \cite[Rem. 6.3 and after Lemma 6.5]{Esnault-Langer-On-a-positive}. 
This conjecture is important in the theory of stratified bundles in positive characteristic; see \cite[Question 3 in the introduction]{Esnault-Langer-On-a-positive} for details. 

Here $C^{(p^\ell)}$ is the base change of $C$ by the $\ell$-th power of the 
absolute Frobenius morphism on $\Spec\,L$ and $V_{C^{(p^\ell)}/L}:C^{(p^\ell)}\to C^{(p^{\ell-1})}$ is the Verschiebung morphism. The abelian variety $\Tr_{L|l_0}(C)$ 
is the $L|l_0$-trace of $C$ (see subsection \ref{trss}). It is an abelian variety over $l_0$ and the variety $\Tr_{L|l_0}(C)_L$
 comes with an injective morphism to $C$. This gives in particular 
 an injective map $\Tr_{L|l_0}(C)(l_0)\to C(L)$. The Lang-N\'eron theorem (see \cite[chap. 6, Th. 2]{Lang-Fund}) asserts that $C(L)/\Tr_{L|l_0}(C)(l_0)$ is a finitely generated group. 
 Thus $\Tr_{L|l_0}(C)(l_0)\subseteq C(L)$ is precisely the subgroup of $C(L)$ consisting of 
 divisible elements (ie elements divisible by any integer).
 
 In the present text, we shall call a point $x_0\in C(L)$ with the property described in 
 Conjecture \ref{ELconj} an {\it indefinitely Verschiebung divisible point.} 
 We shall write $\IVD(C)=\IVD(C,L)\subseteq C(L)$ for the subgroup of 
 indefinitely Verschiebung divisible points. 
 
 We prove:
 
 \begin{theor} Conjecture \ref{ELconj} holds.\label{ELth}\end{theor}

Note that 
Theorem \ref{ELth} has the following consequence, which is of independent interest: 
if $C$ is as in Conjecture \ref{ELconj}, $C$ is ordinary and 
$\Tr_{L^\perf|l_0}(C_{L^\perf})=0$ then 
$$
\bigcap_{j\geq 0}p^j\cdot C(L^\perf)=\Tor^p(C(L^\perf)).
$$
To see this, let $x\in C(L^\perf)$. Let $L_1|L$ be a finite purely inseparable extension, which is a field of definition for $x$. Remember that the multiplication by $p$ endomorphism of $C$ is the composition of the Verschiebung morphism with the relative Frobenius morphism, which is purely inseparable. Also, recall that since $C$ is ordinary, the Verschiebung morphism is (by definition) separable. Note finally that since $\Tr_{L^\perf|l_0}(C_{L^\perf})=0$ we also have \mbox{$\Tr_{L_1|l_0}(C_{L_1})=0$.} In particular, if 
$x\in \bigcap_{j\geq 0}p^j\cdot C(L^\perf)$ then $x$ is an indefinitely Verschiebung divisible element of 
$C(L_1)$ and thus must lie in $\Tor^p(C(L_1))\subseteq\Tor^p(C(L^\perf))$ according to Theorem \ref{ELth}. The inclusion $
\Tor^p(C(L^\perf))\subseteq\bigcap_{j\geq 0}p^j\cdot C(L^\perf)
$ is straightforward. 

{\bf Outline of the paper.} The basic strategy of the paper hinges on  
Lemma \ref{lemcansub} below. This Lemma associates a maximal multiplicative subgroup scheme with any finite flat group scheme of height one over $S$. The existence of this subgroup scheme is not straightforward and follows from an analysis of the Harder-Narasimhan filtration of (a Frobenius twist of) the coLie algebra of the group scheme. This analysis is carried out in subsection \ref{sssHone}. 

One can apply Lemma \ref{lemcansub} to the kernel of the relative Frobenius morphism $F_{\CA/S}:\CA\to\CA^{(p)}$, replace $\CA$ by the resulting quotient and repeat this construction ad infinitum, stopping only when the maximal multiplicative subgroup scheme is trivial. 

It is then a basic (unresolved) question to determine minimal geometric conditions on $\CA$ ensuring that the resulting sequence of semiabelian schemes stops. This  also 
makes sense (and seems important to us) if $k$ is replaced by any perfect field 
of characteristic $p>0$ (not only when $k$ is finite). 

This question turns out to be intimately related to Theorems \ref{THA}, \ref{THB} and \ref{ELth}. To explain why, we shall first quote 
a result, which improves on (and elucidates) Lemma \ref{flem} in the Appendix. 
This result is proven in \cite{Rossler-Selmer}, which builds on the present article. We shall only need 
Lemma \ref{flem} in the present text but for conceptual clarity, we shall 
present the improved result in this outline. 
Let $E\subseteq S$ be the finite set of points $s\in S$ where $\CA_s$ is not an abelian variety. Let $U:=S\backslash E$. We first recall a classical result:
\begin{theor}[Artin-Milne] There is a canonical injective group homomorphism
$$
A^{(p)}(K)/F_{A/K}(A(K))\hookrightarrow\Hom_K(F_K^*(\omega_K),\Omega_{K/k}).
$$
\end{theor} 
Here $F_K$ is the absolute Frobenius endomorphism of $K$ (the $p$-th power map).
See  \cite[III.3.5.6]{AM-Duality} for the proof, which works in a more general setting. 
In \cite{Rossler-Selmer} this is refined as follows:
\begin{theor}[R.]
The image of the Artin-Milne map lies inside the subgroup 
$
\Hom_C(F_S^*(\omega),\Omega_{S/k}(E))
$
of
$
\Hom_K(F_K^*(\omega_K),\Omega_{K/k}).$
\label{RAM}
\end{theor}
Here $F_S$ is the absolute Frobenius endomorphism of $S$. Here we write
$\Omega_{S/k}(E):=\Omega_{S/k}(E)\otimes\CO_S(E)$ and $E$ is understood as a divisor with no multiplicities. 
Theorem \ref{RAM} 
refines Lemma \ref{flem} below (for the knowledgeable reader, in \cite{Rossler-Selmer} it is even proven 
that the image of the Selmer group of the relative Frobenius morphism lies in $
\Hom_C(F_S^*(\omega),\Omega_{S/k}(E))
$). This theorem is proven by providing a geometric interpretation for the Artin-Milne map and analysing its poles, making essential use of Faltings-Chai's semistable compactification of the universal abelian scheme. The existence of this compactification allows us to show that the poles are at most logarithmic, which is in essence the content of Theorem \ref{RAM}. Let us now explain why Theorem \ref{RAM} 
is relevant for Theorem \ref{THA}. Consider eg (b) in Theorem \ref{THA}. 
Suppose that $A(K^\perf)$ is not finitely generated. We have
$$
A(K^\perf)=\bigcup_{i\geq 0}A(K^{p^{-i}})
$$
and by the Lang-N\'eron theorem (see also subsection \ref{trss}) $A(K^{p^{-i}})$ is finitely 
generated. Hence for infinitely many $i\geq 0$, we must have 
$$
A^{(p^{i+1})}(K)/F_{A^{(p^i)}/K}(A(K))\simeq A(K^{p^{-i-1}})/A(K^{p^{-i}})\not=0.
$$
In particular, for infinitely many $i\geq 0$, we must have 
$$
\Hom_C(F_S^{\circ(i+1),*}(\omega),\Omega_{S/k}(E))\not=0
$$ according to Theorem \ref{RAM}. If now the vector bundle $\omega$ were 
ample, this would lead to a contradiction, because if $i$ is large enough 
and $\omega$ is ample then there cannot be any morphism from 
$F_S^{\circ(i+1),*}(\omega)$ to $\Omega_{S/k}(E)$. This was already noticed in the earlier article \cite{Rossler-On-the-group}, where details are given. One can refine this line of reasoning as follows. 
If $\omega$ is not ample and $A$ is ordinary then one can show 
that $\omega$ must have a certain non trivial quotient, which is semistable of degree $0$.  This non trivial quotient turns out to be induced by the maximal 
multiplicative subgroup scheme mentioned above. Calling it $G_\CA$, we may then replace 
$\CA$ by $\CA/G_\CA$. The group  $(\CA/G_\CA)_K(K^\perf)$ will again 
be infinitely generated, since the morphism $A\to (\CA/G_\CA)_K$ has finite kernel. Hence we can repeat the above reasoning for $\CA/G_\CA$ and we obtain 
an infinite sequence of isogenous abelian varieties. The next step in the proof of Theorem \ref{THA} (b) is to show that in this sequence, there are finitely many isomorphism classes. This follows from the fact that the degrees 
of $\omega_\CA$ and $\CA/G_\CA$ are the same and more generally the degrees 
of the Hodge bundles of all the semiabelian schemes in the sequence are the same. 
This is a consequence of a computation involving the cotangent complex 
of the quotient morphism (see Lemma \ref{IDeglem}). 
It then follows from a classical reasoning involving moduli spaces of abelian 
varieties, familiar from Zarhin's proof of the Tate conjecture over function fields, that 
the sequence contains only finitely many isomorphism classes. We can thus 
conclude that, up to isogeny, $A$ contains a non trivial finite endomorphism, whose kernel is multiplicative. The dual of this endomorphism is then separable and this 
shows that $\Tor_p(A^\vee)(K^\sep)$ is infinite (consider the 
kernels of its powers). Since $A^\vee$ is isogenous 
to $A$, we see that $\Tor_p(A)(K^\sep)$ is also infinite. This concludes our outline of the proof of Theorem \ref{THA} (b). 

For Theorem \ref{THA} (a), we consider the quotients of $A$ by finite subgroups of 
$\Tor_p(A)(K^\sep)$ of increasing size. These quotients also run through finitely many isomorphism classes by a similar reasoning and we thus see that 
if $\Tor_p(A)(K^\sep)$ is infinite then, up to isogeny, $A$ is endowed with a separable finite endomorphism. The dual of this endomorphism is then purely inseparable and of degree a positive power of $p$, and if $A^\vee(K)$ is not finite, we may show that $A^\vee(K^\perf)$ is infinitely generated by considering the inverse images of $A(K)$ under the powers of this endomorphism. If now $A^\vee(K^\perf)$ is not finitely generated, neither 
is $A(K^\perf)$, since $A$ and $A^\vee$ are isogenous. This concludes our outline of the proof of Theorem \ref{THA} (a). 

In Theorem \ref{THB}, we start out as in Theorem \ref{THA} (a) and we again obtain, up to isogeny, a separable finite endomorphism of degree a positive power of $p$. The rest of the theorem 
investigates the geometric consequences of the existence of this endomorphism. 
The most interesting consequence is the fact that it implies that 
$\CA$ must be an abelian scheme (if $A$ is geometrically simple). This is (a) in Theorem \ref{THB}). 
The main point here is that the endomorphism extends to an \'etale endomorphism of $\CA$. If $\CA$ had a fibre with a toric part then the endomorphism would 
induce an automorphism of the toric part, because tori only have 
infinitesimal $p$-primary subgroups in characteristic $p$ and these are only \'etale if they are trivial. This fact forces the whole endomorphism to be an automorphism, which is impossible. The proof of (c), (d) and (e) are straightforward and not much more than a rewording of the fact that there are only finitely many isomorphism classes in the set of quotients described above. The proof of (b) 
follows essentially from a variant of the fact that, under (a), the above endomorphism extends to an everywhere \'etale and finite endomorphism of $\CA$. 
This also easily gives a proof of (h).  
The proof of (g) is based on class field theory and the Serre-Tate theory of canonical liftings. First, up to a finite extension, 
the field extension generated by the points of $\Tor_p(A)(K^\sep)$ is everywhere unramified by (a) and (b). 
If $\Tor_p(A)(K^\sep)=\Tor_p(A)(\bar K)$ then a simple application of the Serre-Tate theory of canonical liftings shows 
that $A_{\bar K}$ is the base change of an abelian variety defined over  $\bar k.$
Hence it must be contained in the Hilbert class field of $K$, which is but a constant field extension (ie comes from an extension of $k$), up to a finite extension. So if $\Tor_p(A)(K^\sep)$ is infinite then it is an infinite torsion subset of $A(K\bar k)$, which is finitely generated by the Lang-N\'eron theorem if the trace of $A$ vanishes: contradiction.

We now turn to Theorem \ref{ELth}. Using a height argument due to Raynaud, Esnault and Langer prove in \cite[Th. 6.2]{Esnault-Langer-On-a-positive} that the image of $x_0$ in $C(L)/\Tr_{L|l_0}(C)(l_0)$ is a torsion point under the assumption that $C$ 
has everywhere potential good reduction in codimension one. 
Their argument works as follows. Choose a polarisation on $C$. This induces polarisations on all the $C^{(p^\ell)}$ by base change. A simple computation shows 
that if a point $x\in C(L)$ has a preimage $y$ in $C^{(p)}(K)$ under 
the Verschiebung map then the height of $x$ with respect to the polarisation 
is $p$ times the height of $y$ with respect to the base changed polarisation. 
Now if $C$ has everywhere good reduction in codimension one, there is 
an abelian scheme $\CC$ extending $C$ on an open subset with complement of codimension $\geq 2$ of a normal complete model $V$ of $L$ and the polarisations 
on $C$ and $C^{(p)}$ naturally extend to this open subset. This implies that the heights of $x$ are $y$ (with respect to the polarisations and a choice of ample line bundle on $V$) are integers, because they can then computed in a completely geometric fashion. In particular, the height of $x$ is an integer divisible by $p$. Repeating this argument with $y$, one sees that the height of $x$ is divisible by arbitrarily high powers of $p$ and one concludes that it must vanish. Then the conclusion 
follows from a theorem of Lang (see \cite[Th. 9.15]{Conrad-Trace}). 
The argument described above breaks down in the presence of bad reduction in codimension one because the orders of the component groups of 
the special fibres of the local N\'eron models of the varieties $C^{(p^\ell)}$ increase with $\ell$ if they are not trivial and this introduces denominators in the heights. 

Our approach to Theorem \ref{ELth} is again via the infinite sequence of quotients described at the beginning of the outline. This sequence will effectively 
replace the sequence of the $C^{(p^\ell)}$. It has the advantage over 
the sequence of the $C^{(p^\ell)}$ that it falls inside a bounded family of abelian varieties (see below), making it possible to control the order of the (analogues of the) images 
of the $x_{\ell}$ in the component groups of the N\'eron models. This makes 
a similar height computation possible. The proof is in several steps. 

Step (0). Reduction to the case where $L$ is the function field of 
a smooth and projective curve $B$ over $l_0$. This follows from 
a Bertini type argument - see section \ref{secRMW} in the Appendix. 

Step (1). We consider the images of the $x_\ell$ under the Artin-Milne map. 
A crucial point is that these images must be compatible under the Verschiebung 
morphisms (see diagram \refeq{icd} below) and this constrains the image of $x_1$ under the Artin-Milne map. Using Lemma \ref{flem} (or Theorem \ref{RAM}), the theory of semistable sheaves in positive characteristic and various global results on finite flat group schemes of height one in a global situation proven in section \ref{filtsec}, we show that the image of $x_1$ under the Artin-Milne map must factor through the coLie algebra of the maximal multiplicative subgroup $(\ker\,F_{\CC/B})_\mu$ of $\ker\,F_{\CC/B}$. This implies that
the image of $x_1$ in $(C^{(p)}/(\ker\,F_{\CC/B})^{(p)}_{\mu,L})(L)=
(C^{(p)}/G_{\CC^{(p)},L})(L)$ 
maps to $0$ under the Artin-Milne map. From the definitions, this means that the 
image of $x_0$ in $(C/G_{\CC,L})(L)$ is divisible by $p$ in 
$(C/G_{\CC,L})(L)$. 
Suppose for simplicity that $C$ has a semiabelian model $\CC$ over $B$. We can now repeat this process and we obtain a sequence of purely 
inseparable morphisms $\psi_{i}:\CC\to \CC_i$ of increasing degree, such that 
$\psi_{i,L}(x_0)$ in $C_i$ is divisible by $p^i$ in $C_i(L)$.

Step (2). We choose a polarisation $\phi_{D_0}:C\to C^\vee$. The image 
of $x_0$ under $\phi_{D_0}$ is of course also indefinitely Verschiebung divisible. 
We identify $\phi_{D_0}(x_0)$ with a line bundle $M$ on $C$. 
Since $\phi_{D_0}(x_0)$ is indefinitely Verschiebung divisible, 
there are line bundles $M_i$ on $C^{(p^i)}$ such that $M$ 
is the pull-back of $M_i$ by the morphism $C\to C^{(p^i)}$ arising 
by composing relative Frobenii. The morphism $C\to C^{(p^i)}$ 
factors through $\psi_{i,L}$ by construction. Hence there are line bundles 
$J_i$ on the $C_i$ such that $\psi_{i,L}^*(J_i)=M$. 

Step (3). We now compute the height pairing between 
$x_0$ and $M$. This can easily be seen to equal
the height pairing between $\psi_{i,L}(x_0)$ and $J_i$.
Since 
$\psi_{i,L}(x_0)$ is divisible by $p^i$, we see that the height pairing between 
$x_0$ and $M$ is divisible by $p^i$. If the $\CC_i$ were all abelian schemes 
we could deduce (like Raynaud-Esnault-Langer above) that the height pairing between $x_0$ and $M$ must vanish, because then all the values of the various height pairing would be integral. However, we cannot assume this. 

Step (4). All the $C_i$ are essentially part of a bounded family of abelian varieties over $L$ because the degrees of the Hodge bundles of the $\CC_i$ are all equal 
(see above in the outline). Using this, one can prove that there is an infinite set $I_0\subseteq\mN$ such that if $i\in I_0$ the image of any element 
of $C_i(L)$ in the component groups of the N\'eron model of $C_i$ 
has an order, which is bounded independently of $i$. This follows from Proposition \ref{GUprop} (a) in the appendix. The gist of the argument is that in a bounded family of semiabelian varieties over $B$, it is possible 
to smoothly compactify the generic fibre, up to to normalisation in a finite extension of the function field of the parameter space. This would 
follow from resolution of singularities but in the present situation is a consequence of the work of Mumford, Chai-Faltings and K\"unnemann (see \cite[Th. 4.2]{Ku-Proj}). This means that the abelian varieties in the family almost all have regular compactifications with a bounded number of geometric fibres over $B$. 
This bound is also a bound for the order of the image of a rational point in the component groups of the N\'eron model.

Step (5). In view of Step (4), if we replace $x_0$ by a certain multiple of $x_0$, all the height pairing in sight are integers. Hence the divisibility argument envisaged in Step (3) can be carried out and yields that the height 
pairing of $x_0$ and $M$ vanishes. This pairing is by construction twice the N\'eron-Tate height of $x_0$ with respect to the polarisation $\phi_{D_0}$ and we conclude from a theorem of Lang (op. cit.) that the image of $x_0$ in $C(L)/\Tr_{L|l_0}(C)(l_0)$ is a torsion point. It remains to show that its order is prime to $p$. 

Step (6).  We first show that we may suppose that $\Tr_{\bar L|l_0}(C_{\bar L})=0$. This is not completely 
straightforward, because when one passes to a finite extension in Conjecture \ref{ELconj}, one loses control of part of the 
torsion of $C(L)/\Tr_{L|l_0}(C)(l_0)$. However, although the parasitical torsion subgroup that might appear is 
not known, its exponent only depends on the degree of the extension. This degree can be taken to be the same for all the 
Frobenius twists of $C$ and the information one gathers from this suffices to prove the conjecture, provided one can prove it for a finite extension. Thus we may suppose that $\dim(\Tr_{\bar L|l_0}(C_{\bar L}))=\dim(\Tr_{L|l_0}(C))$ and then, after quotienting by $\Tr_{L|l_0}(C)$, that $\Tr_{\bar L|l_0}(C_{\bar L})=0$. Now recall that the $C_i$ are essentially part of a bounded family of abelian varieties over $L$ (see step (4)). Using this, and the fact that now $\Tr_{\bar L|l_0}(C_{i,\bar L})=0$ for all $i\geq 1$, one can prove that there is an infinite set $I_0\subseteq\mN$ such that if $i\in I_0$, the cardinality of the torsion subgroup of 
$C_i(L)$ is uniformly bounded. This follows from Proposition \ref{GUprop} (b) in the Appendix. To finish the proof of Conjecture \ref{ELth}, suppose that $x_0$ is a non-zero torsion point, which is indefinitely Verschiebung divisible. Since the image of $x_0$ in $C_i(L)$ is divisible by $p^i$, we see that 
the torsion group of $C_i(L)$ has an element of order $p^{i+1}$. This contradicts the above uniformity statement and 
shows that the order of $x_0$ must be prime to $p$. 

The argument to prove the uniformity statement alluded to in Step (6) goes roughly as follows. 
One first notices that the torsion subgroup of a trace free abelian variety coincides with the set of elements of 
vanishing N\'eron-Tate height by the already quoted theorem of Lang. Thus they can be described as the points of a moduli 
space of sections, which is of finite type over $l_0$, at least for those torsion points, whose image in the component 
groups of the N\'eron model of the abelian variety is trivial. Since the abelian variety is trace free, the torsion subgroup is finite and thus this moduli space is finite. Using the uniformity statement in Step (4), we may assume that the torsion points of the $C_i(L)$ have trivial images in 
the components of the corresponding N\'eron models, up to multiplication by a 
fixed integer (independent of $i$ running through an infinite set).  The number of irreducible components of the moduli space 
of each $C_i$ is now uniformly bounded, since the $C_i$ are part of a bounded family. This gives a uniform 
bound for the torsion subgroups of the $C_i(L)$. 

The reader may enjoy the talk \cite{Rossler-IHP-video} as an introduction to parts of the present article.

%A final remark is 
%that in the course of our proof of Theorem \ref{ELth}, we show that to prove Conjecture 
%\ref{ELconj} in general, it is sufficient to prove it in the situation where $L$ has transcendence degree one over $l_0$ (for any algebraically closed field 
%$l_0$ of characteristic $p>0$). See section \ref{secELproof}, reduction (1). 

The structure of the article is as follows. In section \ref{sIR}, we state various intermediate 
results, from which we shall deduce Theorems \ref{THA} and \ref{THB}. 
Theorem \ref{mainthm1} in subsection \ref{ssIGP} is of independent interest and 
is (we feel) likely to be useful for the study of the geometry of (especially ordinary) abelian varieties in 
general. The results in subsection \ref{ssIGP} are deduced from some results in the theory of 
finite flat groups schemes of height one over $S$, most of which follow from the existence of a Harder-Narasimhan filtration on their Lie algebras. These results on finite flat group schemes are proven in 
section \ref{filtsec} and for the convenience of the reader, we included a section (section 
\ref{prolsec}) listing the results on semistable sheaves over curves in 
positive characteristic that we need. To the knowledge of the author, there are 
very few general results on the structure of finite flat group schemes in a global situation 
(eg when the base is not affine) and it seems that it is the first time time that 
the theory of semistability of vector bundles is being used in this context. In \cite{Bost-Dwork} a similar idea is used in characteristic $0$, where it is applied to the study of formal groups over curves (recall that all groups schemes are smooth in characteristic $0$, so the Lie algebras of finite flat group schemes vanish in characteristic $0$). Lemma \ref{SBlem} below (which concerns finite flat group schemes of height one) is inspired by \cite[Lemma 2.9]{Bost-Dwork}. A prototype 
of Lemma \ref{SBlem} can be found in \cite[Lemma 9.1.3.1]{SB-generic} but it is not applied to the study of group schemes there. 
The key results here are the  Lemmata \ref{SBlem} and \ref{lemcansub}, which will 
hopefully lead to further generalisations (eg in the situation when the base scheme is 
of dimension higher than one - in this direction, see \cite[Th. 7.3]{Langer-Generic}). The results in subsection \ref{ssIGKS} do not require 
the theory of semistable sheaves and are based on geometric class field theory, the theory of Serre-Tate canonical liftings and 
on the existence of moduli schemes for abelian varieties. In section \ref{finsec}, we prove
the various claims made in subsection \ref{ssIGP} and in section \ref{pcIGKS} we prove the claims 
made in subsection \ref{ssIGKS}. In section \ref{pTHA}, we prove Theorem \ref{THA} and 
in section \ref{pTHB} we prove Theorem \ref{THB}. In section \ref{secELproof}, we give a proof of Theorem \ref{ELth}. 
The proof of Theorem \ref{ELth} is quite long and uses virtually all the other results proven in this text. 

In his very interesting recent preprint \cite{XY-Pos},  Xinyi Yuan uses some techniques 
which are also used in the present paper. They were discovered independently. His text focusses on the case where the base curve is the projective line.  
In particular, the "quotient process" used in step (2) of the proof of Theorem \ref{ELth} and also in the proof of Theorem \ref{THB} 
also appears (over the projective line) in section 2.2 of \cite{XY-Pos}. Theorem 2.9 of \cite{XY-Pos} overlaps with the proof of Lemma \ref{Slem}.

The prerequisites for this article are algebraic geometry at the level of the EGA, familiarity with the basic theory of finite flat group schemes, as expounded in \cite{Tate-Finite} and a good knowledge of the theory of abelian schemes and varieties, as presented in \cite{Milne-Abelian}, \cite{Mumford-Abelian} and 
\cite{MB-Pinceaux}. We also expect the reader to be familiar with the basic properties of N\'eron models (as in the chapter on basics of \cite{Bosch-Raynaud-Neron}) and to have a working knowledge of Grothendieck topologies.

%Finally, in section \ref{pCEL}, we prove Theorem \ref{allC}, as a straightforward 
%consequence of Theorems \ref{THA} and \ref{THB} and of results already available in the literature.

%Appendix \ref{appendix} 
%provides a proof of a generalisation of some results proven in 
%the earlier article \cite{Rossler-On-the-group} of the author, which are crucial for the proof of Theorem \ref{THA} (b) and Theorem \ref{ELth}.

{\bf Acknowledgments.} 
My warm thanks to the referee for his careful reading and for many suggestions. The article would be much less clear without his help and encouragement. I would like to thank J.-B. Bost for his feedback, especially for pointing out the article 
\cite{Catanese-Dettweiler-Vector}, for suggesting Remark \ref{remCat} and for providing 
\cite[Lemma 2.9]{Bost-Dwork}, whose positive characteristic analogue is technically at the root of the present text. 
Minhyong Kim's article \cite{Kim-Purely} also played a fundamental role in the genesis of the present text; the construction described there pointed me in (what I hope is) the right direction when I started studying purely inseparable points on abelian varieties. I had many interesting discussions with him about his article. 
I am very grateful to J.-F. Voloch for many exchanges on the material of this article and for his remarks on the text and to P. Ziegler 
for many discussions on and around the "full" Mordell-Lang conjecture. Many thanks also to T. Scanlon for his interest and for interesting discussions 
around the group $A(K^\perf).$ Last but not least, many thanks to H\'el\`ene Esnault 
and her student Marco d'Addezio for their interest and for many enlightening discussions around Theorem \ref{ELth}. I also benefitted from  A.-J. de Jong's and F. Oort's vast knowledge; they both very kindly took the time 
to answer some rather speculative messages.

{\bf Notation.} If $X$ is an integral scheme, we write $\kappa(X)$ for the local ring at the generic point of $X$ (which is a field). If $X$ is a scheme of characteristic $p$, we denote the absolute Frobenius 
endomorphism of $X$ by $F_X$. If $f:X\to Y$ is a morphism between two schemes of characteristic $p$ and $\ell>0$, 
abusing language, we denote by $X^{(p^\ell)}$ the fibre product of 
$f$ and $F^{\circ \ell}_Y$, where $F^{\circ \ell}_Y$ is the  $\ell$-th power of the Frobenius 
endomorphism $F_Y$ of $Y$. If $G\to X$ is a group scheme, we write 
$\epsilon_{G/X}:X\to G$ for the zero section of $G$ and 
$$\omega_{G/X}=\omega_G:=\epsilon_{G/X}^*(\Omega_{G/X}).$$ If $X$ is of characteristic $p$, we shall 
write $F_{G/X}:G\to G\twp$ for the relative Frobenius morphism. If in addition $G$ is flat and commutative, we shall write 
$V_{G\twp/X}:G\twp\to G$ for the corresponding Verschiebung morphism; 
we shall write $F^{(n)}_{G/X}:G\to G\twpn$ (resp. $V^{(j)}_{G\twpn/X}:G\twpn\to G^{(p^{n-j})}$) for the composition 
of morphisms $$F_{G^{(p^{n-1})}/X}\circ\dots\circ F_{G/X}$$  
(resp. the composition of morphisms $$V_{G^{(p^{n-j+1})}/X}\circ V_{G^{(p^{n-j+2})}/X}\circ\dots
\circ V_{G^{(p^{n})}/X}$$).  See \cite[Exp. $\rm VII_A$, par. 4, "Frobeniuseries"]{SGA3-1} for the definition of the relative Frobenius morphism and the Verschiebung. If $G$ is finite flat and commutative, we shall write $G^\vee$ for the Cartier dual of $G$.

\section{Intermediate results}

We keep the notations and terminology of the introduction.

\label{sIR}

\subsection{Consequences of infinite generation of $A(K^\perf)$}

\label{ssIGP}

We shall write
$$
\brk_\min(\omega_\CA):=\lim_{\ell\to\infty}\rk((F^{\circ\ell,*}_S(\omega_\CA))_\min)
$$
and
$$
\bar\mu_\min(\omega_\CA):=\lim_{\ell\to\infty}{\deg((F^{\circ\ell,*}_S(\omega_\CA))_\min)
\over p^\ell\cdot \rk((F^{\circ\ell,*}_S(\omega_\CA))_\min)}.
$$
Here $F^{\circ\ell}_S$ is the $\ell$-th power of the 
absolute Frobenius endomorphism of $S$ and $(F^{\circ\ell,*}_S(\omega_\CA))_\min$ is the semistable quotient with minimal slope 
of the vector bundle $F^{\circ\ell,*}_S(\omega_\CA)$. See section \ref{prolsec} for details.
 Our main tool will be the following theorem.

\begin{theor} 
There exists a (necessarily unique) multiplicative subgroup scheme $G_\CA\hookrightarrow\ker\,F_{\CA/S}$, with the following property: if $H$ is a finite, flat, multiplicative group scheme of height one over $S$ and 
$f:H\to\ker\,F_{\CA/S}$ is a morphism of group schemes, then $f$ factors through $G_\CA$. 

If $A$ is ordinary and $\omega_\CA$ is not ample then 
the order of $G_\CA$ is $p^{\brk_\min(\omega_\CA)}$. 

If $\phi:\CA\to\CB$ is a morphism
of smooth commutative group schemes over $S$, then the restriction of $\phi$ to $G_\CA$ factors through 
$G_\CB$. Furthermore, we have  
$\deg(\omega_{\CA})=\deg(\omega_{\CA/G_\CA}).$
\label{mainthm1}
\end{theor}
Here $\CA/G_\CA$ is the "fppf quotient" of $\CA$ by $G$, which is also a smooth commutative group scheme over $S$. See Proposition \ref{propQ} below 
for details.

\begin{rem}\rm Note that $\bar\mu_\min(\omega_\CA)>0$ is equivalent to $\omega_\CA$ being ample (see \cite{Barton-Tensor}). \end{rem}

\begin{rem} \rm Theorem \ref{mainthm1} holds more generally if $k$ is only supposed to be perfect (the proof does not use the fact that $k$ is finite).\end{rem}

\begin{rem}\rm It would be interesting to provide an explicit example of an abelian variety $A$ as in the introduction to this article, such that 
$A$ is 
ordinary, $\CA$ is semiabelian, $\Tr_{\bar K|\bar k}(A_{\bar K})=0$ and $G_\CA\not=0$. 
It should be possible to construct such an example by considering mod $p$ reductions of the abelian variety constructed in 
\cite[Th. 1.3]{Catanese-Dettweiler-Vector}. We hope to return to this question in a later article. The following question is also 
of interest: is there an ordinary abelian variety $A$ as above, such that $A$ has maximal Kodaira-Spencer rank, $\CA$ is semiabelian and 
$G_\CA\not=0$?
\label{remCat}
\end{rem}

%\begin{rem}\rm The following questions are crucial in the context of Theorem \ref{mainthm1}. 
%
%Is there an example of a smooth commutative group scheme over $S$ satisfying the 
%assumptions of Theorem \ref{mainthm1}, such that $\bar\mu_\min(\omega_\CA)=0$ (equivalently: $\omega_\CA$ is not ample) and such that $\Tr_{\bar K|\bar k}(A_{\bar K})=0$? We shall see in Proposition \ref{mainth0} below that a counterexample to Conjecture \ref{conjSZ} would 
%provide such a group scheme.
%
%The characteristic zero analog of this question is the following. Is there an example 
%of a smooth commutative group scheme $\CB$ over a smooth proper curve $T$ over 
%the complex numbers $\mC$, such that 
%
%- the generic fibre $\CB_{\kappa}$ of $\CB$ is an abelian variety (here $\kappa=\kappa(T)$ is the function field of $T$),
%
%- $\Tr_{\bar{\kappa}|\mC}(\CB_{\bar{\kappa}})=0$,
%
%- $\omega_{\CB/T}$ is not ample?
%
%To the author's knowledge, this is an open question. Note that it is known that in this situation, $\omega_{\CB/T}$ is nef 
%(see \cite[Th. 5.2]{Griffiths-Periods-III}). In \cite[Cor. 2.7]{Bost-Germs} a proof of this fact is given, which is close in spirit to the proof of Theorem \ref{mainthm1} (the method is to show that a counterexample to nefness would imply the existence of an algebraisable isotrivial formal subgroup of $\CB$) and it would be very interesting 
%to have a better understanding of this analogy. For instance, Lemma \ref{lemid} seems to be structurally close 
%to a part of Theorem 4.2 in \cite{Bost-Germs}.
%\end{rem}

\begin{prop} Suppose that $A$ is ordinary and that $\CA$ is semiabelian. 
Suppose that $A(K^\perf)$ is not finitely generated. Then 
$G_\CA$ is of order $>1$ and $\CA/G_\CA$ is also semiabelian.
\label{mainth0}
\end{prop}

\begin{prop}  
Suppose that $A$ is ordinary and that $\CA$ is semiabelian over $S$. Suppose that $A(K^\perf)$ is not finitely generated. 

Then there is a finite flat morphism $$\phi:\CA\to\CB$$ where $\CB$ is a semiabelian over $S$ and 
a finite flat morphism $$\lambda:\CB\to\CB$$
such that $\ker(\phi)$ and $\ker(\lambda)$ are multiplicative group schemes and such that the order of 
$\ker(\lambda)$ is $>1$.
\label{mainth}
\end{prop}

%Note that by a theorem of Raynaud (see \cite[chap. 5, §5.2]{Abbes-RSSC} for this and further references), the connected component of 
%the  zero section of the N\'eron model of $A$ over $S$ is necessarily semiabelian if $A(K)[n]\simeq(\mZ/n\mZ)^{2\dim(A)}$ for some 
%$n\geq 3$ such that $(p,n)=1$. 

%\begin{theor}[Voloch]
%Suppose that $A$ is ordinary and that the Kodaira-Spencer class of $A$ over $K$ has maximal rank. 
%Then $\Tor_p(A(K^\sep))=0$.
%\label{volth}
%\end{theor}
%
%See \cite[§4, p. 1093]{Voloch-Dioph-p} for the proof and the precise definition of 
%the Kodaira-Spencer map. 
%In view of Corollary \ref{corUnr}, (****) above is now established.

\subsection{Consequences of infiniteness of $\Tor_p(A(K^\sep))$ or $\Tor_p(A(K^\unr))$}

\label{ssIGKS}

\begin{theor} Suppose that $\Tr_{\bar K|\bar k}(A_{\bar K})=0$. 
Suppose that the action of $\Gal(K^\sep|K)$ on $\Tor_p(A(K^\unr))$ factors 
through $\Gal(K^\sep|K)^\ab$. Then $\Tor_p(A(K^\unr))$ is finite.
\label{classmain}
\end{theor}
Here $\Gal(K^\sep|K)^\ab$ is the maximal abelian quotient of $\Gal(K^\sep|K)$.

\begin{prop} Suppose that $\dim(A)\leqslant 2$ and that $\Tr_{\bar K|\bar k}(A_{\bar K})=0$. 
Then $\Tor_p(A(K^\unr))$ is finite.
\label{classcor}
\end{prop}

\begin{theor}
Suppose that $\Tor_p(A(K^\sep))$ is infinite. Then there is an \'etale $K$-isogeny $$\phi:A\to B$$ where $B$ is an abelian variety over $K$ and there is an \'etale $K$-isogeny $$\lambda:B\to B$$
such that the order of 
$\ker(\lambda)$ is $>1$ and such that the orders of $\ker(\lambda)$ and $\ker(\phi)$ are powers of $p$.
\label{threp}
\end{theor}

\begin{theor}
Suppose that there exists an \'etale $K$-isogeny $\phi:A\to A$, such that 
$\deg(\phi)=p^r$ for some $r>0$. Suppose also that $A$ is a geometrically simple abelian variety and that $\CA$ is a semiabelian scheme.

Then $\CA$ is an abelian scheme and $\phi$ extends to an \'etale (necessarily finite) $S$-morphism $\CA\to \CA$ of group schemes. 
\label{thabet}
\end{theor}

%Proposition \ref{mainth0}, Proposition \ref{mainth} and Corollary \ref{corUnr} will be derived from mild generalisations of results proven in  (see the two appendices for these 
%generalisations) and from
%some results in the theory of finite flat group schemes over curves proven in section \ref{filtsec} below. 
%The proof of Proposition \ref{mainth} also uses the existence of moduli spaces 
%for abelian varieties over finite fields. 
%In section \ref{prolsec}, we recall various results in the theory of semistable sheaves, which 
%will be used in section \ref{filtsec}. In section \ref{finsec}, we give the proof 
%of Proposition \ref{mainth0}, Proposition \ref{mainth} and Corollary \ref{corUnr}, using almost all the lemmata proven in section \ref{filtsec}.
%
%In section \ref{Rossec}, we shall prove Theorem \ref{RosTh}. The method of the proof is a variant 
%of the method used to prove \cite[Th. 1.3]{Rossler-Infinitely}. 
%
%In section \ref{classsec}, we shall prove Theorem \ref{classmain} and Corollary 
%\ref{classcor}. The proof relies on geometric class field theory and the theory of Serre-Tate liftings. 
%
%In section \ref{Bizsec}, we shall prove Theorem \ref{Thbiz} and Theorem \ref{LEconj}, making heavy use of 
%Proposition \ref{mainth}.
%
%In the appendix \ref{appendix}, we prove mild extensions of the main results of \cite{Rossler-On-the-group}, making use 
%of Zarhin's trick to eliminate certain assumptions on polarisations. In the appendix 
%\ref{appendixII}, we prove a straightforward generalisation of a result proven in 
%\cite{Rossler-Infinitely}.

\section{Semistable sheaves on curves}

\label{prolsec}

Let $Y$ be a scheme, which is smooth, projective and geometrically connected of relative dimension one over a field $t_0$. 

Suppose to begin with that $t_0$ is algebraically closed. 

If $V$ is a non zero coherent  locally free sheaf on $Y$, we write as is customary
$$
\mu(V)=\deg(V)/\rk(V)
$$
where $$\deg(V):=\int_Y \c1(V)$$ and 
$\rk(V)$ is the rank of $V$. The quantity 
$\mu(V)$ is called the {\it slope} of $V$. 
Recall that a non zero locally free coherent sheaf $V$ on $Y$ is called 
semistable if for any non zero coherent subsheaf $W\subseteq V$, we have 
$\mu(W)\leqslant\mu(V)$. 
Let $V/Y$ be a non zero locally free coherent sheaf on $Y$. There is a unique filtration 
by coherent subsheaves
$$
0=V_{0}\subsetneq V_1\subsetneq V_2\subsetneq\dots\subsetneq V_{\hn(V)}=V
$$
such that all the sheaves $V_i/V_{i-1}$ ($1\leqslant i\leqslant\hn(V)$) are (locally free and) semistable and 
such that the sequence $\mu(V_i/V_{i-1})$ is strictly decreasing. 
This 
filtration is called the {\it Harder-Narasimhan filtration} of $V$ (shorthand: \HN filtration). One then defines 
$$
V_\min:=V/V_{\hn(V)-1},\,
V_\max(V):=V_1
$$
and
$$\mu_\max(V):=\mu(V_1),\,\mu_\min(V):=\mu(V_\min).
$$
Let now $r\in\mQ$. Suppose that $r\in\{\mu(V_1),\dots,\mu(V/V_{\hn(V)-1})\}$. Let $i(r)\in\mN$ be the unique natural number such that $\mu(V_{i(r)}/V_{i(r)-1})=r$. We shall write 
$$
V_{=r}:=V_{i(r)}/V_{i(r)-1}
$$
and
$$
V_{\geq r}:=V_{i(r)}.
$$
We shall also write
$$
V_{>r}:=V_{j(r)}
$$
where $j(r)\in\mN$ is the largest natural number such that $\mu(V_{j(r)}/V_{j(r)-1})>r$.

One basic property of semistable sheaves that we shall use repeatedly 
is the following. If 
$V$ and $W$ are non zero coherent locally free sheaves on $Y$ and $\mu_\min(V)>\mu_\max(W)$ then 
$\Hom_Y(V,W)=0$. This follows from the definitions. 

See  \cite[chap. 5]{Brenner-Herzog-Villamayor-Three} (for instance) for all these notions. 

If $V$ is a non zero coherent locally free sheaf on $Y$ and $t_0$ has positive characteristic, 
we say that $V$ is {\it Frobenius semistable} if $F^{\circ r,*}_Y(V)$ is semistable for all $r\in\mN$.  The terminology {\it strongly semistable} also appears in the literature. 

\begin{theor} Let $V$ be a non zero coherent locally free sheaf on $Y$. 
%The sequence $\mu_\min(F^{\circ\ell,*}_S(V))/p^\ell$ 
%(resp. $\mu_\max(F^{\circ\ell,*}_S(V))/p^\ell$) is 
%decreasing (resp. increasing). 
There is an $\ell_0=\ell_0(V)\in\mN$ such that the quotients of the Harder-Narasimhan filtration of $F^{\circ\ell_0,*}_Y(V)$ are 
all Frobenius semistable. 
\label{SSth}
\end{theor}
\beginProof
See eg \cite[Th. 2.7, p. 259]{Langer-Semistable}.
\endProof

Theorem \ref{SSth} shows in particular that the following definitions : 
$$
\bar\mu_\min(V):=\lim_{\ell\to\infty}\mu_\min(F^{\circ\ell,*}_Y(V))/p^\ell,
$$
$$
\bar\mu_\max(V):=\lim_{\ell\to\infty}\mu_\max(F^{\circ\ell,*}_Y(V))/p^\ell,
$$
$$
\brk_\min(V):=\lim_{\ell\to\infty}\rk((F^{\circ\ell,*}_Y(V))_\min),
$$
and
$$
\brk_\max(V):=\lim_{\ell\to\infty}\rk((F^{\circ\ell,*}_Y(V))_\max).
$$
make sense if $V$ is a non zero locally free and coherent sheaf on $Y$.

{\it Suppose now that $t_0$ is only perfect} (not necessarily algebraically closed). If $V$ is a non zero coherent sheaf 
on $Y$, then we shall 
write $\mu(V):=\mu(V_{\bar t_0})$ and we shall 
say that $V$ is semistable if $V_{\bar t_0}$ is semistable. 
The \HN filtration of $V_{\bar t_0}$ is invariant under $\Gal(\bar t_0|t_0)$ by unicity and by 
a simple descent argument, we see that there is a unique filtration by coherent 
subsheaves
$$
V_0\subsetneq V_1\subsetneq V_2\subsetneq\dots\subsetneq V_{\hn(V)}
$$
such that 
$$
V_{0,\bar t_0}\subsetneq V_{1,\bar t_0}\subsetneq V_{2,\bar t_0}\subsetneq\dots\subsetneq V_{\hn(V),\bar t_0}
$$
is the \HN filtration of $V_{\bar t_0}$. We then define as before
$$\mu_\max(V):=\mu(V_1)$$
and $$\mu_\min(V):=\mu(V/V_{\hn(V)-1}).$$ Notice that we have 
$\mu_\max(V)=\mu_\max(V_{\bar t_0})$ and $\mu_\min(V)=\mu_\min(V_{\bar t_0})$. 

Notice that if 
$V$ and $W$ are non zero coherent locally free coherent sheaves on $Y$ and $\mu_\min(V)>\mu_\max(W)$ then 
we still have $\Hom_Y(V,W)=0$, since there is a natural 
inclusion $$\Hom_Y(V,W)\subseteq \Hom_{Y_{\bar t_0}}(V_{\bar t_0},W_{\bar t_0}).$$

If $t_0$ has positive characteristic, we shall 
say that $V$ is Frobenius semistable if $V_{\bar t_0}$ is Frobenius semistable. 
Since Frobenius morphisms commute with all morphisms, this is equivalent 
to requiring that $F^{r,*}_Y(V)$ is semistable for all $r\in\mN$ (with our extended 
definition of semistability).

We can now extend the range of the terminology introduced above:
$$
V_\max:=V_1,\,
V_\min:=V/V_{\hn(V)-1},
$$
$$
\bar\mu_\min(V):=\lim_{\ell\to\infty}\mu_\min(F^{\circ\ell,*}_Y(V))/p^\ell,\,
\bar\mu_\max(V):=\lim_{\ell\to\infty}\mu_\max(F^{\circ\ell,*}_Y(V))/p^\ell,
$$
$$
\brk_\min(V):=\lim_{\ell\to\infty}\rk((F^{\circ\ell,*}_Y(V))_\min),\,
\brk_\max(V):=\lim_{\ell\to\infty}\rk((F^{\circ\ell,*}_Y(V))_\max).
$$
Note that we have $\bar\mu_\min(V)=\bar\mu_\min(V_{\bar t_0})$, 
$\bar\mu_\max(V)=\bar\mu_\max(V_{\bar t_0})$, $\brk_\min(V)=\brk_\min(V_{\bar t_0})$, 
$\brk_\max(V)=\brk_\max(V_{\bar t_0})$
as expected.

If $V$ is a non zero coherent locally free coherent sheaf on $Y$ such that all the quotients of the 
\HN filtration of $V$ are Frobenius semistable, we shall say that 
$V$ has a Frobenius semistable \HN filtration. Note that by Theorem \ref{SSth} above, for 
any non zero coherent locally free coherent sheaf $V$ on $Y$, the sheaf $F^{\circ r,*}(V)$ 
has a Frobenius semistable \HN filtration for all but finitely many $r\in\mN$.

The following simple lemma will also prove very useful. It was suggested by J.-B. Bost.

\begin{lemma}
Let $V$ and $W$ be coherent locally free sheaves on $Y$. 
Suppose that $\mu(V)=\mu(W)$ and that $\rk(V)=\rk(W)$. Let $\phi:V\to W$ be a 
monomorphism of $\CO_Y$-modules. Then $\phi$ is an isomorphism.
\label{keyinjlem}
\end{lemma}
\beginProof We may suppose that $V$ and $W$ are of positive rank, otherwise the lemma is tautologically true. 
Let $M:=\det(W)\otimes\det(V)^\vee$. The assumptions imply that $\deg(M)=0$. Let 
$\det(\phi)\in H^0(Y,M)$ be the section induced by $\phi$. The zero scheme $Z(\det(\phi))$ of $\det(\phi)$ is a torsion sheaf since $\det(\phi)$ is non zero at the generic point of $Y$ and the length of $Z(\det(\phi))$ is equal to the degree of $M$ so $Z(\det(\phi))$ must be empty. In other words, 
$M$ is the trivial sheaf and $\det(\phi)$ is a constant non zero section of $M$. In particular, 
$\phi$ is an isomorphism.\endProof

\section{Finite flat group schemes over curves}

{The terminology of this section is independent of the introduction.}

\label{filtsec}

\subsection{Quotients by proper flat group schemes}

Let $Y$ be a noetherian scheme. 
Let $G$ be a commutative strongly quasiprojective flat group scheme over $Y$. 
See \cite[8.2, p. 211]{Bosch-Raynaud-Neron} for the definition of strong quasi-projectivity. 
Note that if $Y$ is regular then $G$ is strongly quasiprojective over $Y$ if it is 
quasiprojective over $Y$. 

Suppose that $H$ is a closed subgroup scheme of $G$, which is proper and flat over $Y$. The $Y$-scheme $G$ 
(resp. $H$) defines a functor $\underbar G$ (resp. $\underbar H$) from the category of 
$Y$-schemes to the category of abelian groups. Both functors are fppf sheaves 
by a classical result of Grothendieck. We may thus form the 
quotient $\underbar G/\underbar H$ of $\underbar G$ and $\underbar H$ in the 
category of fppf sheaves. 

The following proposition describes the quotient construction that we use in this text.

\begin{prop}
The fppf sheaf $\underbar G/\underbar H$ is representable by 
a group scheme $G/H$ over $Y$, which is also strongly quasiprojective. The natural morphism 
$q:G\to G/H$ is proper and faithfully flat and makes $G$ into an $H_{G/H}$-torsor 
over $G/H$. 
\label{propQ}
\end{prop}
\beginProof See \cite[Th. 8.12, p. 220]{Bosch-Raynaud-Neron}.\endProof

Note that if $G$ is semiabelian and $Y$ is normal then $G$ is quasiprojective 
over $Y$ (combine  \cite[VI.3.1]{MB-Pinceaux} with \cite[XI.1.4]{Raynaud-Faisceaux}). 
In particular if $Y$ is regular and $G$ is semiabelian then $G$ is strongly quasiprojective over $Y$. 

\subsection{The $\HN$-filtration on the Lie algebra of a finite flat group scheme of height one}

Let $S$ be a smooth, projective and  
geometrically connected curve over a perfect field 
$k$. Suppose that ${\rm char}(k)=p>0$.

\label{sssHone}

The following preliminary lemma will be very useful.

\begin{lemma}
Let $G$ be a finite flat commutative group scheme over $S$. 
Let $T\to S$ be a flat, radicial and finite 
morphism and let $\phi:H\hookrightarrow G_T$ be a closed 
subgroup scheme, which is finite, flat and 
multiplicative. Then there is a finite flat closed subgroup scheme 
$\phi_0:H_0\hookrightarrow G$, such  
that $\phi_{0,T}\simeq\phi$. 
\label{backlem}
\end{lemma}
\beginProof 
Taking Cartier duals, we get a morphism 
$$
\phi^\vee:G_T^\vee\to H^\vee.
$$
Notice that $H^\vee$ is \'etale over $T$, since 
$H$ is multiplicative. By radicial invariance of \'etale morphisms, 
there is a finite flat group scheme $J_0\to S$, such that 
$J_{0,T}\simeq H^\vee$. Notice also that 
the morphism $\phi^\vee$ is given by a section of the first 
projection
$$
G_T^\vee\times_T H^\vee\to G_T^\vee
$$
and since $H^\vee$ is \'etale over $T$, the image of this 
section is open and closed (see \cite[Cor. 3.12]{Milne-Etale}). Since the projection morphism 
$$
G_T^\vee\times_T H^\vee\to G^\vee\times_S J_0
$$
is also radicial, this open set comes from a unique open subset 
of $G\times_S J_0$ and this open subset defines an open and closed 
subscheme of $G^\vee\times_S J_0$, which is isomorphic to 
$G^\vee$ via the first projection. Hence the morphism 
$\phi^\vee$ comes from a unique morphism 
$G^\vee\to J_0$. Taking the Cartier dual of this morphism 
gives the morphism $\phi_0$.\endProof

Recall that a commutative finite flat group scheme $\psi:G\to S$ over $S$ is said to be {\it of height one} if 
$F_{G/S}=\epsilon_{G/S}\circ\psi$. 
Recall also that a (sheaf in) commutative $p$-Lie algebras (resp. $p$-coLie) algebras $V$ over $S$ 
 is a coherent locally free sheaf $V$ on $S$ together with a morphism of 
$\CO_S$-modules  
$F^*_S(V)\to V$ (resp. $V\to F_S^*(V)$). A morphism of commutative $p$-Lie (resp. $p$-coLie) algebras $V\to W$ is a 
morphism of $\CO_S$-modules from $V$ to $W$ satisfying an evident 
compatibility condition. There is a covariant functor $\Lie(\cdot)$ (resp. contravariant functor $\coLie(\cdot)$) from the category of commutative finite flat group schemes of height one over $S$ to the category of commutative $p$-Lie (resp. $p$-coLie) algebras , which sends 
a group scheme $G$ over $S$ to $\Lie(G):=\epsilon_{G/S}^*(\Omega_{G/S})^\vee$ 
(resp. $\coLie(G):=\epsilon_{G/S}^*(\Omega_{G/S})$, together 
with the morphism $$\Lie(V_{G\twp/S}):=(V_{G\twp/S}^*)^\vee:F_S^*(\Lie(G))=\Lie(G\twp)\to\Lie(G)$$

(resp. 
$$\coLie(V_{G\twp/S}):=V_{G\twp/S}^*:\coLie(G)\to F_S^*(\coLie(G\twp))=
\coLie(G\twp)$$)

Here $(V_{G\twp/S}^*)^\vee$ (resp. $V_{G\twp/S}^*$) is the dual of the pull-back morphism $V_{G\twp/S}^*$ (resp. is the pull-back morphism) on differentials induced by the Verschiebung morphism $V_{G\twp/S}$. 

The category of sheaves in commutative $p$-Lie algebras is tautologically 
antiequivalent to the category of sheaves in commutative $p$-coLie algebras. 

It can be shown that $\Lie$ is an equivalence of additive categories (see \cite[Expos\'e VIIA, rem. 7.5]{SGA3-1}). In particular, a 
sequence of finite flat group schemes of height one
$$
0\to G'\to G\to G''\to 0
$$
is exact if and only if the sequence 
$$
0\to\Lie(G')\to\Lie(G)\to\Lie(G'')\to 0
$$ is a sequence of commutative $p$-Lie algebras.  
Furthermore, we have
$$
\order(G)=p^{\rk(\Lie(G))}
$$
(see \cite[Proof of Th., p. 139, par. 14]{Mumford-Abelian}.)

\begin{lemma}
Let $\phi:V\to W$ be a morphism of commutative $p$-Lie algebras.
Then the image $\Im(\phi)$ (resp. the kernel $\ker(\phi)$) of $\phi$ as a morphism of $\CO_S$-modules is endowed with a unique structure of commutative $p$-Lie algebra, 
such that the morphism $\Im(\phi)\to W$ (resp. $\ker(\phi)\to V$) 
is a morphism of commutative $p$-Lie algebras.
\label{lemti}
\end{lemma}
\beginProof Left to the reader.\endProof

If $\phi:V\to W$ is an injective morphism of commutative $p$-Lie algebras, we 
shall say that $\Im(\phi)$ is a subsheaf in commutative $p$-Lie algebras. 
Beware that in this situation, the arrow $\phi$ might have no cokernel in 
the category of commutative $p$-Lie algebras. So in particular, 
$\Im(\phi)$ might not correspond to a subgroup scheme. On the other hand, if 
the quotient of $\CO_S$-modules $W/\Im(\phi)$ is locally free, then 
$W/\Im(\phi)$ can be endowed with an evident commutative $p$-Lie algebra structure, making it 
into a cokernel of $W$ by $\Im(\phi)$ in 
the category of commutative $p$-Lie algebras. In that case, $\Im(\phi)$ corresponds to a subgroup scheme.

We shall say that a finite flat commutative group scheme $G$ of height one (or its associated commutative $p$-Lie algebra) is {\it biinfinitesimal} if 
the associated morphism $F^*_S(\Lie(G))\to \Lie(G)$ is nilpotent. To say that $F^*_S(\Lie(G))\to \Lie(G)$ is nilpotent 
means that for some $n\geq 1$, the composition
$$
F^{\circ n,*}_S(\Lie(G))\to F^{\circ (n-1),*}_S(\Lie(G))\to\dots\to F^*_S(\Lie(G))\to \Lie(G)\to 0
$$
vanishes. We notice without proof that if 
$$
0\to G'\to G\to G''\to 0
$$
is an exact sequence of commutative finite flat group schemes, then 
$G'$ and $G''$ are biinfinitesimal if and only if $G$ is biinfinitesimal. Note also that 
a finite flat commutative group scheme $G$ of height one is multiplicative iff the associated morphism $F^*_S(\Lie(G))\to \Lie(G)$ is an isomorphism. This implies that if $G_1$ and $G_2$ are finite flat 
group schemes of height one over $S$, where $G_1$ is biinfinitesimal and $G_2$ is multiplicative then there 
are no non-zero morphisms of group schemes from $G_1$ to $G_2$ and also no non-zero morphisms of group schemes from $G_2$ to $G_1$.

We inserted the following alternative proof of a special case of Lemma \ref{backlem} to show the mechanics of $p$-Lie algebras at work in a simple situation. 

{\it Second proof of Lemma \ref{backlem} when $G$ is of height one and $T$ is smooth.} 

We may assume that $T\simeq S$ and that $T\to S$ is a power $F^{\circ n}_S$ of 
$F_S$. By induction on $n$, we are reduced to prove the statement for $n=1$. 

We are given a commutative diagram with exact rows and columns

\begin{equation*}
\xymatrix{
 & 0\ar[d]&\\
0\ar[r] & F_T^*(\Lie(H))\ar[rr]^{\small F_T^*(\Lie(\phi))}\ar[d]^{\Lie(V_{H/T})} & &F_T^*(\Lie(G)_T)\ar[d]^{\Lie(V_{G_T/T})}\\
0\ar[r] & \Lie(H)\ar[rr]^{\Lie(\phi)}\ar[d] & & \Lie(G)_T\\
& 0 &}
\end{equation*}
With the above reductions in place, this gives 
a commutative diagram with exact rows and columns
\begin{equation*}
\xymatrix{
 & 0\ar[d]&\\
0\ar[r] & F_S^*(\Lie(H))\ar[rr]^{\small F_S^*(\Lie(\phi))}\ar[d]^{\Lie(V_{H/S})} && F^{\circ 2,*}_S(\Lie(G))
\ar[d]^{F_S^*(\Lie(V_{G/S})}\\
0\ar[r] & \Lie(H)\ar[rr]^{\Lie(\phi)}\ar[d] & &F^{*}_S(\Lie(G))\\
& 0 & &}
\end{equation*}
Now consider the commutative diagram
\begin{equation*}
\xymatrix{
 & F_S^*(\Lie(H))\ar[rr]^{\small F_S^*(\Lie(\phi))}\ar[d]^{\Lie(V_{H/S})}
\ar[drr] && F^{\circ 2,*}_S(\Lie(G))
\ar[d]^{F_S^*(\Lie(V_{G/S})}\\
 & \Lie(H)\ar[drr]\ar[rr]^{\Lie(\phi)} & &F^{*}_S(\Lie(G))\ar[d]^{\Lie(G)}\\
&  & &\Lie(G)}
\end{equation*}
where the diagonal arrows are defined so that the diagram becomes commutative. 
The labelling of the arrows shows that the upper triangle is the base change by $F_S$ of the lower triangle. Hence the image of $\Lie(\phi)$ is the base change by 
$F_S$ of the image of $\Lie(H)$ in $\Lie(G)$, since $\Lie(V_{H/S})$ is an isomorphism. So $H_0$ can be defined as the group scheme of height one 
associated with the image of $\Lie(H)$ in $\Lie(G)$.\endProof

%The next lemma is parallel to Lemma \ref{lemRHO}. It is a refinement of Lemma \ref{lemRHO} if the underlying group scheme is of height one. 

\begin{lemma}  Let $V$ be a sheaf in commutative $p$-Lie algebras $V$ over $S$. 
Suppose that the \HN filtration 
$$
0=V_0\subsetneq V_1\subsetneq V_2\subsetneq\dots\subsetneq V_{\hn(V)}=V
$$
of $V$ is Frobenius semistable. Then for any $V_i$ such that $\mu_\min(V_i)\geq 0$, 
$V_i$ is a subsheaf in commutative $p$-Lie algebras $V$ over $S$. If 
$\mu_\min(V_i)>0$ then $V_i$ is biinfinitesimal.
\label{SBlem}
\end{lemma}
\beginProof  For the first statement, consider the morphism $\phi:F^*_S(V_i)\to V$ given by the composition of 
the inclusion $F^*_S(V_i)\to F^*_S(V)$ with the morphism $F^*_S(V)\to V$ given by the commutative $p$-Lie algebra structure. We have to check that the image of 
$\phi$ lies in $V_i$. The composition of 
$\phi$ with the quotient morphism $V\to V/V_i$ gives a morphism 
$F^*_S(V_i)\to V/V_i$ and it is equivalent to check that this morphism 
vanishes. Now compute
$$\mu_\min(F^*_S(V_i))=p\cdot\mu(V_i/V_{i-1})$$
 and 
 $$
 \mu_\max(V/V_i)=\mu(V_{i+1}/V_i)<\mu(V_i/V_{i-1})
 $$
 and thus $\mu_\min(F^*_S(V_i))>\mu_\max(V/V_i)$. We conclude that
 $\Hom_S(F^*_S(V_i),V/V_i)=0$ (see the discussion after Theorem \ref{SSth}) which concludes the proof of the first statement.
 To prove the second statement, it is sufficient by the remarks preceding the lemma to show that $V_i/V_{i-1}$ is biinfinitesimal for all indices $i$ such that $\mu(V_i/V_{i-1})>0$. 
 By the above computation, we have 
 $$\mu_\min(F^*_S(V_i/V_{i-1}))=\mu(F^*_S(V_i/V_{i-1}))=p\cdot \mu(V_i/V_{i-1})$$ 
 and thus $\mu_\min(F^*_S(V_i/V_{i-1}))>\mu(V_i/V_{i-1})$. Again, this implies that 
 $\Hom_S(F^*_S(V_i/V_{i-1}),V_i/V_{i-1})=0$, showing that $V_i/V_{i-1}$ is biinfinitesimal. 
\endProof 

\begin{rem}\rm As explained in the introduction, a characteristic $0$ analog of Lemma \ref{SBlem} can be found in \cite[Lemma 2.9]{Bost-Dwork}. See also \cite[Lemma 9.1.3.1]{SB-generic}, where 
a variant of a special case of Lemma \ref{SBlem} is proven under the assumption 
that $p$ is sufficiently large.\end{rem}

%\begin{lemma}
%Let $H$ be a finite flat group scheme of height one over $S$ and 
%suppose given an exact sequence 
%$$
%0\to H_\ord \to H\to H_\binf\to 0
%$$
%of finite flat group schemes such that $H_\ord$ is multiplicative and 
%$H_\binf$ is biinfinitesimal. Then the sequence is uniquely determined by 
%$H$ up to isomorphism of exact sequences.
%\end{lemma}

\begin{lemma}
Let $G$ be a commutative finite flat group scheme of height one over $S$ and 
suppose given an exact sequence 
$$
0\to G_\binf \to G\to G_\mu\to 0
$$
of finite flat group schemes such that $G_\mu$ is multiplicative and 
$G_\binf$ is biinfinitesimal. Then the sequence splits and this splitting is unique.
\label{lemsplit}
\end{lemma}
\beginProof
Consider the commutative diagram with exact rows and columns

\hskip3cm
\xymatrix{
0\ar[r] & \ker(\Lie(V^{(n)}_{G_\binf\twpn/S}))\ar[r]\ar[d]^{\simeq} & \ker(\Lie(V^{(n)}_{G\twpn/S}))\ar[r]\ar[d] & 0\ar[d]\\
0\ar[r] & F_S^{\circ n,*}(\Lie(G_\binf))\ar[r]\ar[d]^{=0} & F_S^{\circ n,*}(\Lie(G))\ar[r]\ar[d]& F_S^{\circ n,*}(\Lie(G_\mu))\ar[r]\ar[d]^{\simeq} & 0\\
0\ar[r] & \Lie(G_\binf)\ar[r] & \Lie(G)\ar[r]& \Lie(G_\mu)\ar[r] & 0
}

where $n\geq 0$ is chosen so that $V^{(n),*}_{G_\binf\twpn/S}=0$. Then the image of 
the arrow $$F_S^{\circ n,*}(\Lie(G))\to \Lie(G)$$ splits the bottom sequence. For the unicity of the splitting, 
note that for any two splittings $\sigma_1,\sigma_2$ of the bottom sequence the morphism 
$\sigma_1-\sigma_2:\Lie(G_\mu)\to \Lie(G)$ of vector bundles factors through the image of $\Lie(G_\binf)$. 
It thus defines a morphism of vector bundles $\Lie(G_\mu)\to \Lie(G_\binf)$, which is by construction a morphism 
of $p$-Lie algebras. Such a morphism must vanish (see the discussion after Lemma \ref{lemti}). 
Thus $\sigma_1=\sigma_2$.\endProof

%\begin{lemma}
%Suppose that $S$ is a Dedekind scheme. 
%Let $H$ be a finite flat group scheme of height one over $S$ and 
%let 
%$$
%0\to h' \to H_\eta\to h''\to 0
%$$
%be an exact sequence of finite flat group schemes over $\eta$. Then 
%this sequence extends uniquely to an exact sequence
%$$
%0\to H'\to H\to H''\to 0
%$$
%over $S$.
%\end{lemma}

\begin{lemma}
Let $G$ be a commutative finite flat group scheme of height one over $S$. 
Suppose that $\Lie(G)$ is Frobenius semistable of slope $0$. Let  
$n\geq 0$ be such that  $\rk(\ker(V^{(n),*}_{G\twpn/S}))$ is maximal. Then there is a canonical decomposition
$$
G\twpn\simeq H_\binf\times_S H_\mu
$$
where $H_\binf$ (resp. $H_\mu$) is a biinfinitesimal (resp. multiplicative) finite flat group scheme 
over $S$.
\label{lemtwistsplit}
\end{lemma}
\beginProof
Consider the commutative diagram

\hskip3cm
\xymatrix{
0\ar[r] & F_S^{\circ n,*}(\ker(\Lie(V^{(n)}_{G\twpn/S})))\ar[r]\ar[d]^{\sim} & F_S^{\circ n,*}(\ker(\Lie(V^{(n)}_{G\twpn/S})))\ar[r]\ar[d]& 0\ar[d]& \\
0\ar[r] & F_S^{\circ n,*}(\ker(\Lie(V^{(n)}_{G\twpn/S})))\ar[r]\ar[d]^{=0} & F_S^{\circ(2n),*}(\Lie(G))\ar[r]\ar[d]& F_S^{\circ n,*}(W)\ar[d]\ar[r]& 0\\
0\ar[r] & \ker(\Lie(V^{(n)}_{G\twpn/S}))\ar[r]& F_S^{\circ n,*}(\Lie(G))\ar[r]& W\ar[r]& 0
}

where $n\geq 0$ is such that $\rk(\ker(\Lie(V^{(n)}_{G\twpn/S})))$ is maximal 
and $W$ is the image of $\Lie(V^{(n)}_{G\twpn/S})$. The two bottom rows and the two leftmost columns in this diagram 
are exact by construction. Furthermore the map $F_S^{(n),*}W\to W$ is a monomorphism for otherwise $\rk(\ker(\Lie(V^{(n)}_{G\twpn/S})))$ is not maximal. The diagram thus has exact rows and columns. Since the second row gives 
a surjection 
$$
F_S^{\circ(2n),*}(\Lie(G))\to F_S^{\circ n,*}(W)
$$
we have $\mu_\min(F_S^{\circ n,*}(W))\geq 0$. Also, since the second column gives an injection 
$$
F_S^{\circ n,*}(W)\hookrightarrow F_S^{(n),*}(\Lie(G))
$$
we have $\mu_\max(F_S^{\circ n,*}(W))\leq 0$. Thus $F_S^{\circ n,*}(W)$ is  
of slope $0$. Thus $W$ is also of slope $0$. Hence by Lemma \ref{keyinjlem}, 
the monomorphism $$F_S^{\circ n,*}(W)\to W$$ is an isomorphism. Now we see that 
the image of 
the morphism $F_S^{\circ(2n),*}(\Lie(G))\to F_S^{\circ n,*}(\Lie(G))$ splits the bottom sequence. 
\endProof

\begin{lemma} Let $G$ be a finite flat commutative group scheme of height one over $S$. 
There exists a (necessarily unique) multiplicative subgroup scheme $G_\mu\hookrightarrow G$, 
such that if $H$ is a multiplicative subgroup scheme of height one over $S$ and 
$f:H\to G$ is a morphism of group schemes, then $f$ factors through $G_\mu$. 
Furthermore, for any $n\geq 0$, we have $(G_\mu)\twpn=(G\twpn)_\mu$. 
If $G$ is multiplicative over a dense open subset of $S$ and 
$\Lie(G)$ has Frobenius semistable \HN filtration then $\Lie(G)=\Lie(G)_{\leq 0}$ and $G_\mu$ corresponds to 
the subgroup scheme associated with $\Lie(G)_{=0}$.
\label{lemcansub}
\end{lemma}
\beginProof
In view of Lemma \ref{backlem}, we may replace $G$ by $G^\twpn$ for any $n\geq 0$ and in particular suppose that $\Lie(G)$ has a Frobenius semistable \HN filtration. Let $f:H\to G$ be a morphism of group schemes and 
consider the corresponding map
$$
\Lie(f):\Lie(H)\to\Lie(G).
$$
Since $H$ is multiplicative, $\Lie(H)$ is Frobenius semistable of slope $0$ (this is a consequence of Theorem \ref{SSth}). Thus the  image of 
$\Lie(f)$ lies in $\Lie(G)_{\geq 0}$. According to Lemma \ref{SBlem} there is an exact sequence of $p$-Lie algebras
$$
0\to \Lie(G)_{>0}\to \Lie(G)_{\geq 0}\stackrel{\pi}{\to} \Lie(G)_{=0}\to 0
$$
and we may assume according to Lemma \ref{lemtwistsplit} that there is a splitting 
$$
\Lie(G)_{=0}\simeq\Lie(G)_{=0,\binf}\oplus\Lie(G)_{=0,\mu}
$$
of $\Lie(G)_{=0}$ into multiplicative and biinfinitesimal part (we might have to twist $G$ some more for this). The inverse 
image of $\Lie(G)_{=0,\mu}$ by $\pi$ gives a $p$-Lie subalgebra 
$\pi^*(\Lie(G)_{=0,\mu})$ of $\Lie(G)_{\geq 0}$. This gives an exact sequence
$$
0\to \pi^*(\Lie(G)_{=0,\mu})\to \Lie(G)_{\geq 0}\to \Lie(G)_{=0,\binf}\to 0
$$
Since $\Lie(H)$ is multiplicative, the image of $\Lie(H)$ in $\Lie(G)_{=0,\binf}$ 
vanishes and thus the image of $\Lie(H)$ lies in $\pi^*(\Lie(G)_{=0,\mu})$. 
On the other hand by Lemma \ref{lemsplit} and Lemma \ref{SBlem}, we have again a canonical decomposition
$$
\pi^*(\Lie(G)_{=0,\mu})_\mu\oplus \pi^*(\Lie(G)_{=0,\mu})_\binf
$$
into multiplicative and binfinitesimal part and thus the image of 
$\Lie(f)$ lies in $\pi^*(\Lie(G)_{=0,\mu})_\mu$. Now $\pi^*(\Lie(G)_{=0,\mu})_\mu$ 
is a multiplicative $p$-Lie subalgebra of $\Lie(G)$ and it defines the required 
subgroup scheme.

If $G$ is multiplicative over an open subset of $S$ then 
we have an injection
$$
F_S^{\circ n,*}(\Lie(G))\hookrightarrow\Lie(G)
$$
(obtained by composition) for any $n\geq 0$ and thus if $\Lie(G)$ has Frobenius semistable \HN filtration then we must have $\Lie(G)=\Lie(G)_{\leq 0}$. Secondly the morphism \mbox{$
F_S^*(\Lie(G))\hookrightarrow\Lie(G)
$}
then induces an injection $$F_S^*(\Lie(G)_{=0})\hookrightarrow\Lie(G)_{=0}$$ and since both source 
and target in this map have the same rank and the same slope, we deduce 
from Lemma \ref{keyinjlem} that this map must be an isomorphism. Thus $\Lie(G)_{=0}$ is multiplicative and by the explicit construction above, it is associated 
with $G_\mu$. \endProof

%\begin{rem}\rm Note that if $G$ is multiplicative over an open subset of $S$ then Lemma \ref{EMlem} shows that $G_{\rho(G)-1}$ is a subgroup of multiplicative type and 
%thus if $G$ is also of height one, there is an embedding $G_{\rho(G)-1}\hookrightarrow 
%G_\mu$. It is likely that this is an isomorphism but we have not checked this.\end{rem}

\begin{rem}\rm Note that the "connected \'etale" decomposition of $G^\vee_K$ (see the beginning of \cite{Tate-Finite}) gives a canonical exact sequence of group schemes 
$$
0\to (G^\vee_K)_{\rm inf}\to G^\vee_K\to (G^\vee_K)_\et\to 0
$$
over $K$, where  $(G^\vee_K)_{\rm inf}$ is an infinitesimal group scheme and $(G^\vee_K)_\et$ is an \'etale group scheme over $K$. The group scheme $(G^\vee_K)_\et$ corresponds to 
a representation of $\Gal(K^\sep|K)$ into a finite $p$-group $E$ and one might be tempted 
to think that $G_\mu$ is the Cartier dual of the group scheme corresponding to  the largest unramified quotient of $E$, ie the largest quotient of $E$, such that the action of 
$\Gal(K^\sep|K)$ factors through the fundamental group $\pi_1(S)$. This not so, however. Indeed, consider 
a finite flat commutative group scheme $G$ of height one, which is such that $\bar\mu_\max(\Lie(G))<0$. 
Then $G_\mu=0$ and for any finite flat base change $S'\to S$, we also have $(G_{S'})_\mu=0$. 
On the other hand $(G^\vee_K)_\et$ will become constant (and hence entirely unramified) after a finite separable field extension $K'|K$ . 
\end{rem}

%\begin{cor}
%Let $G$ be a finite flat commutative group scheme of height one over $S$. Suppose 
%that $S$ a smooth and proper curve over $k$. Suppose that $\Lie(G)$ has a Frobenius semistable HN filtration. Then for any $\alpha\geq 0$, 
%$\Lie(G)_{\geq \alpha}$ is a $p$-Lie subalgebra of $\Lie(G)$. 
%For any $n$ such that $\rk(\ker(V^{(n),*}_{\Lie(G)_{=0}}))$ is maximal, there is a canonical splitting
%$$
%\Lie(F_S^{n,*}(G))_{=0}\simeq \Lie(F_S^{n,*}(G))_{=0,\mu}\oplus 
%\Lie(F_S^{n,*}(G))_{=0,\binf}
%$$
%and the exact sequence
%$$
%0\to\Lie(F_S^{n,*}(G))_{>0}\to \Lie(F_S^{n,*}(G))_{\geq 0}\to \Lie(F_S^{n,*}(G))_{=0}\to 0
%$$
%splits canonically over $\Lie(F_S^{n,*}(G))_{=0,\mu}$. 
%\label{corimp}
%\end{cor}
%
%\begin{cor}
%Let $n$ be such that $\Lie(F_S^{n,*}(G))$ has Frobenius semistable HN filtration. 
%Then the canonical subgroup scheme 
% $\Lie(F_S^{n,*}(G))_{=0,\mu}$ of $\Lie(F_S^{n,*}(G))_{=0}$ descends to a unique 
% subgroup scheme of $\Lie(G).$
% \end{cor}

\subsection{Quotients of semiabelian schemes by finite flat multiplicative group schemes}

Let $S$ be a smooth, projective and  
geometrically connected curve over a perfect field 
$k$.

\begin{lemma}
Let $\CA\to S$ be a semiabelian scheme. Suppose that 
there is an open dense subset $U\subseteq S$, such that 
$\CA_U\to U$ is an abelian scheme. Suppose 
that $G\hookrightarrow\CA$ is a finite, flat, closed 
subgroup scheme. Then the quotient 
scheme $\CA/G$ is also a semiabelian scheme and $(\CA/G)_U\to U$ is 
an abelian scheme. 
\label{BrionLem}
\end{lemma}
\beginProof Since the quotient morphism $q:\CA\to\CA/G$ is faithfully flat, the group scheme 
$\CA/G$ also has geometrically regular fibres (and is flat). Hence $\CA/G$ is smooth over $S$. 
Over $U$, its fibres are proper since the quotient morphism is also proper and they are thus abelian varieties. In other other words, $(\CA/G)_U\to U$ is 
an abelian scheme.  Now let $s\in S$. Since $(\CA/G)_s$ is smooth, we know by the 
Barsotti-Chevalley theorem (see \cite[Th. 10.25, p. 157]{Milne-Alg}) that 
$(\CA/G)_s$ sits in the middle of an exact sequence
\begin{equation}
0\to E_1\to (\CA/G)_s\to A_1\to 0
\label{qes}
\end{equation}
where $A_1$ is an abelian variety over $s$ and $E_1$ is a connected affine algebraic group variety over $s$.  The subgroup variety $E_1$ is maximal among connected affine subgroup varieties of $(\CA/G)_s$ (see \cite[Th. 10.5, p. 153 and proof]{Milne-Alg} and 
\cite[Th. 10.24, p. 156]{Milne-Alg}). Finally it has the form $E_1=T_1\times_s U$, where $T_1$ is a torus and $U$ is a connected unipotent group variety (see \cite[chap. 10.(i), p. 161]{Milne-Alg}). When we write that the sequence \refeq{qes} is exact, we mean 
that the third morphism is faithfully flat and that its kernel is $E_1$. 

By assumption, the corresponding presentation for $\CA_s$ is 
$$
0\to T\to\CA_s\to A_0\to 0
$$
where $T$ is a torus and $A_0$ is an abelian variety, both over $s$. 
%The scheme theoretic image of $T$ in $A_1$ is the $0$-section, since there are no non-trivial homomorphism 
%from a torus into an abelian variety (see ). Hence $q_s|_{T}$ factors through $E_1$. 
%We are thus given a commutative diagram
%
%\xymatrix{
%0\ar[r] & T\ar[r]\ar[d] & \CA_s\ar[r]\ar[d] & A_0\ar[r]\ar[d] & 0\\
%0\ar[r] & E_1\ar[r] & (\CA/G)_s\ar[r] & A_1\ar[r] & 0
%}

Let $D$ be the identity component of the closed subgroup scheme $q^{-1}_s(U\times 0)$ of $\CA_s$ (see \cite[Prop. 1.14]{Milne-Alg} for details). 
Since $s$ is perfect the closed subscheme $D_\red$ of $D$ is a closed subgroup scheme of $D$ 
(see \cite[Cor. 1.25, p. 24]{Milne-Alg}). Moreover $D$ and hence $D_\red$ is affine, since $q_s$ is finite. Since $T$ is the maximal connected affine subgroup variety of $\CA_s$, 
we see that $D_\red$ must be contained in $T$. However, every closed subgroup scheme of 
a multiplicative group over $s$ is multiplicative (see \cite[8.1, Exp. IX]{SGA3-2}) and thus $D_\red$ is multiplicative. Thus $D_\red$ is contained in the kernel of the morphism $q^{-1}_s(U\times 0)\to U\times 0$ (because there are no non trivial 
morphisms between multiplicative and unipotent algebraic groups - see 
\cite[Cor. 15.18, p. 255]{Milne-Alg}). Now notice that $q^{-1}_s(U\times 0)(\bar s)/
D(\bar s)$ is a finite set (see \cite[Prop. 1.14, p. 21]{Milne-Alg}). 
On the other hand $q_s(D(\bar s))=\{0\}$ by the above so $U(\bar s)$ must be finite. 
Since $U$ is smooth, it must thus be trivial. 
This shows that $(\CA/G)_s$ is an extension of an abelian variety by a torus. Since $s\in S$ 
was arbitrary, we see that $\CA/G$ is a semiabelian scheme.\endProof

\begin{lemma}
Let $G\to S$ be a finite flat group scheme 
of multiplicative type. Then there is a finite \'etale 
morphism $T\to S$ such that $G_T$ is a diagonalisable group scheme.
\label{Slem}
\end{lemma}
\beginProof See \cite[Exp. IX, Intro.]{SGA3-2}. \endProof

\begin{lemma}
Let $\CA\to S$ be a smooth commutative group scheme. Suppose 
that $G\hookrightarrow\CA$ is a finite, flat, closed 
subgroup scheme, which is multiplicative. Then 
$$
\deg(\omega_{\CA})=\deg(\omega_{\CA/G})
$$
\label{IDeglem}
\end{lemma}
\beginProof
By Lemma \ref{Slem}, we may assume that $G$ is diagonalisable. 
In particular, we may assume that there is a finite group scheme $G_0\to\Spec(k)$ such that 
$G_{0,S}\simeq G$. Let $\CB:=\CA/G$. Let $f:\CA\to S$ and $g:\CB\to S$ be the 
structural morphisms and let $\pi:\CA\to \CB$ be the quotient morphism. 
The triangle of cotangent complexes associated with the morphisms 
$\pi$, $g$ and $f$ gives an exact sequence
\begin{equation}
0\to \CH_1(\CTC(\pi))\to \pi^*(\Omega_g)\to\Omega_f\to\Omega_\pi\to 0
\label{CTS}
\end{equation}
where $\CTC(\pi)$ is the cotangent complex of $\pi$ and 
$\CH_1(\CTC(\pi))$ is its first homology sheaf. Now $\pi$ makes 
$\CA$ into a torsor over $\CB$ and under $G_\CB$. Hence there 
is a faithfully flat morphism $T\to \CB$ (for instance, we may take 
$T=\CA$), such that $\CA_T\simeq (G_\CB)\times_\CB T$. In particular we have 
$$\Omega_{\pi_T}\simeq \Omega_{G_0/k,T}$$ and $$\CH_1(\CTC(\pi_T))\simeq 
\CH_1(\CTC(G_0/k))_T$$
because the homology sheaves of the cotangent complex of 
$G_0$ over $k$ are flat (since they are $k$-vector spaces).

On the other hand, since $T\to \CB$ is flat, we have
$$\Omega_{\pi_T}\simeq\Omega_{\pi,T}$$
and
$$
\CH_1(\CTC(\pi_T))\simeq \CH_1(\CTC(\pi))_T
$$
Finally, notice that $\Omega_{G_0/k,T}$ and $\CH_1(\CTC(G_0/k))_T$ are flat 
and thus by flat descent, the sheaves $\CH_1(\CTC(\pi))$ and 
$\Omega_\pi$ are flat (in other words: locally free). Hence the sequence
\begin{equation}
0\to \epsilon_{\CA/S}^*(\CH_1(\CTC(\pi)))\to \epsilon_{\CB/S}^*(\Omega_g)\to\epsilon_{\CA/S}^*(\Omega_f)\to\epsilon_{\CA/S}^*(\Omega_\pi)\to 0
\label{CTS1}
\end{equation}
is also exact. Furthermore, we then have 
$$
\epsilon_{\CA/S}^*(\CH_1(\CTC(\pi)))\simeq \CH_1(\CTC(G_0/k))_S
$$
and 
$$
\epsilon_{\CA/S}^*(\Omega_\pi)\simeq\Omega_{G_0/k,S}
$$
and thus the sheaves $\epsilon_{\CA/S}^*(\CH_1(\CTC(\pi)))$ and 
$\epsilon_{\CA/S}^*(\Omega_\pi)$ are trivial sheaves. In particular, 
we have $\deg(\epsilon_{\CA/S}^*(\CH_1(\CTC(\pi))))=\deg(\epsilon_{\CA/S}^*(\Omega_\pi))=0$ and by the additivity of $\deg(\cdot)$, we deduce from 
the existence of the sequence \refeq{CTS1} that 
$
\deg(\omega_{\CA})=\deg(\omega_{\CA/G}).
$
\endProof

\begin{rem}\rm The computation of the cotangent complex made in the proof of Lemma \ref{Slem} is in essence also contained in  
\cite[Prop. 1.1]{Ekedahl-Surfaces} (but the assumptions made there are not quite the right ones for us).\end{rem}

\section{Proofs of the claims made in subsection \ref{ssIGP}}

\label{finsec}

{\it We now use the terminology of the introduction.} So let $k$ be a \underline{finite} field of characteristic $p>0$ and let 
$S$ be a smooth, projective and  geometrically connected curve over $k$. Let $K:=\kappa(S)$ be its function field. Let $A$ be an abelian variety of dimension $g$ over $K$. Fix an algebraic closure $\bar K$ of $K$. Let $K^\perf\subseteq\bar K$ be the maximal purely inseparable extension of $K$ and let $K^\unr\subseteq K^\sep$ be the maximal separable 
extension of $K$, which is unramified above every place of $K$. Finally, we let 
$\CA$ be a smooth commutative group scheme over $S$ such that $\CA_K=A$. 

{\it Proof of Theorem \ref{mainthm1}.} Recall the statement: there exists a (necessarily unique) multiplicative subgroup scheme $G_\CA\hookrightarrow\ker\,F_{\CA/S}$, with the following property: if $H$ is a multiplicative, finite and flat group scheme of height one over $S$ and 
$f:H\to\ker\,F_{\CA/S}$ is a morphism of group schemes, then $f$ factors through $G_\CA$. If $A$ is ordinary and $\omega_\CA$ is not ample then 
the order of $G_\CA$ is $p^{\brk_\min(\omega_\CA)}$. If $\phi:\CA\to\CB$ is a morphism
of smooth commutative group schemes over $S$, then the restriction of $\phi$ to $G_\CA$ factors through 
$G_\CB$. Furthermore, we have  
$\deg(\omega_{\CA})=\deg(\omega_{\CA/G_\CA}).$

In spite of its lengthy statement, the proof Theorem \ref{mainthm1} readily follows from Lemma \ref{lemcansub} and Lemma \ref{IDeglem}. More precisely, 
we simply have to define $G_\CA:=(\ker\,F_{\CA/S})_\mu$ in the notation of 
Lemma \ref{lemcansub}. The equality $\deg(\omega_{\CA})=\deg(\omega_{\CA/G_\CA})$ now follows from Lemma \ref{IDeglem}.
 
 {\it Proof of Proposition \ref{mainth0}.} Recall the assumptions of Proposition \ref{mainth0}: $A$ is ordinary, $\CA$ is semiabelian and 
 $A(K^\perf)$ is not finitely generated. We have to prove that $G_\CA$ is of order $>1$ and that 
 $\CA/G_\CA$ is also semiabelian. 
 
We know that $\bar\mu_\min(\omega_{\CA/S})\geq 0$ by 
Lemma \ref{lemcansub} and since $A(K^\perf)$ is not finitely generated, we 
know by Theorem \ref{theorA1}  that $\bar\mu_\min(\omega_{\CA/S})=0$. 
Proposition \ref{mainth0} now follows from Theorem \ref{mainthm1} and Lemma 
\ref{BrionLem}.
 
 {\it Proof of Proposition \ref{mainth}.} Recall the assumptions of Proposition \ref{mainth}: $A$ is ordinary, $\CA$ is semiabelian over $S$ and $A(K^\perf)$ is not finitely generated. 
We have to prove that there a finite flat morphism $$\phi:\CA\to\CB$$ where $\CB$ is a semiabelian over $S$ and 
a finite flat morphism $$\lambda:\CB\to\CB$$
such that $\ker(\phi)$ are $\ker(\lambda)$ are multiplicative group schemes and such that the order of 
$\ker(\lambda)$ is $>1$.
 
 Consider now $\CA_1:=\CA/G_\CA$. By Lemma \ref{BrionLem}, 
 the group scheme $\CA_1$ is also semiabelian and of course $A_1:=\CA_{1,K}$ is 
 also an ordinary abelian variety. We also have that $A_1(K^\perf)$ is not 
 finitely generated, since the natural map $A(K^\perf)\to A_1(K^\perf)$
  has finite kernel. Finally, the quotient morphism is 
  $\CA\to\CA_1$ is finite, flat, with multiplicative kernel and 
  $G_{\CA}$ is non trivial by Proposition \ref{mainth0}. 
  
  Repeating the above procedure for $\CA_1$ in place of $\CA$ and 
  continuing this way, we obtain an infinite sequence of semiabelian 
  schemes over $S$
\begin{equation}
\CA{\to}\CA_1{\to}\CA_2\to\dots
\label{InfSeq}
\end{equation}
where all the connecting morphisms are finite, flat, of degree $>1$
 and with multiplicative kernel. Applying Lemma \ref{IDeglem}, we see that 
$$
\deg(\omega_\CA)=\deg(\omega_{\CA_1})=\deg(\omega_{\CA_2})=\dots
$$
Let now $K'$ be a finite separable extension of $K$ such that 
$A(K)[n]\simeq(\mZ/n\mZ)^{2\dim(A)}$ for some 
$n\geq 3$ such that $(p,n)=1$. Let $S'$ be the normalisation of $S$ in $K'$. After base-change, 
we obtain an infinite sequence of semiabelian schemes over $S'$
\begin{equation}
\CA_{S'}{\to}\CA_{1,S'}{\to}\CA_{2,S'}\to\dots
\label{InfSeqp}
\end{equation}
and applying a theorem of Zarhin (see \cite[Th. 3.1]{Rossler-Infinitely} for a statement, 
explanations and further references), we conclude that in the sequence 
\refeq{InfSeqp}, there are only finitely many isomorphism classes of 
semiabelian schemes over $S'$. On the other hand, applying a basic finiteness 
result in Galois cohomology proven by Borel and Serre (see \cite[par. 3, p. 69]{Parshin-Zarhin-Finiteness}), we can now conclude that 
in the sequence 
\refeq{InfSeq}, there are also only finitely many isomorphism classes of 
semiabelian schemes over $S$.

Hence there are integers $j>i\geq 0$ 
and an isomorphism $$I:\CA_i\simeq\CA_j$$ over $S$. Letting  
$\phi:\CA\to\CA_i$ be the constructed morphism and letting $\lambda$ be the constructed morphism $\CA_i\to\CA_j$ composed with $I^{-1}$, 
we can now conclude 
the proof of Proposition \ref{mainth}.

%{\it Proof of Corollary \ref{corUnr}}. Recall the assumptions of Corollary \ref{corUnr}: $A$ is ordinary, $\CA$ is semiabelian over $S$ and $\Tor_p(A(K^\unr))$ is finite. 
%We have to prove that $A(K^\perf)$ is finitely generated.
%
%Suppose that $A(K^\perf)$ is not finitely generated.
%To prove Corollary \ref{corUnr}, we only have to prove that 
%$\Tor_p(A(K^\unr))$ is infinite.  
%
%Taking the dual of the morphism $\lambda_K$, we obtain a $K$-morphism 
%$$
%\lambda^{\vee}_K:A_i^\vee\to A_i^\vee
%$$
%such that $\ker(\lambda^{\vee}_K)$ is a finite \'etale group scheme over $K$, which extends to a finite \'etale group scheme of $p$-power order $>1$ over $S$. Hence for any $r\geq 0$, we 
%have 
%$$
%\ker(\lambda^{\vee,\circ r}_K)(\bar K)\subseteq \ker(\lambda^{\vee,\circ r}_K)(K^\unr)
%\subseteq\Tor_p(A^\vee_i(K^\unr))
%$$
%and since $r$ can be taken to be arbitrarily large, we see that 
%$\Tor_p(A^\vee_i(K^\unr))$ is infinite. Since $A_i^\vee$ is $K$-isogenous 
%to $A$, we conclude that $\Tor_p(A(K^\unr))$ is also infinite. 

%\section{Proofs of Theorem \ref{TamTh} and Corollary \ref{TamCor}}
%
%\label{Tamsec}

\section{Proofs of the claims made in subsection \ref{ssIGKS}}

\label{pcIGKS}

We start with the proof of Theorem \ref{classmain}. We recall the statement:

{\it Suppose that $\Tr_{\bar K|\bar k}(A_{\bar K})=0$. 
Suppose that the action of $\Gal(K^\sep|K)$ on $\Tor_p(A(K^\unr))$ factors 
through $\Gal(K^\sep|K)^\ab$. Then $\Tor_p(A(K^\unr))$ is finite.}

For the proof, let $L|K$ be the maximal subextension of $K^\unr|K$, which is Galois with abelian 
Galois group. Since $S$ is geometrically integral, $K\otimes_k\bar k$ is a 
field and $L$ contains a subfield 
isomorphic to $K\otimes_k\bar k$ (note that $\bar k=k^\sep$ and that 
$\Gal(\bar k|k)\simeq\widehat{\mZ}$, which is an abelian group). 
Furthermore, geometric class field theory (see eg \cite[Cor. 1.3]{Szamuely-Corps}) tells us that 
$\Gal(L\,|\,K\otimes_k\bar k)$ is a finite group. In particular, the 
field $L$ is finitely generated (as a field) over $\bar k$, since 
 $K\otimes_k\bar k$ is finitely generated over $\bar k.$
Now suppose to obtain a contradiction that $\Tor_p(A(K^\unr))$ were infinite. 
By assumption, we have 
$$
\Tor_p(A(K^\unr))\subseteq \Tor_p(A(L))
$$
Thus $\Tor_p(A(L))$ is infinite as well. By the Lang-N\'eron theorem, 
this implies that $$\Tr_{L|\bar k}(A_L)\not=0,$$ contradicting the first assumption.

We now turn to the proof of Proposition \ref{classcor}. We recall the statement:

{\it Suppose that $\dim(A)\leqslant 2$ and that $\Tr_{\bar K|\bar k}(A_{\bar K})=0$. 
Then $\Tor_p(A(K^\unr))$ is finite.}

For the proof, notice that if  $\Tor_p(A(K^\unr))$ is infinite then we have
$$
\bigcap_{\ell\geq 0}p^\ell\cdot\Tor_p(A(K^\unr))\not=0
$$
This follows from the fact that for each $n\geq 0$, the set 
$$
\{x\in\Tor_p(A(K^\unr))\,|\,p^n\cdot x=0\}
$$
is finite (the details are left to the reader). Let $G\subseteq \bigcap_{\ell\geq 0}p^\ell\cdot(\Tor_p(A(K^\unr)))$ 
be the subgroup of elements annihilated by the multiplication by $p$ map. 

If $G=0$ then there the conclusion holds, because then $\bigcap_{\ell\geq 0}p^\ell\cdot(\Tor_p(A(K^\unr)))=0$ and thus $\Tor_p(A(K^\unr))$ is finite by the above remark.

Suppose now that 
$\#G=p$. Then $\bigcap_{\ell\geq 0}p^\ell\cdot \Tor_p(A(K^\unr))$ is infinite and it is a union of cyclic groups of $p$-power order (use the classification theorem for finite abelian groups). Thus the action of $\Gal(K^\sep|K)$ on $\bigcap_{\ell\geq 0}p^\ell\cdot(\Tor_p(A(K^\unr)))$ 
factors through $\Gal(K^\sep|K)^\ab$. But this contradicts 
Theorem \ref{classmain} and thus we must have $\#G>p$. If $\#G>p$ then by the assumption that $\dim(A)\leq 2$, we see 
that we must have $\#G=p^2$ and thus the inclusions
$$
\Tor_p(A(K^\unr))\subseteq \Tor_p(A(K^\sep))\subseteq\Tor_p(A(\bar K))
$$
are both equalities. In particular, $A$ is an ordinary abelian surface. Let now $s\in S$ be a closed point such that 
$\CA_s$ is an ordinary abelian variety over $s$. Let $W:=\Spec(\widehat{\CO_{S,s}^\sh})$
 be the spectrum of the completion of the strict henselisation of the local ring at $s$ and write 
 $\widehat{K^\sh_s}$ for the fraction field of $\widehat{\CO_{S,s}^\sh}.$
 The abelian scheme $\CA_W\to W$ gives rise to an element $e$ of 
 $$
 \Hom_{\mZ_p}(T_p(\CA_{\bar s}(\bar s))\otimes T_p(\CA_{\bar s}^\vee(\bar s)),\widehat{\CO_{S,s}^\sh}^*).
 $$
 Here $T_p(\CA_{\bar s}(\bar s))$ and $T_p(\CA_{\bar s}^\vee(\bar s))$ are the $p$-adic Tate modules of 
 $\CA_{\bar s}$ and $\CA_{\bar s}^\vee$ respectively and $\widehat{\CO_{S,s}^\sh}^*$
 is the group of multiplicative units of $\widehat{\CO_{S,s}^\sh}$. The element $e$ is called the Serre-Tate pairing associated with 
 $\CA_W$. See \cite{Katz-Serre-Tate} for the construction of this pairing. We have $e=0$ if and only if $\CA_W\simeq\CA_{\bar s}\times_{\bar s}W$. Furthermore, the fact that 
 $$
 \Tor_p(\CA(W))=\Tor_p(A(\widehat{K^\sh_s}))=\Tor_p(A(K^\unr))=\Tor_p(A(\overline{\widehat{K^\sh_s})})
 $$
 in our situation shows that $e=0$. This follows directly from the definition of the Serre-Tate pairing in the ordinary case (see the definition of the morphism $"p^n"$ in 
 \cite[p.151]{Katz-Serre-Tate}). 
 Thus we have $\CA_W\simeq\CA_{\bar s}\times_{\bar s}W$ and in particular 
 $\Tr_{\bar K|\bar k}(A_{\bar K})\not=0$ by Proposition \ref{CRDprop} (c) below. This contradicts one of our assumptions. We conclude that $G=0$, so that the conclusion must hold. 

We shall now prove Theorem \ref{threp}. Recall the statement:

{\it Suppose that $\Tor_p(A(K^\sep))$ is infinite. Then there is an \'etale $K$-isogeny $$\phi:A\to B$$ where $B$ is an abelian variety over $K$ and there is an \'etale $K$-isogeny $$\lambda:B\to B$$
such that the order of 
$\ker(\lambda)$ is $>1$ and such that the orders of $\ker(\lambda)$ and $\ker(\phi)$ are powers of $p$.}

For the proof, note that in \cite[Th. 1.4]{Rossler-Infinitely}, this statement is proven under the supplementary assumption that there exist 
$n\in\mZ$, such that $(n,p)=1$ and $n>3$ and such that $A[n](\bar K)\simeq(\mZ/n\mZ)^{2\dim(A)}$. 
Using \cite[par. 3, "Finiteness Theorem for Forms", p. 69]{Parshin-Zarhin-Finiteness} in the proof, it can be seen that this assumption is not necessary. A completely parallel argument is described in 
the proof of Proposition \ref{mainth}. We leave the details to the reader. 

We now turn to the proof of Theorem \ref{thabet}. Recall the statement:

{\it Suppose that there exists an \'etale $K$-isogeny $\phi:A\to A$, such that $\deg(\phi)$ is strictly larger than $1$ and that 
$\deg(\phi)=p^r$ for some $r>0$. Suppose also that $A$ is a geometrically simple abelian variety and that $\CA$ is a semiabelian scheme.

Then $\CA$ is an abelian scheme and $\phi$ extends to an \'etale $S$-morphism $\CA\to \CA$ of group schemes.}

For the proof, notice first that by a result of Raynaud (see \cite[IX, Cor. 1.4, p. 130]{Raynaud-Faisceaux}), the morphism $\phi$ extends uniquely to an $S$-morphism 
$\bar\phi:\CA\to\CA$ of group schemes. Since $\bar\phi$ is \'etale 
over $K$, we have an exact sequence of 
coherent sheaves 
$$
0\to\bar\phi^*(\Omega_{\CA/S})\to \Omega_{\CA/S}
$$
on $\CA$. Let $\sigma\in H^0(\CA,\det(\bar\phi^*(\Omega_{\CA/S}))^\vee\otimes\det(\Omega_{\CA/S}))$ be the corresponding section. Since $\sigma_K\in H^0(A,\det(\phi^*(\Omega_{A/K}))^\vee\otimes\det(\Omega_{A/K}))$ has an empty zero-scheme, the zero scheme 
$Z(\sigma)$ is supported on a finite number of closed fibres of $\CA$. Hence there exists 
a finite number $P_1, \dots P_n$ of closed point of $S$, such that 
$Z(\sigma)=\coprod_{i=1}^n n_i\CA_{P_i}$ (as Weil divisors) for some 
$n_i\geqslant 0$. On the other hand, the Weil divisor $Z(\sigma)$ is rationally equivalent 
to $0$, since $\det(\phi^*(\Omega_{\CA/S}))^\vee\otimes\det(\Omega_{\CA/S})\simeq 
\det(\Omega_{\CA/S})^\vee\otimes\det(\Omega_{\CA/S})\simeq\CO_\CA$. 
Now notice that the morphism $p^*:\Pic(S)\to\Pic(\CA)$ of Picard groups is injective, 
because it is split by the map $\epsilon_{\CA/S}^*:\Pic(\CA)\to \Pic(S)$. Hence 
the Weil divisor $\coprod_{i=1}^n n_i P_i$ is rationally equivalent to $0$ on $S$, which 
implies that $n_i=0$ for all $i=1,\dots n$. In other words, we have 
$Z(\sigma)=\emptyset$ and thus the morphism $\bar\phi^*(\Omega_{\CA/S})\to \Omega_{\CA/S}$ is an isomorphism. By \cite[III, Prop. 10.4]{Hartshorne-Algebraic}, this implies that $\bar\phi$ is \'etale. 

Let now $s\in S$ be a closed point such that $\CA_s$ has a presentation 
$$
0\to G\stackrel{\iota}{\to}\CA_s\to A_0^0\to 0
$$
where $G$ is a torus over $s$ of dimension $d>0$ and $A_0^0$ is an abelian variety over $s$. 
The morphism $\bar\phi_s|_G:G\to \CA_s$ factors through $G$, since there is no non-constant 
$s$-morphism $G\to A_0^0$. Call $\gamma: G\to G$ the resulting morphism. The morphism $\gamma$ is 
\'etale. Indeed, we have a commutative diagram
\begin{center}
\hskip0.9cm
\xymatrix{
\gamma^*(\iota^*(\Omega_{\CA_s/s}))\ar[r]\ar[d]^{\sim} & \gamma^*(\Omega_{G/s})\ar[r] & \Omega_{G/s}
\ar[d]^{=}\\
\iota^*(\bar\phi_s^*(\Omega_{\CA_s/s}))\ar[r] & \iota^*(\Omega_{\CA_s/s})\ar[r] &\Omega_{G/s}.
}
\end{center}
and in the lower row of this diagram all the arrows are surjective. Thus 
the arrow $$\gamma^*(\Omega_{G/s})\to \Omega_{G/s}$$ must also be surjective and hence 
an isomorphism. Since $G$ is smooth over $\kappa(s)$, we conclude that 
$\gamma$ is smooth by \cite[III, Prop. 10.4]{Hartshorne-Algebraic}. 
In particular $\gamma$ is faithfully flat, because it is a morphism 
of group schemes and $G$ is connected (see eg \cite[SGA 3.1, Exp. IV-B, Cor. 1.3.2]{SGA3-1}). 
Now recall that there is a $K$-morphism $\psi:A\to A$ such that 
$\psi\circ\phi=[p^{\deg(\phi)}]_A$ (because finite commutative group schemes over $K$ are annihilated by 
their order; see \cite[Theorem (Deligne), p. 4]{Oort-Tate-Group}). The morphism $\psi$ extends uniquely 
to $\bar\psi:\CA\to \CA$ and thus by unicity, we have 
$\bar\psi\circ\bar\phi=[p^{\deg(\phi)}]_\CA.$ In particular, $\ker(\gamma)$ 
is a closed subscheme of $\ker([p^{\deg(\phi)}]_G)$. Since 
$\ker([p^{\deg(\phi)}]_G)$ is an infinitesimal group scheme and $\gamma$ is \'etale, we see that $\ker(\gamma)=0$ 
(since  $\ker(\gamma)$ is \'etale over $s$). Thus $\gamma$ is an isomorphism. 

Now choose a $\bar s$-isomorphism $G_{\bar s}\simeq \mG_m^{d}$ (here $\bar s$ if the 
spectrum of the algebraic closure of $\kappa(s)$). 
The  morphism $\gamma_{\bar s}$ is described by a matrix $M\in\GL_d(\mZ)$ 
(because the group scheme dual to $G_{\bar s}$ is the diagonalisable group scheme over 
$\bar s$ associated with $\mZ^d$). Hence there exists a monic polynomial $P(x)\in\mZ[x]$, such that 
$P(0)=\pm 1$ and such that $P(\gamma_{\bar s})=0$. 

Finally, choose a prime $l\not=p$. Let $\widehat{\CO}_s^\sh$ be the completion of the strict henselisation of the local ring of $S$ at $s$. Let $\widehat{K}^\sh_s$ be the fraction field of $\widehat{\CO}^\sh_s$ and let $j\in\mN$. The closed subgroup scheme 
$G_{\bar s}[l^j]$ of $G_{\bar s}$ extends uniquely to a finite and \'etale subgroup scheme $\wt{G}_{l^j}$ of $\CA_{\widehat{\CO}_s^\sh}$ over $\widehat{\CO}_s^\sh$. See \cite[Th. 3.6 and Th. 3.6 bis]{SGA3-2}. Furthermore the natural map $\wt{G}_{l^j}(\widehat{\CO}_s^\sh)\to G_{\bar s}[l^j](\bar s)$ 
is a bijection, since $\widehat{\CO}_s^\sh$ is strictly henselian and $\wt{G}_{l^j}$ is \'etale (see \cite[Prop. I.4.4]{Milne-Etale}). Hence $P(\phi)(\wt{G}_{l^j}(\widehat{\CO}_s^\sh))=0$. On the other hand, the image of the group $\bigcup_{j\in\mN}\wt{G}_{l^j}(\widehat{\CO}_s^\sh)$ in $A_{\widehat{K}^\sh_s}$ is dense, because $A$ is geometrically simple and the group $\bigcup_{j\in\mN}\widehat{G}_{l^j,\bar s}(\CO_s^\sh)$ is infinite. 
Hence $P(\phi)=0$ and since $P(0)=\pm 1$, we see that $\phi$ is an automorphism, which is a contradiction.\endProof

%The Tate module 
%$T_l(G_{\bar s})(\kappa(\bar s))$ is naturally a submodule of $T_l(\CA_{K_s^\sh})(K_s^\sh)$, where 
%$K_s^\sh$ is the fraction field of the strict henselisation of the local ring of $S$ at $s$ (it is 
%the "toric part" of the Tate module). 
%Hence the endomorphism $$P(T_l(\phi))\in\End_{\mZ_l}(T_l(\CA_{K_s^\sh})(K_s^\sh))$$ 
%vanishes. Since any infinite group of $l$-primary torsion points of $A(\bar K)$ is dense 
%in $A_{\bar K}$ (because $A$ is geometrically simple), this implies that $P(\phi)=0$. 
%Hence $\phi$ is an automorphism, which is a contradiction. 

\section{Proof of Theorem \ref{THA}}

\label{pTHA}

Recall the statement: 

{\it {\rm (a)} Suppose that $A$ is geometrically simple. If $A(K^\perf)$ is finitely generated and of \mbox{rank $>0$} then $\Tor_p(A(K^\sep))$ is a finite group.

{\rm (b)} Suppose that $A$ is an ordinary (not necessarily simple) abelian variety. If $\Tor_p(A(K^\sep))$ is a finite group then $A(K^\perf)$ is finitely generated.}

We shall need the following

\begin{lemma}
Let $B$ be an abelian variety over $K$ and let $\gamma:B\to B$ be a $K$-isogeny such that $\deg(\phi)>1$. Suppose that 
$B$ is geometrically simple. Let $H\subseteq A(\bar K)$ be a finitely generated subgroup. Then the 
set 
$$
\bigcap_{r\geq 0}\gamma^{\circ r}(H)
$$
is a finite group.
\label{wlem}
\end{lemma}
\beginProof (of Lemma \ref{wlem}) Let $G:=\bigcap_{r\geq 0}\gamma^{\circ r}(H).$ 
Let $F:=G/\Tor(G)$ be the quotient of $G$ by its torsion subgroup. We may suppose without restriction of generality  that $\rk(G)>0$ for otherwise the lemma is proven. Since $\gamma$ is a group homomorphism, 
we have $\gamma(\Tor(G))\subseteq \Tor(G)$ and thus $\gamma$ gives rise to a group homomorphism 
$F\to F$ that we also denote by $\gamma.$ By construction, we have $\gamma(F)=F$ and thus 
$\gamma:F\to F$ is a bijection, since $F$ is a finitely generated free $\mZ$-module. Let 
$$P(t):=t^n+a_{n-1}t^{n-1}+\dots +a_1 t+a_0\in\mZ[t]$$
be the characteristic polynomial of $\gamma:F\to 
F$. We have $P(\gamma)=0$ by the Cayley-Hamilton 
theorem and since $\gamma$ is an automorphism, 
we have have $$P(0)=a_0=\pm 1=\det(\gamma).$$
Hence
$$
(-a_0)^{-1}\cdot (\gamma^{\circ,n-1}+a_{n-1}\cdot \gamma^{\circ,n-2}+\dots a_1\cdot\Id_{F})
$$
is the inverse of $\gamma:F\to F$. Now let $\widetilde{\gamma}$ be the $K$-group scheme homomorphism
$$
\widetilde{\gamma}:=(-a_0)^{-1}\cdot (\gamma^{\circ,n-1}+a_{n-1}\cdot\gamma^{\circ,n-2}+\dots a_1\cdot\Id_B)
$$
from $B$ to $B$. Suppose first that the morphism of $K$-group schemes $\widetilde{\gamma}\circ\gamma-\Id_B$ is not the zero 
morphism. Then it is surjective, because $B$ is simple. Furthermore the group $G$ is dense in $B_{\bar K}$, since 
$B$ is geometrically simple. Thus the group $(\widetilde{\gamma}\circ\gamma-\Id_B)(G)$ is dense in 
$B_{\bar K}$. On the other hand, by construction $(\widetilde{\gamma}\circ\gamma-\Id_B)(G)\subseteq\Tor(G)$. Since 
$\Tor(G)$ is a finite group, it is not dense in $B_{\bar K}$ and thus we deduce that $\widetilde{\gamma}\circ\gamma-\Id_B$ must be the zero morphism. 
Hence $\gamma$ is invertible (with inverse $\widetilde{\gamma}$), which contradicts the assumption that 
$\deg(\gamma)>1$. We conclude that we cannot have $\rk(G)>0$ and thus $G=\Tor(G)$ is a finite group.\endProof

For the proof of Theorem \ref{THA} (a), suppose first that $\Tor_p(A(K^\sep))$ is not a finite group. Then 
by Theorem \ref{threp}, there exists an abelian variety $B$ over $K$, which is $K$-isogenous to $A$ and 
which carries an \'etale $K$-endomorphism $B\to B$, whose degree is $>1$ and is 
a power of $p$. The dual of $B$ hence carries an isogeny $\phi$, which is purely inseparable (because the dual of a finite \'etale 
group scheme over a field is an infinitesimal group scheme) and thus we have 
$$
B^\vee(K^\perf)=\bigcap_{r\geq 0}\phi^{\circ r}(B^\vee(K^\perf))
$$
By Lemma \ref{wlem}, $B^\vee(K^\perf)$ is thus either finite or not finitely generated and the same 
holds for $A$, since $A$ is isogenous to $B^\vee$. This proves (a).

We now turn to the proof of statement (b). Note that by Grothendieck's semiabelian reduction theorem, 
there is a finite and separable extension $K_1|K$ such that $A_{K_1}$ extends to a semiabelian scheme over the normalisation $S_1$ of $S$ in $K_1$. The scheme $S_1$ might not be 
geometrically connected over $k$ but there a finite extension $k_1$ of $k$, such that the connected components of $S_{1,k_1}$ are geometrically connected. We choose one of these connected components, say $S_2$. The extension of function fields 
corresponding to the morphism $S_2\to S$ is separable by construction so we may (and do)
assume that $\CA$ is semiabelian to begin with. Suppose that $A(K^\perf)$ is not finitely generated and that $A$ is ordinary. Then by 
Proposition \ref{mainth}, there is an abelian variety $B$ over $K$, which is $K$-isogenous to $A$ and 
which carries a $K$-isogeny $B\to B$, whose kernel is a multiplicative group scheme 
of order $>1$. The dual $\phi$ of this isogeny is an \'etale isogeny 
of $B^\vee$, which has degree $p^r$ for some $r>0$. Thus $\Tor_p(B^\vee(K^\sep))$ is an 
infinite group and the  same 
holds for $A$, since $A$ is isogenous to $B^\vee$. This proves (b).

\section{Proof of Theorem \ref{THB}}

\label{pTHB}

Recall the statement:

{\it Suppose that $\CA$ is a semiabelian scheme and that $A$ is a geometrically simple abelian 
variety over $K$.  
If $\Tor_p(A(K^\sep))$ is infinite, then 
\begin{itemize}
\item[{\rm (a)}] $\CA$ is an abelian scheme;
\item[{\rm (b)}]  there is $r_A\geq 0$ such that $p^{r_A}\cdot \Tor_p(A(K^\sep))\subseteq\Tor_p(A(K^\unr))$;
\end{itemize}
Furthermore, there is
\begin{itemize}
\item[{\rm (c)}]  an abelian scheme $\CB$ over $S$;
\item[{\rm (d)}]  a generically \'etale $S$-isogeny $\CA\to \CB$, whose degree is a power of $p$;
\item[{\rm (e)}] an \'etale $S$-isogeny $\CB\to \CB$ whose degree is $>1$ and is a power of $p$.
\end{itemize}
Finally
\begin{itemize}
\item[{\rm (f)}]  if $A$ is ordinary then the Kodaira-Spencer rank of $A$ is not maximal;
\item[{\rm (g)}]  if $\dim(A)\leqslant 2$ then $\Tr_{\bar K|\bar k}(A_{\bar K})\not=0$.
\item[{\rm (h)}] for all closed points $s\in S$, the $p$-rank of $\CA_s$ is $>0$.
\end{itemize}
}

Proof of (a): note that by Theorem \ref{threp}, the abelian variety $A$ is isogenous to an abelian variety $B$ over $K$, which is endowed with an \'etale endomorphism of degree a positive power of $p$. 
Since $A$ extends to a semiabelian scheme over $S$ so does $B$. This is 
a consequence of a theorem of Grothendieck (see \cite[5.]{Abbes-RSSC} for a nice presentation). Thus, by Theorem 
\ref{thabet} we see that $B$ extends to an abelian scheme $\CB$ over $S$. Using the criterion of N\'eron-Ogg-Shafarevich (see \cite{ST-NOS}), we see that $A$ also extends to an abelian scheme over $S$. 
By the uniqueness of semiabelian models (see \cite[IX, Cor. 1.4, p. 130]{Raynaud-Faisceaux}), this extension must be $\CA$ and thus $\CA$ is an abelian scheme. 

Proof of (b): Let $H:=\Gal(K^\sep|K^\unr)$. For $i\geq 0$, let $G_i:=A(K^\sep)[p^i]$. The group $G_i$ 
is the group of $K$-rational points of an \'etale finite group scheme 
$\underbar G_i$ over $K$, which is naturally a closed subgroup scheme of $A$. Let
$A_i:=A/\underbar G_i$ and for $i\leq j$ let $\phi_{i,j}:A_i\to A_j$ be the natural morphism. 
Let $\CA_i$ be the connected component of the zero section of the N\'eron model of $A_i$ over $S$. By (a) and the criterion of N\'eron-Ogg-Shafarevich (see \cite{ST-NOS}), this is an abelian scheme. Furthermore, 
by \cite[IX, Cor. 1.4, p. 130]{Raynaud-Faisceaux} the morphisms 
$\phi_{i,j}$ extend to morphisms $\bar\phi_{i,j}:\CA_i\to\CA_j$ and we have the classical exact sequence
$$
\bar\phi_{i,j}^*(\Omega_{\CA_j/S})\to\Omega_{\CA_i/S}\to\Omega_{\bar\phi_{i,j}}\to 0.
$$
Now the morphism $\bar\phi_{i,j}^*(\Omega_{\CA_j/S})\to\Omega_{\CA_i/S}$ is injective over the generic point of $\CA_i$, because $\phi_{i,j}=\bar\phi_{i,j,K}$ is smooth by construction. 
On the other hand both $\bar\phi_{i,j}^*(\Omega_{\CA_j/S})$ and $\Omega_{\CA_i/S}$ are locally free and thus it follows that $\bar\phi_{i,j}^*(\Omega_{\CA_j/S})\to\Omega_{\CA_i/S}$ is also injective. 
Hence we have an exact sequence
\begin{equation}
0\to\bar\phi_{i,j}^*(\Omega_{\CA_j/S})\to\Omega_{\CA_i/S}\to\Omega_{\bar\phi_{i,j}}\to 0.
\label{SA}
\end{equation}
Let $\pi_i:\CA_i\to S$ be the structural morphism. We have a functorial isomorphism 
$$\Omega_{\CA_i}\simeq\pi_i^*(\pi_{i,*}(\Omega_{\CA_i/S}))$$ and thus there is 
a coherent sheaf $T_{i,j}$ on $S$, which is a torsion sheaf, such that 
$\pi_i^*(T_{i,j})\simeq \Omega_{\bar\phi_{i,j}}$ and the sequence \refeq{SA} is the pull-back by $\pi_i^*$ of 
a sequence
$$
0\to \pi_{j,*}(\Omega_{\CA_j/S})\to \pi_{i,*}(\Omega_{\CA_i/S})\to T_{i,j}\to 0
$$
and in particular
$$
\deg_S(\pi_{j,*}(\Omega_{\CA_j/S}))+\deg_S(T_{i,j})=\deg_S(\pi_{i,*}(\Omega_{\CA_i/S})).
$$
Now recall that $\deg_S(\pi_{i,*}(\Omega_{\CA_i/S}))\geqslant 0$ for all $i\geq 0$ (see \cite[V, Prop. 2.2, p. 164]{FC-Degen}). 
Thus, for $i=0,1,\dots$ the sequence $\deg_S(\pi_{i,*}(\Omega_{\CA_i/S}))$ is a non-increasing sequence 
of natural numbers. Hence for large enough $i$, say $i_0$, it reaches its minimum. 
We conclude that $T_{i_0,j}=0$ for $j>i_0$, so that the morphism 
$\bar\phi_{i_0,j}$ is \'etale and finite. Now $\phi_{0,i_0}(G_{j}(K^\sep))$ lies by construction in the 
kernel of $\phi_{i_0,j}$. Thus 
$$
\phi_{0,i_0}(G_{j}(K^\sep))\subseteq A_{i_0}(K^\unr)
$$
when $j>i_0$. In other words, for any $x\in G_{j}(K^\sep)$ and any 
$\gamma\in H$, we have $$\gamma(x)-x\in G_{i_0}(K^\sep).$$ 
In particular, we have
$$
\gamma(p^{i_0}\cdot x)=p^{i_0}\cdot\gamma(x)=p^{i_0}\cdot x
$$
In particular, since $j>i_0$ was arbitrary, we see that 
$$
\gamma(p^{i_0}\cdot x)=p^{i_0}\cdot x
$$
for all $x\in \Tor_p(A(K^\sep))$ and all $\gamma\in H$. Setting $r_A=i_0$ concludes the proof of (b).

Proof of the existence statements (c), (d), (e): this is a consequence of (a) and Theorems \ref{threp} and \ref{thabet}.

Proof of (f): this is contained in a theorem of J.-F. Voloch; see \cite[Proposition on p. 1093]{Voloch-Dioph-p}.

Proof of (g): this is a consequence of (b) and Proposition \ref{classcor}.

Proof of (h): This follows from (a) and (e).

%\section{Proof of Theorem \ref{allC}}
%
%\label{pCEL}
%
%Note first that the fact that $A$ has property (P) implies that $\Tr_{\bar K|\bar k}(A_{\bar K})=0$. Hence 
%(a) immediately follows from Theorem \ref{THB} and from Th. B in \cite{Poonen-Voloch-BM}. Statement (c) follows immediately from 
%Theorem \ref{THB} and Theorem \ref{THA} (b).  Statement (c) follows immediately 
%from Theorem \ref{THB} and Theorem \ref{THA} (b), from Th. 2.2 in \cite{Ghioca-Moosa-Division} and from 
%the main result of \cite{Udi-ML}.

\section{Proof of Theorem \ref{ELth}}

\subsection{The trace of an abelian variety over a function field: basic facts}

\label{trss}

Let $E$ be an abelian over a field $F$. Let $F_0\subseteq F$ be a subfield. 

The $F|F_0$ trace $(\Tr_{F|F_0}(E),\lambda)$ (if it exists) of $E$ over $F_0$ is an abelian variety $\Tr_{F|F_0}(E)$ over 
$F_0$ together with a homomorphism $\lambda:\Tr_{F|F_0}(E)_F\to E$ 
of abelian varieties over $F$. They have the following universal property. For any abelian $E_0$ over $F_0$ and a homomorphism 
$\phi:E_{0,F}\to E$ of abelian varieties, there is a unique morphism 
$\wt{\phi}:E_{0,F}\to\Tr_{F|F_0}(E)_F$ such that $\phi=\lambda\circ\wt{\phi}$. 
This means that $\Tr_{F|F_0}(E)$ and $\lambda$ are uniquely determined if they exist. 

Here are some known facts about $\Tr_{F|F_0}(E)$. Before 
stating them, we record the fact for any finite morphism of abelian 
varieties $f:E'\to E$ over $F$, the natural morphism 
$E'/\ker(f)\to E$ is a closed immersion. Here $E'/\ker(f)$ is the quotient 
described in Proposition \ref{propQ}. To see this, consider that the morphism 
$E'/\ker(f)\to E$ is by definition a monomorphism of fppf sheaves over $F_0$ and 
hence a monomorphism of schemes. On the other hand, it is proper and of finite type and thus a closed immersion (see \cite[EGA IV.4, 18.12.6]{EGA}) for this). We shall call $\Im(f)$ the abelian variety 
$E'/\ker(f)$ viewed as an abelian subvariety of $E$.

The field extension $F|F_0$ is called primary (resp. regular) if the algebraic closure of $F_0$ in 
$F$ is purely inseparable over $F_0$ (if $F_0$ is algebraically closed in $F$ and $F$ is separable over $F_0$). Note that if $F$ is the function field of a smooth and geometrically 
integral variety over $F_0$ then $F|F_0$ is regular.

\begin{prop}$[${\rm see \cite[Th. 6.4 and Th. 6.12]{Conrad-Trace}}$]$

{\rm (a)} If $F|F_0$ is primary then the $F|F_0$ trace $(\Tr_{F|F_0}(E),\lambda)$ of $E$ over $F_0$ exists and the kernel of $\lambda$ is finite over $F$. 

{\rm (b)} If $F|F_0$ is regular then the kernel of the morphism $\lambda$ is connected 
and so is its Cartier dual.

{\rm (c)} If $F_1|F$ and $F|F_0$ are primary extensions then 
$(\Tr_{F|F_0}(E)_{F_1},\lambda_{F_1})$ is an $F_1|F_0$-trace of $E_{F_1}$.

{\rm (d)} We have $\Tr_{F|F_0}(A/\Im(\lambda))=0$. 

\label{CRDprop}
\end{prop}
 
We also recall the {\it Lang-N\'eron theorem} (see \cite[Th. 7.1]{Conrad-Trace} and 
\cite[chap. 6, Th. 2]{Lang-Fund}): if $F|F_0$ is a finitely generated regular 
extension then the quotient group $E(F)/\Tr_{F|F_0}(E)(F_0)$ is finitely generated. 
Here $\Tr_{F|F_0}(E)(F_0)$ is viewed as a subgroup of $E(F)$ via $\lambda$ and the natural 
base change map from $F_0$ to $F$.

\subsection{The proof}

\label{secELproof}

We now  use the notations of 
Conjecture \ref{ELconj}. 

Let $\lambda:\Tr_{L|l_0}(C)\to C$ be the canonical morphism. We write $C/\Im(\lambda)$ for the quotient of $C$ by $\Im(\lambda)$ in the sense of Proposition \ref{propQ}. 

We begin with the
\begin{prop} If $\IVD(C/\Im(\lambda),L)\subseteq\Tor^p((C/\Im(\lambda))(L))$ then 
$\IVD(C,L)\subseteq\Tor^p(C(L))$.
\label{traceless}
\end{prop}

For the proof of Proposition \ref{traceless}, we shall need the following 

\begin{lemma}
Let $N$ be a finite flat infinitesimal group scheme over a field $J$ of characteristic $p$. 
There is a finite field extension $J'|J$ such that for any $n\geq 0$ and 
any element $\alpha\in H^1(J,N\twpn)$, the image $\alpha_{J'}$ 
of $\alpha$ in $H^1(J',N\twpn_{J'})$ vanishes.
\label{vanlem}
\end{lemma}
Here $H^1(J,N\twpn)$ is the first cohomology group of $N\twpn$ viewed as a sheaf in the fppf topology. 
More concretely, it is the group of isomorphism classes of torsors of $N\twpn$ over $J$.
In the following proof, we shall write $J^{p^{-m}}\subseteq\bar J$ for the subfield of $\bar J$ consisting of 
elements of the form $x^{p^{-m}}$, where $x\in J$.  
\beginProof (of Lemma \ref{vanlem}) 
First suppose that $N$ 
has a filtration by finite closed subgroup schemes, whose 
quotients are isomorphic to either $\alpha_{p,J}$ or 
$\mu_{p,J}$.  Let $m\geq 0$ be the number of non vanishing quotients. 
We shall prove by induction on $m$ that the image of $\alpha$ in 
$H^1(J^{p^{-m}},N\twpn)$ vanishes for all $n\geq 0$ (under the supplementary assumption on $N$), for any field $J$ of characteristic $p$. If $m=0$ the statement holds tautologically, so we shall suppose that it holds for $1,\dots,m-1$. Let 
$$
0\to F_1\to N_{J_1}\to F_2\to 0
$$
be a presentation of $N$ where $F_2$ is isomorphic to either $\alpha_{p,J}$ 
or $\mu_{p,J}$ and $F_1$ has a filtration as above, whose number of non vanishing quotients is $\leq m-1$. This induces exact sequences
$$
0\to H^1(J^{p^{-1}},(F_{1,J^{p^{-1}}})\twpn)\to H^1(J^{p^{-1}},(N_{J^{p^{-1}}})\twpn)\to H^1(J^{p^{-1}},(F_{2,J^{p^{-1}}})\twpn)
$$
and
$$
0\to H^1(J^{p^{-m}},(F_{1,J^{p^{-m}}})\twpn)\to H^1(J^{p^{-m}},(N_{J^{p^{-m}}})\twpn)\to H^1(J^{p^{-m}},(F_{2,J^{p^{-m}}})\twpn)
$$
(observe that $H^0(J^{p^{-m}},(F_{2,J^{p^{-m}}})\twpn)=0$ since $F_2$ is infinitesimal). 
Since $F_2\twpn$ is of height one, the image of $\alpha$ in $H^1(J^{p^{-1}},(F_{2,J^{p^{-1}}})\twpn)$ 
vanishes by \cite[Lemma III.3.5.7]{Milne-Arith}. The element $\alpha$ is thus the image of 
an element $\beta\in H^1(J^{p^{-1}},(F_{1,J^{p^{-1}}})\twpn)$. By the inductive hypothesis, the image 
of $\beta$ in $H^1(J^{p^{-m}},(F_{1,J^{p^{-m}}})\twpn)$ vanishes and thus 
the image of $\alpha$ in $H^1(J^{p^{-m}},(N_{J^{p^{-m}}})\twpn)$ vanishes, proving the claim.  

Now according to \cite[par. 2.4, p. 28]{Groth-Barsotti} there is a finite extension $J_1$ of $J$ such that $N_{J_1}$ has a filtration by finite closed subgroup schemes, whose 
quotients are isomorphic to either $\alpha_{p,J_1}$ or 
$\mu_{p,J_1}$. This extension will by construction also work for all the 
group schemes $N\twpn$ and the number of non vanishing 
quotients of all the 
group schemes $N\twpn_{J_1}$ is constant, say it is $m$. Hence the extension $J':=J_1^{p^{-m}}$ has the required property.\endProof
 
\beginProof (of Proposition \ref{traceless}). Now suppose that $\IVD(C/\Im(\lambda),L)\subseteq\Tor^p((C/\Im(\lambda))(L))$. 
We want to show that $\IVD(C,L)\subseteq\Tor^p(C(L))$.  

Write 
$$\lambda\twpn:\Tr_{L|l_0}(C)\twpn\to C\twpn$$ for the base change of $\lambda$ by $F_L^{\circ n}$. We have an exact sequence 
$$
0\to\Im(\lambda)(L)\to C(L)\to (C/\Im(\lambda))(L)
$$
and we have $(C/\Im(\lambda))\twpn\simeq C\twpn/\Im(\lambda\twpn)$. 
Let now 
$$
x_0\in C(L), x_1\in C\twp(L),x_2\in C^{(p^2)}(L),\dots
$$
be a sequence of points such $V_{C\twp/L}(x_1)=x_0$, $V_{C^{(p^2)}/L}(x_2)=x_1$ etc. 
Then we know from the above supposition that the image of $x_n$ in $(C\twpn/\Im(\lambda\twpn))(L)$ is a prime to $p$ torsion point for all $n\geq 0$. In particular, the order $m$ of the image of $x_n$ in $(C\twpn/\Im(\lambda\twpn))(L)$ is 
independent of $n$, because the degree of the Verschiebung is always a power of $p$. Let $m$ be the order of $x_0$ (and hence of all the $x_n$). Then $m\cdot x_n\in \Im(\lambda\twpn)(L)$ 
for all $n$ and thus $m\cdot x_0$ is indefinitely Verschiebung divisible in $\Im(\lambda)(L)$ 
(because the Verschiebung morphism commutes with morphisms of 
commutative group schemes). It now suffices to prove that $m\cdot x_0$ is of finite and prime to $p$ order in $\Im(\lambda)(L)$. Hence, we may and do assume that the morphism 
$\lambda:\Tr_{L|l_0}(C)\to C$ is a surjection. 

Now $\lambda$ is also finite and 
purely inseparable by \cite[Th. 6.12]{Conrad-Trace} and it is thus a bijection. We are now given infinitely many $L$-morphisms 
$$
\dots(\lambda\twpn)^*(x_n)\to\dots\to(\lambda\twp)^*(x_1)\to\lambda^*(x_0)
$$
where $(\lambda\twpn)^*(x_n)$ is the base change by $\lambda\twpn$ of $x_n$ viewed as a closed subscheme of $C\twpn$. The $L$-scheme $(\lambda\twpn)^*(x_n)$ is a  torsor under the group scheme 
$(\ker\,\lambda)\twpn\simeq \ker\,\lambda\twpn$ and according to Lemma \ref{vanlem}, there is a finite 
extension $L'$, which splits all the $(\lambda\twpn)^*(x_n)$. We thus obtain an indefinitely Verschiebung divisible point $x'_0$ in $\Tr_{L|l_0}(C)(L')$, whose 
image in $C(L')$ is $x_0$. 
Now $\Tr_{L|l_0}(C)_{L'}$ is by definition the base change to $L'$ of an 
abelian variety over $l_0$; so  we are reduced to showing Theorem \ref{ELth} for abelian varieties $C$ 
that arise by base-change from $l_0$. Lemma \ref{isocase} below thus concludes the proof.\endProof

\begin{lemma} We have $\IVD(C,L)\subseteq\Tor^p(C(L))$ if $C\simeq C_{0}\times_{l_0} L$, where $C_0$ is an abelian variety over $l_0$. 
\label{isocase}
\end{lemma}
\beginProof (of Lemma \ref{isocase}) 
By \cite[Th. 6.2 and afterwards]{Esnault-Langer-On-a-positive}
 there is an $m\geq 1$ so that $m\cdot x_0\in C_0(l_0)$. Since $l_0$ is algebraically closed, this implies that $x_0\in C_0(l_0)$, concluding the proof.\endProof
 
\beginProof (of Theorem \ref{ELth}.)

We begin with a couple of reductions.

(1) {\it We may assume in the statement of Theorem \ref{ELth} that $L$ is the function field of a smooth and proper curve 
$B$ over $l_0$.} 

Using Proposition \ref{propQ} and Proposition \ref{CRDprop} (d), 
we see that when carrying out reduction (1), we may assume that $\Tr_{L|l_0}(C)=0$. Reduction (1) now follows from a standard spreading out argument together with Proposition 
\ref{redMW} in the Appendix. Here one could probably appeal instead 
to Hilbert's irreducibility theorem (as in \cite[chap. 9, cor. 6.3]{Lang-Fund}) but 
for lack of an adequate reference in the case of function fields, we prefer to use Proposition 
\ref{redMW}. 

(2) {\it We may assume in the statement of Theorem \ref{ELth}
that $\dim(\Tr_{\bar L|l_0}(C_{\bar L}))=\dim(\Tr_{L|l_0}(C))$.}

To see this, suppose for the space of this paragraph that we know that Theorem \ref{ELth} is true in general under restrictions (1) and (2). 
Let $L'|L$ be a finite extension such that $\dim(\Tr_{L'|l_0}(C_{L'}))$ is maximal 
among all finite extensions of $L$. In particular we then have $\dim(\Tr_{L'|l_0}(C_{L'}))=\dim(\Tr_{\bar L|l_0}(C_{\bar L}))$.  According to Proposition \ref{CRDprop} (c), we may assume that $L'|L$ is separable. 
Replacing $L'$ by the Galois closure of $L'$ over $L$, we may even suppose that $L'|L$ is Galois. 
Let $y_0\in C(L)$ be an indefinitely Verschiebung divisible element. Suppose $y_0\not=0$. 
Applying our assumptions  to 
$C_{L'}$ and to the normalisation $C'$ of $C$ in $L'$, we see that the image of $y_0$ in 
$C_{L'}(L')$ is indefinitely Verschiebung divisible. Thus for some integer $m_{y_0}$, which 
is prime to $p$, the element $m_{y_0}\cdot y_0$ is divisible in the group $C_{L'}(L')$. 
Now there is a natural group homomorphism $u:C_{L'}(L')\to C(L)$ (the trace) given 
by the formula
$$
u(z)=\sum_{\sigma\in\Gal(L'|L)}\sigma(z)
$$
Hence $m_{y_0}\cdot u(y_0)=m_{y_0}\cdot[L':L]\cdot y_0$ is divisible in the group $C(L)$ and hence 
$$m_{y_0}\cdot[L':L]\cdot y_0\in\Tr_{L|l_0}(C)(l_0).$$ Now if the order of the image of $y_0$ in 
$C(L)/\Tr_{L|l_0}(C(l_0))$ is prime to 
$p$ then we are done. Otherwise, we may (and do) replace $y_0$ by a multiple such that the image in $C(L)/\Tr_{L|l_0}(C(l_0))$ of $y_0$ is a non-zero  element of order $p$. In the rest of the argument, we shall derive a contradiction 
from the existence of this element. Let $i\geq 1$. Let $y_i\in C^{(p^i)}(L)$ be such that $V^{(i)}_{C^{(p^i)}/L}(y_1)=y_0$. 
The variety $$(C^{(p^i)})_{L'}=(C_{L'})^{(p^i)}\equiv C_{L'}^{(p^i)}$$ also has the property that 
$\dim(\Tr_{L'|l_0}(C^{(p^i)}_{L'}))=\dim(\Tr_{\bar L|l_0}(C^{(p^i)}_{\bar L}))$ since $C^{(p^i)}$ is isogenous to $C$ over $L$. Hence, repeating the above reasoning, there is an integer $m_{y_i}$, which is prime to 
$p$, such that $m_{y_i}\cdot[L':L]\cdot y_i\in\Tr_{L|l_0}(C^{(p^i)}(L)$. 
Now according to Proposition \ref{CRDprop} (c), the natural morphism $\Tr_{L|l_0}(C)^{(p^i)}_L\to C^{(p^i)}$ obtained
by base change under $F_C^{\circ i}$ from the morphism $\Tr_{L|l_0}(C)_L\to C$ makes 
$\Tr_{L|l_0}(C)^{(p^i)}$ into the trace of $C^{(p^i)}$. Thus the map
$V^{(i)}_{C^{(p^i)}/L}(\bar L)$ (resp. $F^{(i)}_{C/L}(\bar L)$) induces a surjective map 
$$C^{(p^i)}(\bar L)/\Tr_{L|l_0}(C)^{(p^i)}(l_0)\to
C(\bar L)/\Tr_{L|l_0}(C)(l_0)$$ (resp. a bijective map 
$$C(\bar L)/\Tr_{L|l_0}(C)(l_0)\to C^{(p^i)}(\bar L)/\Tr_{L|l_0}(C)^{(p^i)}(l_0)$$). 
Since $V^{(i)}_{C^{(p^i)}/L}(\bar L)\circ F^{(i)}_{C/L}(\bar L)=p^i$, we see that the order of 
$y_i$ in $$
C^{(p^i)}(L)/\Tr_{L|l_0}(C^{(p^i)})(l_0)\subseteq C^{(p^i)}(\bar L)/\Tr_{L|l_0}(C^{(p^i)})(l_0)$$ is $p^{i+1}.$
This is a contradiction if $i$ is chosen large enough so that $p^i$ is not a divisor of 
$[L':L]$. We conclude that the order of the image of $y_0$ in 
$C(L)/\Tr_{L|l_0}(C)(l_0)$ is prime to 
$p$ and this concludes reduction step (2).

We now assume that we are given an abelian variety $C$ over $L$ and that $C$ satisfy the assumptions of \ref{ELth} as well as (1) and (2). 

Let as before $\lambda:\Tr_{L|l_0}(C)_L\to C$ be the canonical morphism.
According to Proposition \ref{traceless}, it will be 
sufficient to prove that $\IVD(C/\Im(\lambda),L)\subseteq\Tor^p((C/\Im(\lambda))(L)).$ 
By Theorem \ref{CRDprop} (d), we have $\Tr_{L|l_0}(C/\Im(\lambda))=0$ and since we work under 
supplementary assumption (2), we even have $\Tr_{\bar L|l_0}(C/\Im(\lambda))=0.$ 
Thus we may replace $C$ by $C/\Im(\lambda)$ and assume from now on that 
$\Tr_{\bar L|l_0}(C)=0.$ Finally, since we have $\Tr_{\bar L|l_0}(C)=0$, we may 
replace without restriction of generality replace $L$ by a finite extension $L'$ and 
$C$ by its normalisation $C'$ in $L'$. We may thus assume that there is an integer $m\geq 3$, with $(m,p)=1$ and such that 
$C[m]\simeq(\mZ/m\mZ)^{2\dim(C)}$ and $C^\vee[m]\simeq(\mZ/m\mZ)^{2\dim(C^\vee)}$.

By a theorem of Raynaud (see \cite[Prop. 5.10]{Abbes-RSSC}), the connected component of the N\'eron model 
of $C$ will then be a semiabelian scheme. We call it $\CC$.

Now suppose as in the statement of Conjecture \ref{ELconj} that we are given points 
$x_\ell\in C^{(p^\ell)}(L)$ and suppose that for all $\ell\geqslant 1$, we have $V_{C^{(p^\ell)}/L}(x_{\ell})=x_{\ell-1}$. We want to show that $x_0\in\Tor^p(C(L))$.

By Lemma \ref{flem} and the discussion preceding it we have a canonical map  
\begin{equation}
\alpha:C\twp(L)\to \Hom_B(\omega_{\CC\twp},\Omega_{B/l_0}(E))
\label{canimpinj}
\end{equation}
such that $\alpha(x)=0$ iff $x\in F_{C/L}(C(L))$. Here $E=E(\CC)$ is the reduced divisor, which is the union of the closed point $b\in B$ such that $\CC_b$ is not proper over $\kappa(b)$. Note that 
we have $E(\CC)=E(\CC\twp)=E(\CC^{(p^2)})=\dots$.  
The  map $\alpha$ is naturally compatible with isogenies (we skip the verification) 
and so there is an infinite commutative diagram

\
\begin{equation}
\xymatrix{
C\twp(L)\ar[r] & \Hom_B(\omega_{\CC\twp},\Omega_{B/l_0}(E))\\
C^{(p^{2})}(L)\ar[r]\ar[u]^{V_{C^{(p^{2})}/L}} & \Hom_B(\omega_{\CC^{(p^{2})}},\Omega_{B/l_0}(E))
\ar[u]^{V^*_{\CC^{(p^{2})}/B}}\\
\vdots\ar[u]\ar[r] & \vdots\ar[u]
}
\label{icd}
\end{equation}

 Remember that we have 
$$
\omega_{\CC^{(p^{n})}}\simeq F^{\circ n,*}_B(\omega_{\CC}).
$$
Now choose $n_1\geq 1$ so that 

- $\omega_{\CC^{(p^{n_1})}}$ has a Frobenius semistable \HN filtration;

- $(\omega_{\CC^{(p^{n_1})}})_{=0}\simeq (\omega_{\CC^{(p^{n_1})}})_{=0,\binf}\oplus 
(\omega_{\CC^{(p^{n_1})}})_{=0,\mu}$ splits into a biinfinitesimal and a multiplicative commutative coLie-algebra (see Lemmata \ref{SBlem} and \ref{lemtwistsplit}).

Note that if some $n_1\geq 1$ has the two above properties, than any higher 
$n_1$ will as well (by definition for the first property and tautologically for the second one). 

Choose $n_2>n_1$  so that 

(I) the image of the map 
$$
V^{(n_2-n_1),*}_{\CC^{(p^{n_2})}/B}:\omega_{\CC^{(p^{n_1})}}\to \omega_{\CC^{(p^{n_2})}}
$$ 
lies in $(\omega_{\CC^{(p^{n_2})}})_{\geq 0}\simeq F_B^{\circ (n_2-n_1),*}((\omega_{\CC^{(p^{n_1})}})_{\geq 0});$

(II) the image of the map of coLie algebras
$$
V^{(n_2-n_1),*}_{\CC^{(p^{n_2})}/B}:(\omega_{\CC^{(p^{n_1})}})_{=0}\to 
F_B^{\circ (n_2-n_1),*}((\omega_{\CC^{(p^{n_1})}})_{=0})=(\omega_{\CC^{(p^{n_2})}})_{=0}
$$
is $F_B^{\circ (n_2-n_1),*}((\omega_{\CC^{(p^{n_1})}})_{=0,\mu})$. Note that this 
is possible because the biinfinitesimal part of $(\omega_{\CC^{(p^{n_1})}})_{=0}$ 
will be sent to $0$ by sufficiently many composed Verschiebung morphisms (by definition). 

Note that under (I) for any $n_3>n_2$  the image of the map 
$$
V^{(n_3-n_2),*}_{\CC^{(p^{n_3})}/B}:(\omega_{\CC^{(p^{n_2})}})_{\geq 0}\to \omega_{\CC^{(p^{n_3})}}
$$ 
and hence of the map 
$$
V^{(n_3-n_1),*}_{\CC^{(p^{n_3})}/B}:\omega_{\CC^{(p^{n_1})}}\to \omega_{\CC^{(p^{n_3})}}
$$
automatically lies in $(\omega_{\CC^{(p^{n_3})}})_{\geq 0}\simeq F_B^{\circ (n_3-n_1),*}((\omega_{\CC^{(p^{n_1})}})_{\geq 0})$.

Choose $n_3>n_2$ so that 

(III) the map 
$$
\omega_{\CC^{(p^{n_3})}}\to\Omega_{B/l_0}(E)
$$
 given by $x_{n_3}$ factors through its quotient 
$(F^{\circ n_3,*}_B(\omega_{\CC}))_{\leq 0}\simeq F^{\circ (n_3-n_1),*}_B((\omega_{\CC^{(p^{n_1})}})_{\leq 0})$;

(IV) the image of the map $$
V^{(n_3-n_2),*}_{\CC^{(p^{n_3})}/B}:F_B^{\circ (n_2-n_1),*}((\omega_{\CC^{(p^{n_1})}})_{=0})\to
F_B^{\circ (n_3-n_2),*}((\omega_{\CC^{(p^{n_2})}})_{=0})
$$
is 
$F_B^{\circ (n_3-n_2),*}((\omega_{\CC^{(p^{n_2})}})_{=0,\mu})\simeq 
F_B^{\circ (n_3-n_1),*}((\omega_{\CC^{(p^{n_1})}})_{=0,\mu})$.
 
Now we shall exploit the compatibility between the morphism
$$
\omega_{\CC^{(p^{n_1})}}\stackrel{c(x_{n_1})}{\to}\Omega_{B/k}(E)
$$
induced by $x_{n_1}$ and the morphism
$$
\omega_{\CC^{(p^{n_3})}}\stackrel{c(x_{n_3})}{\to}\Omega_{B/k}(E)
$$
induced by $x_{n_3}$. According to the diagram \refeq{icd}, this compatibility gives the 
equality 
$$
c(x_{n_3})\circ V^*_{\CC^{(p^{n_3-n_1})}/B}=c(x_{n_1}).
$$
In other words the composition of morphisms 
$$
\omega_{\CC^{(p^{n_1})}}\stackrel{V^*_{\CC^{(p^{n_3-n_1})}/B}}{\to} \omega_{\CC^{(p^{n_3})}}\stackrel{c(x_{n_3})}{\to}\Omega_{B/k}(E)
$$
is $c(x_{n_1})$. Furthermore, in view of (I) and (III) the map 
$c(x_{n_1})$ factors as follows: 
$$\hskip-1.5cm
\omega_{\CC^{(p^{n_1})}}\stackrel{V^{(n_3-n_1),*}_{\CC^{(p^{n_3})}/B}}{\to}F_B^{\circ (n_3-n_1),*}((\omega_{\CC^{(p^{n_1})}})_{\geq 0})\stackrel{}{\to} 
F_B^{\circ (n_3-n_1),*}((\omega_{\CC^{(p^{n_1})}})_{=0})\to
F_B^{\circ (n_3-n_1),*}((\omega_{\CC^{(p^{n_1})}})_{\leq 0})\to \Omega_{B/k}(E)
$$
and by (I) the map 
$$
\omega_{\CC^{(p^{n_1})}}\stackrel{V^{(n_3-n_1),*}_{\CC^{(p^{n_3})}/B}}{\to} 
F_B^{\circ (n_3-n_1),*}((\omega_{\CC^{(p^{n_1})}})_{=0})
$$
factors as follows 
$$
\hskip-1cm
\omega_{\CC^{(p^{n_1})}}\stackrel{V^{(n_2-n_1),*}_{\CC^{(p^{n_1})}/B}}{\to} 
F_B^{\circ (n_2-n_1),*}((\omega_{\CC^{(p^{n_1})}})_{\geq 0})\to 
F_B^{\circ (n_2-n_1),*}((\omega_{\CC^{(p^{n_1})}})_{=0})\stackrel{V^{(n_3-n_2),*}_{\CC^{(p^{n_1})}/B}}{\to}
F_B^{\circ (n_3-n_1),*}((\omega_{\CC^{(p^{n_1})}})_{=0})
$$
and thus by (IV) and (II) the image of this last map is precisely
$F_B^{\circ (n_3-n_1),*}((\omega_{\CC^{(p^{n_1})}})_{=0,\mu})$. 

We have thus constructed a multiplicative quotient of the $p$-coLie algebra  
$\omega_{\CC^{(p^{n_1})}}$. On the other hand the $p$-coLie algebra $\omega_{\CC^{(p^{n_1})}}$ is the $p$-coLie algebra of the finite flat group scheme  
$\ker\,F_{\CC^{(p^{n_1})}/B}$. By the equivalence of 
categories recalled in subsubsection \ref{sssHone}, this quotient corresponds to 
a multiplicative subgroup scheme of $\ker\,F_{\CC^{(p^{n_1})}/B}$. By Lemma \ref{lemcansub}, this subgroup scheme embeds in the canonical largest multiplicative subgroup scheme $(\ker\,F_{\CC^{(p^{n_1})}/B})_\mu$ of $\ker\,F_{\CC^{(p^{n_1})}/B}$ (in fact, it coincides with it, but we shall not need this). 
Finally note that $$(\ker\,F_{\CC^{(p^{n_1})}/B})_\mu\simeq((\ker\,F_{\CC/B})_\mu)^{(p^{n_1})},$$ by the last part of Lemma \ref{lemcansub}.

Let $G:=(\ker\,F_{\CC/B})_\mu$. Note that $G=G_\CC$ in the notation of Theorem \ref{mainthm1}. Now consider the quotient $\CC_1:=\CC/G$ (which is a semiabelian scheme by \ref{BrionLem}) and let $\psi_1:\CC\to\CC_1$ be the quotient morphism. The 
point $x_{n_1}$ and its image $y_{n_1}$ in $\CC_1(L)$ give a commutative diagram

\hskip4cm
\xymatrix{
0 & F_B^{\circ (n_3-n_1),*}((\omega_{\CC^{(p^{n_1})}})_{=0,\mu})\ar@/^6pc/[dd]\\
\omega_{G^{(p^{n_1})}}\ar[u]\ar[ur] & F_B^{\circ (n_3-n_1),*}(\omega_{\CC^{(p^{n_1})}})\ar[d]\ar[u]\ar[l] \\
\omega_{\CC^{(p^{n_1})}}\ar[r]^{c(x_{n_1})}\ar[u]\ar[ur]& \Omega_{B/k}(E)\ar[d]^{\simeq}\\
\omega_{\CC_1^{(p^{n_1})}}\ar[r]^{c(y_{n_1})}\ar[u]^{\psi_1^*}& \Omega_{B/k}(E)
}

where the left column is an exact sequence and $c(y_{n_1})$ is the morphism induced by $y_{n_1}$. 

{\it Thus $c(y_{n_1})$ vanishes.} In particular, $y_{n_1}$ lies in the image 
of $F_{\CC_1^{(p^{n_1-1})}/B}(\CC_1^{(p^{n_1-1})}(L))$. Using the fact that $$[p]_{\CC_1^{(p^{n_1-1})}}=
V_{\CC_1^{(p^{n_1})}/B}\circ F_{\CC_1^{(p^{n_1-1})}/B},$$ we conclude that $y_{n_1-1}$ has 
a $p$-th root in $\CC_1^{(p^{n_1-1})}(L)$. Hence $y_0$ also has a $p$-th root in 
$\CC_1(L)$. Now since $G$ is independent of $x_0$, 
we conclude that the image of any indefinitely Verschiebung divisible 
point of $C(L)$ in $\CC_1(L)$ has a $p$-th root. Since $G$ is compatible with twists, we also see that for any $n\geq 0$ the image of any indefinitely Verschiebung divisible 
point of $C\twpn(L)$ in $\CC_1\twpn(L)$ has a $p$-th root. From this, by an elementary combinatorial consideration, we see that the image of any indefinitely Verschiebung divisible 
point of $C(L)$ in $\CC_1(L)$ has a $p$-th root, which is indefinitely Verschiebung divisible.

 By the discussion above, {\bf the image of $\IVD(C)$ in $\CC_1(L)$ lies in $p\cdot\IVD(\CC_{1,L})$.} This is the crucial fact that the rest of the proof will exploit.

Let $\CC_1:=\CC/G_\CC,\,\CC_2/G_{\CC_1},\dots$ be the sequence of smooth commutative group schemes 
obtained by successively quotienting by the canonical subgroup schemes described in 
Theorem \ref{mainthm1}.  Note that 
all the $\CC_i$ are semiabelian by Lemma \ref{BrionLem}. We shall denote by $\psi_i$ the morphism 
$\CC\to\CC_i$ obtained by composition. Write $C_i:=\CC_{i,L}$ for convenience. 

Let $m_{00}$ be an integer such that $m_{00}\cdot x_0=:v_0$ extends to an element $\wt{v}_0$ of $\CC(B)$. 

%Our strategy to finish the proof is now the following. Since $x_0\in\IVD(C)$, 
%we know that $\psi_{i,L}(x_0)$ is divisible by $p^i$ in $C_i(L)$ by Proposition  .  
%Up to replacing 
%$x_0$ by one of its multiple, we may (and do) assume that $x_0$ extends to an element 
%$\wt{x}_0$ 
%of $\CC(S)$. We will show that the N\'eron-Tate height of $x_0$ with respect to a polarisation is an integer and we shall 
%use the fact that $\psi_i(x_0)$ is divisible by $p^i$ in $C_i(L)$ to show that this integer 
%has arbitrarily large divisors and thus it must vanish. From this, we can conclude 
%using a result of Lang (see below), which shows that a point with vanishing N\'eron-Tate height 
%must be a torsion element, up to a contribution of the trace. 

% Hence $\psi_{i,L}(x_0)$ extends to the element $\psi_i(\wt{x}_0)$ of $\CC_i(S)$. 

Now let $D_0$ be a line bundle on $C$. We suppose that 
$[-1]_C^*(D_0)\simeq D_0$ (ie $D_0$ is symmetric), where $[-1]_C$ is the inversion 
morphism given by the group scheme structure of $C$ over $B$. 
We also suppose that $D_0$ is a relatively ample line bundle. 
If $x\in\CC(B)$, write $\tau_x:\CC\to\CC$ for the translation by $x$ morphism. 
We use the same notation for $x\in C(L)$. 

Now consider the isogeny $\phi_{D_0}:C\to C^\vee$ from $C$ to its dual abelian variety,
 which is induced by $D_0$ (this is the polarisation induced by $D_0$). 
 Since $v_0\in\IVD(C)$, we also have \mbox{$\phi_{D_0}(v_0)\in\IVD(C^\vee)$,} since relative Frobenius morphisms are naturally compatible with morphisms of abelian varieties. The point $\phi_{D_0}(v_0)$ corresponds to the 
 line bundle $$M=\tau_{v_0}^*(D_0)\otimes D_0^\vee$$ on $C$ (see 
 \cite[III.13]{Mumford-Abelian}). 
 Since the morphism dual to the Verschiebung morphism is the relative 
 Frobenius morphism (this is very often the definition of the Verschiebung), we see that the fact that $\phi_{D_0}(v_0)\in\IVD(C^\vee)$ 
 translates to the fact that there exist line bundles $M_i$ on 
 $C^{(p^i)}$ for all $i\geq 1$, such that $$F_{C/L}^*(M_1)\simeq M,\,F_{C^{(p)}/L}^*(M_2)\simeq M_1,\,F_{C^{(p^2)}/L}^*(M_3)\simeq M_2,\,\dots$$ 
 Since $\psi_i$ factors by construction through $F_{C^(p^{i-1})/L}\circ F_{C^{(p^{i-2})}/L}\circ\dots 
 \circ F_{C/L},$ we see that for each $i\geq 1$, there is a line bundle $J_i$ on $C_i$
such that $\psi_{i,L}^*(J_i)\simeq M$. 

Now recall that ${D_0}$ extends uniquely (up to isomorphism) to a line bundle 
$\CD_0$ on $\CC$, if we require $D_0$ to be trivial along the unit section of $\CC$ 
(see \cite[Prop. 2.6, p. 21]{MB-Pinceaux}). 
Similarly the line bundle $M$ extends uniquely (up to isomorphism) to a line bundle $\CM$ on $\CC$ with the same property. We shall write $\CJ_i$ for the line bundle 
similarly associated with $J_i$ on $\CC_i$. 
Notice that by unicity, we have $\psi_i^*(\CJ_i)\simeq\CM$. 

We shall now make a height computation. We shall need the 

\begin{lemma}
Let $\CW$ be a line bundle on $\CC$, which is trivial when restricted to the unit section and 
such that $W_L$ is algebraically equivalent to $0$. Let $x\in\CC(B)$. 
Then $\deg(x^*(\CW))$ is the N\'eron-Tate height pairing of $x_L\in C(L)$ and 
$\CW_L$.
\label{lemNT}
\end{lemma}
\beginProof This follows from \cite[III.3.2 and 3.3]{MB-Pinceaux} and the 
definition of polarisations. \endProof

We shall also need the crucial

\begin{prop}
{\rm (a)} There exists a constant $m_0\in\mN^*$ and an infinite set $I_0\subseteq \mN^*$ such that for any $i\in I_0$ and any
$P\in C_i(L)$, the element 
$m_0\cdot P$ extends to an element of $\CC_i(B)$.

{\rm (b)} There is a constant $c_0\in\mN^*$ and an infinite set $I_0\subseteq \mN^*$ such that for any $i\in I_0$ and any 
$P\in\Tor(C_i(L))$ we have $c_0\cdot P=0$.
\label{Uprop}
\end{prop}

We shall prove this proposition later, using Proposition \ref{GUprop} in the Appendix.

%Up to replacing 
%$x_0$ by one of its multiple, we may (and do) assume that $x_0$ extends to an element 
%$\wt{x}_0$ 
%of $\CC(S)$. We will show that the N\'eron-Tate height of $x_0$ with respect to a polarisation is an integer and we shall 
%use the fact that $\psi_i(x_0)$ is divisible by $p^i$ in $C_i(L)$ to show that this integer 
%has arbitrarily large divisors and thus it must vanish. From this, we can conclude 
%using a result of Lang (see below), which shows that a point with vanishing N\'eron-Tate height 
%must be a torsion element, up to a contribution of the trace. 
% Hence $\psi_{i,L}(x_0)$ extends to the element $\psi_i(\wt{x}_0)$ of $\CC_i(S)$.

Let $i\in I_0$. For the next computation, recall that $\psi_{i,L}(v_0)$ is divisible by $p^i$ in 
$C_i(L)$. Let $z_i$ be an element of $C_i(L)$ such that 
$p^i\cdot z_i=\psi_{i,L}(v_0)$. According to Proposition \ref{Uprop} (a),  $m_0\cdot z_i$ extends 
to an element $u_i$ of $\CC_i(B)$. By construction, we have $p^i\cdot u_i=
m_0\cdot \psi_{i}(\wt{v}_0)$. 
 We compute
\begin{eqnarray*}
\deg(([m_0](\wt{v}_0))^*(\CM))&=&\deg(([m_0](\wt{v}_0))^*(\psi_i^*(J_i)))=\deg(([m_0](\psi_i(\wt{v}_0)))^*(J_i))\\&=&
\deg(([p^i](u_i))^*(J_i))=\deg(u_i^*([p^i]^*(J_i)))\\&=&\deg(u_i^*(J_i^{\otimes p^i}))=
p^i\cdot\deg(u_i^*(J_i)).
\end{eqnarray*}
Here $[m_0]$ refers to the multiplication by $m_0$ morphism (in particular $[m_0](\wt{v}_0)=m_0\cdot \wt{v}_0$). 
Now suppose for contradiction that $\deg(([m_0](\wt{v}_0))^*(\CM))\not=0$. 
If we choose $i$ large enough so that $p^i$ is not a divisor of 
$\deg(([m_0](\wt{v}_0))^*(\CM))$ then we obtain a contradiction. 
Thus $\deg(([m_0](\wt{v}_0))^*(\CM))=0$. We may also compute
$$
\deg(([m_0](\wt{v}_0))^*(\CM))=\deg(\wt{v}_0^*([m_0]^*(\CM)))=
\deg(\wt{v}_0^*(\CM^{\otimes m_0}))=m_0\cdot \deg(\wt{v}_0^*(\CM)).
$$
In particular, by Lemma \ref{lemNT}, the N\'eron-Tate height pairing 
of $v_0$ and $M$ vanishes. Now notice that $M$ is by definition 
the image of ${v}_0$ under the polarisation induced by the symmetric ample 
line bundle $D_0$. Hence the N\'eron-Tate pairing of ${v}_0$ and $M$ is twice the N\'eron-Tate height of ${v}_0$ with respect to the polarisation induced by $D_0$. In particular, the N\'eron-Tate height of ${v}_0$ with respect to $D_0$ vanishes. 
By a theorem of Lang (see \cite[Th. 9.15]{Conrad-Trace}) we conclude that the image of 
${v}_0$ in $C(L)$ is an element of finite order. 
Thus the image of $x_0$ in $C(L)$ is also an element of finite order.

Now we show that $x_0\in\Tor^p(C(L))$. For contradiction, suppose that $x_0\not\in\Tor^p(C(L))$. 
We thus may (and do) replace $x_0$ by one of its multiples and suppose that $p\cdot x_0=0$ and $x_0\not=0$. 
We know that $\psi_{i,L}(x_0)$ is divisible by $p^i$ in $C_i(L)$ and since $\psi_{i,L}$ is injective we conclude that 
there is an element of order $p^{i+1}$ in $C_i(L)$ for all $i\geq 1$. Thus contradicts Proposition \ref{Uprop} (b) so we are done.

{\bf Proof of Proposition \ref{Uprop}}. 
We need some preliminaries on moduli spaces of abelian varieties. Let $n\geq 3$ and $g\geq 1$. 
We shall choose particular values for $g$ and $n$ later. 

Let ${\bf A}_{g,n}$ be the functor  from the category of locally noetherian $\mF_p$-schemes to the category of 
sets, such that
\begin{eqnarray*}
&&{\bf A}_{g,n}(B)=\{\,\textrm{\rm isomorphism classes of the following objects :}\\
&&\textrm{principally polarized abelian schemes over $B$ endowed}\\
&&\textrm{with 
a symplectic isomorphism $(\mZ/n\mZ)^{2g}_{B}\simeq\CA[n]$}\,\}
\end{eqnarray*}

D. Mumford proves (see \cite{Mumford-GIT}) that  the functor ${\bf A}_{g,n}$ is representable 
by a scheme, which is separated and of finite type over $\mF_p$. We shall also denote 
this scheme by  ${\bf A}_{g,n}$. 

Furthermore, in \cite[V, 2., Th. 2.5]{FC-Degen}, C. Chai and G. Faltings prove that there exists
\begin{itemize}
\item a scheme $\bar{\bf A}_{g,n}$ (resp. ${\bf A}_{g,n}^*$), which is proper over $\mF_p$;
\item an open immersion ${\bf A}_{g,n}\hookrightarrow\bar {\bf A}_{g,n}$ 
(resp. an open immersion ${\bf A}_{g,n}\hookrightarrow{\bf A}_{g,n}^*$);
\item a semiabelian scheme $\CU$ over $\bar{\bf A}_{g,n}$, such that $\CU_{{\bf A}_{g,n}}$ is isomorphic 
to the universal abelian scheme over ${\bf A}_{g,n}$.
\item a morphism $\bar\pi:\bar{\bf A}_{g,n}\to{\bf A}_{g,n}^*$ compatible 
with the above open immersions of ${\bf A}_{g,n}$;
\item a line bundle $\omega^0$ on ${\bf A}_{g,n}^*$, which is ample 
and such that $\bar\pi^*(\omega^0)=\omega_{\CU/\bar{\bf A}_{g,n}}$.
\end{itemize}
Now write $Z:=B\times_{l_0}{\bf A}_{g,n,l_0}^*$. Recall that the Hilbert scheme $\Hilb(Z/l_0)$ 
is a scheme, which represents the functor
$$
T\mapsto\{\textrm{closed subschemes of $Z_T$, which are 
proper and flat over $T$}\}
$$
from the category of locally noetherian scheme $T$ over $l_0$ to the category of sets. 
It is locally of finite type over $l_0$ (see \cite{FGA-221}). 

Let $\Phi\in\mQ[\lambda]$ be a polynomial 
with rational coefficients and $L_0/Z$ an ample line bundle. By definition, the $l_0$-scheme $\Hilb_\Phi(Z/l_0)$ represents the functor
\begin{eqnarray*}
&&T\mapsto\{\textrm{closed subschemes $W$ of $Z_T$, which are 
proper and flat over $T$}\\
&&\textrm{and such that $\chi(W_t,L_{0,W_t}^{\otimes\lambda})=\Phi(\lambda)$ for all $\lambda\in\mN$ and all $t\in T$}\}\\
\end{eqnarray*}
\vskip-0.8cm from the category of locally noetherian scheme $T$ over $l_0$ to the category of sets. Here $W_t$ is the fibre at $t\in T$ of the morphism $W\to T$ and $L_{0,W_t}$ is the pull-back of $L$ to 
$W_t$ by the natural morphism $W_t\to Z$. The symbol $\chi(\cdot)$ refers to the Euler characteristic. By definition
$$
\chi(W_t,L_{W_t}^{\otimes\lambda})=\sum_{r\geqslant 0}(-1)^r\dim_{\kappa(t)}H^r(W_t,L_{W_t}^{\otimes\lambda}).
$$
(this is called the Hilbert polynomial of $W_t$ with respect to $L_{W_t}$). 
It is shown in \cite{FGA-221}, that $\Hilb_\Phi(Z/l_0)$ is projective over $l_0$ (as a consequence of the 
projectivity of $Z$). Notice that by construction, we have a disjoint union
$$
\Hilb(Z/l_0)=\coprod_{\Phi\in\mQ[\lambda]}\Hilb_\Phi(Z/l_0)
$$ 
Finally, it is shown in 
\cite[part II, 5.23]{Fantechi-FGA} that the functor $\Mor_{l_0}(B,{\bf A}^*_{g,n})$ from 
locally noetherian $l_0$-schemes $T$ to the category of sets, such that 
$$
\Mor_{l_0}(B,{\bf A}^*_{g,n,l_0})(T)=\{\textrm{$T$-morphisms from $B_T$ to ${\bf A}^*_{g,n,T}$}\}
$$
is representable by an open subscheme of $\Hilb(Z/l_0)$. More precisely, 
the natural transformation of functors
$$
\textrm{$T$-morphism $f$ from $B_T$ to ${\bf A}^*_{g,n,T}$}\mapsto\textrm{graph of $f$}
$$
is represented by an open immersion 
$$\Mor_{l_0}(B,{\bf A}^*_{g,n})\hookrightarrow\Hilb(B\times_{l_0}{\bf A}_{g,n,l_0}^*/l_0).$$ 
Let now $D$ be an ample line bundle on $B$.
We choose $L_0$ to be the line bundle $D\boxtimes\omega^0_{l_0}$ 
on $Z=B\times_{l_0}{\bf A}_{g,n,l_0}^*$. 

Recall that the Hodge bundles of the $\CC_i$ all have the same degree by Lemma \ref{IDeglem}. 
Let $d_0:=\deg(\omega_{\CC/B})$ be this common degree. Our aim is to use this to show that all the $C_i$ embed in a bounded family of abelian varieties and apply Proposition 
\ref{GUprop}.

Notice that $C_i[m]\simeq(\mZ/m\mZ)^{2\dim(C_i)}$ and $C^\vee_i[m]\simeq(\mZ/m\mZ)^{2\dim(C^\vee_i)}$. 
Indeed, since $\psi_{i,L}$ is purely inseparable, it induces an isomorphism $C[m]\to C_i[m]$ and thus 
$C_i[m]\simeq(\mZ/m\mZ)^{2\dim(C_i)}$ by (iv) above. For the isomorphism $C^\vee_i[m]\simeq(\mZ/m\mZ)^{2\dim(C^\vee_i)}$, notice that the dual morphism $\psi^\vee_{i,L}:C^\vee_i\to C^\vee$ is separable 
(because its kernel is the Cartier dual of a multiplicative group scheme) and of order a power of $p$. 
Hence, since $(p,m)=1$ it also induces an isomorphism $C^\vee_i[m]\to C^\vee[m]$ (we leave the details to the reader).  

Now let $E_i:=(C_i\times_L C_i^\vee)^4$ . By Zarhin's trick (see 
\cite[IX.1.1]{MB-Pinceaux}) $E_i$ carries a principal polarisation. 
Furthermore, by the last paragraph, we also have $E_i[m]\simeq(\mZ/m\mZ)^{2\dim(E_i)}$. 
Notice also that the identity component of the N\'eron model of $C_i$ is semiabelian, since $\CC_i$ is semiabelian. 
Hence the identity component of the N\'eron model of $C_i^\vee$ is also semiabelian, since 
$C_i^\vee$ is isogenous to $C_i$ (see 
\cite[Prop. 5.8 (4)]{Abbes-RSSC} for a neat presentation). Since the formation of the N\'eron model is compatible 
with products, we conclude that the identity component $\CE_i$ of the N\'eron model of $E_i$ is also 
semiabelian. We also see $\CE_i|_{B\backslash E(\CC)}$ is an abelian scheme 
over $B\backslash E(\CC)$ (where $E(\CC)$ is as in \refeq{canimpinj}). Finally, we have $\deg(\omega_{\CE_i/B})=8\cdot d_0$ by 
\cite[V.3, Lemma 3.4, p. 166]{FC-Degen}.

Let now $g=\dim(E_i)=8\cdot\dim(C)$ and $n=m$. 
By definition, $E_i$ is associated with an $l_0$-morphism $\Spec\,L\to {\bf A}_{g,n,l_0}$. 
By the valuative criterion of properness, this morphism extends to a morphism 
$\phi_i:B\to {\bf A}^*_{g,n,l_0}$ (resp. to a morphism 
$\bar\phi_i:B\to\bar{\bf A}_{g,n,l_0}$). By unicity, we have 
$\bar\pi\circ\bar\phi_i=\phi_i$. Thus, since semiabelian extensions are unique 
(see \cite[IX, Cor. 1.4, p. 130]{Raynaud-Faisceaux}), 
we have $\phi_i^*(\omega^0_{l_0})\simeq\omega_{\CE_i/B}$. The morphism $\phi_i$ is by definition 
associated with an element of $\Mor_{l_0}(B,{\bf A}^*_{g,n,l_0})(l_0)$. 
We can now compute the Hilbert polynomial of the graph $\Gamma_{\phi_i}$ of $\phi_i$ with respect to the line bundle 
$L_0$:
\begin{eqnarray}
\chi(\Gamma_{\phi_i},L_0^{\otimes\lambda})&=&:Q(\lambda)=\chi(B,(D\otimes\phi_i^*(\omega^0_{l_0}))^{\otimes\lambda})=\deg_B((D\otimes\phi_i^*(\omega^0_{l_0}))^{\otimes\lambda})+1-g(B)\nonumber\\&=&
\lambda\cdot\deg_B(D\otimes\phi_i^*(\omega^0_{l_0}))+1-g(B)\nonumber\\
&=&\lambda\cdot\deg_B(D\otimes\omega_{\CE_i/B})+1-g(B)\nonumber\\
&=&
\lambda\cdot\deg_B(D)+\lambda\cdot\deg_B(\omega_{\CE_i/B})+1-g(B)\nonumber\\
&=&
\lambda\cdot\deg_B(D)+\lambda\cdot 8\cdot d_0+1-g(B).
\label{creq}
\end{eqnarray}
Here $g(B)$ is the genus of $B$. The second equality is justified by the Riemann-Roch theorem on $B$. We thus see that the Hilbert polynomial $Q(\lambda)$ of the graph 
of $\phi_i$ with respect to $L_0$ is $Q(\lambda)$ is independent of $i$. 
Thus the element of $\Mor_{l_0}(B,{\bf A}^*_{g,n,l_0})(l_0)$ corresponding to 
$\CE_i$ lies in the scheme 
$$
\Mor_{l_0}(B,{\bf A}^*_{g,n,l_0})(l_0)\cap \Hilb_{Q(\lambda)}(B\times_{l_0}{\bf A}_{g,n,l_0}^*/l_0)
$$
which is of finite type over $l_0$ by the above discussion. 
We now let $Y$ be the Zariski closure in $\Mor_{l_0}(B,{\bf A}^*_{g,n,l_0})(l_0)\cap \Hilb_{Q(\lambda)}(B\times_{l_0}{\bf A}_{g,n,l_0}^*/l_0)$ of the set all the elements of 
$(\Mor_{l_0}(B,{\bf A}^*_{g,n,l_0})(l_0)\cap \Hilb_{Q(\lambda)}(B\times_{l_0}{\bf A}_{g,n,l_0}^*/l_0))(l_0)$ which correspond to some $\phi_i$ ($i\geq 0$). 
Finally we let $H_{00}$ be some irreducible component of $Y$, which 
meets infinitely many such points. Let $\eta_{00}:=\kappa(H_{00})$. By construction, we have an $H_{00}$-morphism 
$$
B\times_{l_0}H_{00}\to{\bf A}^*_{g,n,H_{00}}
$$
which sends $(B\backslash E(\CC))_{\eta_{00}}$ into 
${\bf A}_{g,n,\eta_{00}}\subseteq{\bf A}^*_{g,n,\eta_{00}} $ (because by construction, $(B\backslash E(\CC))_x$ is sent into 
${\bf A}_{g,n,x}$ for a dense sent of points $x\in H_{00}$). 
Let  
$$
\gamma_0:B_{\eta_{00}}\to{\bf A}^*_{g,n,\eta_{00}}
$$
be the induced morphism over $\eta_{00}$. Now recall that there is a proper morphism 
\mbox{$\bar\pi:\bar{\bf A}_{g,n}\to {\bf A}_{g,n}^*$.} By the valuative criterion of 
properness, there is a unique $\eta_{00}$-morphism $\gamma:B_{\eta_{00}}\to 
\bar{\bf A}_{g,n,\eta_{00}}$ such that $\bar\pi_{\eta_{00}}\circ\gamma=\gamma_0$. The morphism $\gamma$ extends over 
an open subset $H_0$ of $H_{00}$, yielding an $H_0$-morphism 
$$
\wt{\gamma}:B\times_{l_0}H_{0}\to\bar{\bf A}_{g,n,H_{0}}.
$$
Replacing $H_{0}$ by one of its open subsets, we may suppose that $H_{0}$ is normal. 
Let now $\CB_0$ be the base change of $\CU$ by $\wt{\gamma}$. 
A theorem of Moret-Bailly (see \cite[VI.3.1]{MB-Pinceaux}) together with a result of Raynaud (\cite[XI.1.4]{Raynaud-Faisceaux}) then shows that $\CB_0$ can be endowed 
with a relatively ample line bundle, which is symmetric and trivial along the zero section. 
Let also $t_0:=l_0$, $C:=B\times_{l_0}H_0$. If we now apply Proposition 
\ref{GUprop} (a) with this choice of $H_0, t_0, C$ and $\CB_0$, we reach the conclusion that there is 
an infinite set $I_0\subseteq\mN^*$ and a constant $n_0$, such that 
for $i\in I_0$, and any $P\in E_i(L)$, the element $n_0\cdot P$ extends to 
an element of $\CE_i(L)$. Since $\CC_i$ is a direct factor of $\CE_i$, we may replace $E_i$ (resp. $\CE_i$) by $C_i$ (resp. $\CC_i$) in the last sentence. This proves (a), with $m_0=n_0$. 
For (b), note that $\Tr_{L|l_0}(E_i)=0$ (since $E_i$ is a product of abelian varieties isogenous to $C$) and apply Proposition 
\ref{GUprop} (b) to the same situation.\endProof

\appendix 

\section{Rational points in families}

The terminology 
of this section is independent of the terminology of the rest of the article and its appendices.

Let $t_0$ be an algebraically closed field. Let $H_0$ be an integral scheme of finite type over $t_0$. Let 
$\pi:C\to H_0$ be a smooth 
curve over $H_0$, with geometrically connected fibres.  
Let $\CB_0$ be a semiabelian scheme over $C$. Suppose that 
there exists a line bundle $L$ on $\CB_0$, which is ample relatively to $C$, symmetric and trivial along the zero section. 
Let $\eta_0:=\kappa(H_0)$ and let $\lambda_0:=\kappa(C)$. 
Note that $\lambda_0$ lies over $\eta_0$ via $\pi$ and that $\lambda_0$ is also the generic 
point of $C_{\eta_0}$ viewed as a subset of $C$. We suppose that $\CB_{0,\lambda_0}$ is 
an abelian variety over $\lambda_0$. 

In the next proposition, we shall need the following lemma, which is well known from the theory of minimal models of curves. 

\begin{lemma}
Let $\phi:X\to Y$ be a morphism of smooth varieties over $t_0$.  
Suppose also that there is a dense open set $Y_1\subseteq Y$, such that 
$\phi|_{Y_1}:\phi^{-1}(Y_1)\to Y_1$ is smooth. Denote by 
$X^\sm$ the maximal open subscheme of $X$, such that $\phi|_{X^\sm}\to Y$ is 
smooth. 

Let $\sigma\in X(Y)$ be a section of $\phi$. Then 
$\sigma\in X^\sm(Y)\subseteq X(Y)$.
\label{lemLiu}
\end{lemma}

\beginProof 
See \cite[Ex. 4.3.25]{Liu-AGAC}.
\endProof

\begin{prop}

{\rm (a)} There is a natural number $n_0$ and a dense open set $V\subseteq H_0$ with the following properties. For any $x\in V(t_0)$, $\CB_{0,\kappa(C_x)}$ is an abelian variety and 
for any $P_x\in \CB_0(\kappa(C_x))$, the point $n_0\cdot P_x\in \CB_0(\kappa(C_x))$ extends to an element of $\CN(\CB_{0,\kappa(C_x)})^0(C_x)$.

{\rm (b)} Suppose that $C$ is proper over $H_0$. Suppose that there is a set $T_0\subseteq H_0(t_0)$, which is dense in $H_0$ and 
such that for any $x\in T_0$ we have $\Tr_{\kappa(C_x)|t_0}(\CB_{0,\kappa(C_x)})=0.$ 
Then there is a dense open set $V\subseteq H_0$ and a natural number $b_{0}$ such that 
for all $x\in V(t_0)$, we have $\#\Tor(\CB_0(\kappa(C_x)))\leq b_0$. 
\label{GUprop}
\end{prop}
Here $\CN(\CB_{0,\kappa(C_x)})^0$ is the connected component of the identity of the N\'eron model $\CN(\CB_{0,\kappa(C_x)})$ of 
$\CB_{0,\kappa(C_x)}$ over $C_x$.

\beginProof We start with (a). 
We shall write $\bar\eta_0$ for an algebraic closure of $\eta_0$. 
Consider the semiabelian scheme $\CB_{0,\bar\eta_0}$ over $C_{\bar\eta_0}$. 
According to \cite[Th. 4.2]{Ku-Proj}, there is an open immersion
\begin{equation}
\CB_{0,\bar\eta_0}\hookrightarrow S_1
\label{I1}
\end{equation} of $C_{\bar\eta_0}$-schemes, with the following properties: $S_1$ is a regular scheme, which is projective 
over $C_{\bar\eta_0}$ and the open immersion 
$\CB_{0,\bar\eta_0}\hookrightarrow S_1$ is an isomorphism when restricted 
to the open subset of $C_{\bar\eta_0}$ over which $\CB_{0,\bar\eta_0}$ is an abelian scheme. In particular $S_1$ is smooth over $\bar\eta_0$, since 
$\bar\eta_0$ is perfect. There is a finite field 
extension $\eta\to\eta_0$ and a morphism 
\begin{equation}
\CB_{0,\eta}\to S
\label{I2}
\end{equation} of $C_{\eta}$-schemes, which is model of \refeq{I1}. 
By flat descent, the morphism $\CB_{0,\eta}\to S$ is also an open immersion and 
$S$ is also smooth over $\eta$ and projective over $C_\eta$. Again by flat descent $\CB_{0,\eta}\to S$ is an isomorphism when restricted 
to the open subset of $C_{\eta}$ over which $\CB_{0,\eta}$ is an abelian scheme. 

We now let 
$g:H\to H_0$ be the normalisation of $H_0$ in $\eta$. Slightly abusing notation, we also denote by $\eta$ the 
generic point of $H$. Note that $g$ is a finite morphism (see eg 
\cite[IV.7.8]{EGA}). We let $\CB$ be the semiabelian scheme on 
$C_H$ obtained by base change and we let $\lambda$ be the 
generic point of $C_H$. Again $\lambda$ lies over $\eta$ via the second projection and is also the generic point of the $C_\eta$. 
By an elementary constructibility argument, there is a non empty open 
set $U\subseteq H$ and an open immersion 
$$
\CB_{C_U}\hookrightarrow\wt{S}
$$
of $C_U$-schemes, where $\wt{S}$ is smooth over $U$ and projective over $C_U.$ 
Furthermore, we may assume that there is an open subset $U'\subseteq C_U$, 
which surjects onto $U$, with the property that $\CB_{U'}$ is an abelian 
scheme over $U'$ and that the induced morphism $\CB_{U'}\hookrightarrow\wt{S}_{U'}$ is an isomorphism.

Let $N_0$ be the supremum of the set of values of the function, which 
associates with any $q\in C_U$ the number of geometric irreducible 
components of the fibre $\wt{S}_q$ of $\wt{S}$ over $q$. This function is 
constructible (see \cite[IV.9.7.9]{EGA}) and so $N_0$ is finite. 

Now let $y\in U(t_0)$. By construction $\CB_{C_y}$ is then a generically abelian 
semiabelian scheme over $C_y$. 
We have a canonical $C_y$-morphism $f:(\wt{S}_{C_y})^\sm\to 
\CN(\CB_{\kappa(C_y)})$ by the definition of the N\'eron model. 
Let $P_y\in\CB(\kappa(C_y))$. The section $P_y$ extends uniquely to a element 
of $(\wt{S}_{C_y})^\sm(C_y)$ by the valuative criterion of properness and Lemma \ref{lemLiu}. It also extends uniquely to an element of $\CN(\CB_{\kappa(C_y)})(C_y)$ by the definition of the N\'eron model. By unicity, these two extensions are compatible with the morphism $f$. 
%Let $m$ be the smallest natural number such that 
%$m\cdot P_y$ extends to an element of \CN(\CB_{\kappa(C_y)})^0$. 
Let $s\in C_y(t_0)$. 
Since the number of irreducible components of $(\wt{S}_{C_y})^\sm_s$ is 
$\leq N_0$, we see that the images of the multiples $P_y,2\cdot P_y,\dots$ 
of $P_y$ in $\CN(\CB_{\kappa(C_y)})(s)$ are contained in at most $N_0$ components 
of $\CN(\CB_{\kappa(C_y)})_s$. Hence the order of the image of 
$P_y$ in the component group of $\CN(\CB_{\kappa(C_y)})_s$ is $\leq N_0$. 
Since $s$ was arbitrary, we see that $N_0!\cdot P_y$ extends to an element of $\CN(\CB_{\kappa(C_y)})^0(C_y)$. 
Note also (for use in (b) below) that since $\CB_{C_y}$ is semiabelian, $\CN(\CB_{\kappa(C_y)})^0(C_y)$ naturally identifies with $\CB_{C_y}$ by the unicity of semiabelian extensions.

Finally let $V$ be the open set $H_0\backslash g(H\backslash U)$. 
By construction, we have $g^{-1}(V)\subseteq U$. Thus 
every point of $V(t_0)$ lifts to a point of $U(t_0)$ (since $g$ is finite) and we see that 
$V$ has the required properties.

For the proof of (b) we first let $U$ be as in the proof of (a).  
We let $\underline{\Sec}^0_{U}(\CB_{C_U}/C_U)$ the functor from locally noetherian 
$U$-schemes $T$ to sets, such that 
$$
 \underline{\Sec}^0_{U}(\CB_{C_U}/C_U)(T)=\{\textrm{sections $\sigma$ of $\CB_{C_T}\to C_T$ 
 such that $\deg((\sigma^*(L))_{C_t})=0$ for all $t\in T$}\}.
 $$
 As $\CB_{C_U}$ is quasi-projective over $U$, this functor is representable by a scheme ${\Sec}^0_{U}(\CB_{C_U}/C_U)$ 
 of finite type over $U$. See eg \cite[Ex. before 5.6.3]{Nitsure-Hilb}. See 
 the proof of Proposition \ref{Uprop} for a similar construction. We leave the details to the reader.  Now let $x\in g^{-1}(T_0)\cap U$. We have an identification 
 \begin{eqnarray*}
 {\Sec}^0_{U}(\CB_{C_U}/C_U)_x(t_0)&=&{\Sec}^0_{x}(\CB_{C_x}/C_x)(t_0)\\&=&\{P\in \CB_{C_x}(C_x)\,|\,\textrm{the N\'eron-Tate height of $P$ with respect to 
 $L_{\CB_{C_x}}$ vanishes}\}
 \end{eqnarray*}
 See \cite[III.3.2 and 3.3]{MB-Pinceaux}. 
 Since $\Tr_{\kappa(C_x)|t_0}(\CB_{0,\kappa(C_x)})=0$, a theorem of Lang (see 
 \cite[Th. 9.15]{Conrad-Trace}) implies that 
 $
 {\Sec}^0_{U}(\CB_{C_U}/C_U)_x(t_0)$ consists of torsion sections. Furthermore, by the Lang-N\'eron theorem, 
 $
 {\Sec}^0_{U}(\CB_{C_U}/C_U)_x(t_0)$ is finite. Hence $
 {\Sec}^0_{U}(\CB_{C_U}/C_U)_x$ is quasi-finite. Since quasi-finiteness is a constructible property (see \cite[IV.9.6.1 (vii)]{EGA}) 
 and $g^{-1}(T_0)\cap U$ is dense in $U$ (because $g$ is finite and $T_0$ is dense in $H_0$), this implies that the scheme $
 {\Sec}^0_{U}(\CB_{C_U}/C_U)$ is quasi-finite over an open subset of $U$. Now replace 
 $U$ by one of its open subschemes so that ${\Sec}^0_{U}(\CB_{C_U}/C_U)$ becomes quasi-finite over $U$. 
 Let $b_{00}$ be an upper bound for the cardinality of the fibres of 
 ${\Sec}^0_{U}(\CB_{C_U}/C_U)\to U$. 
 Using (a), we conclude that we have $$\#(n_0\cdot \Tor(\CB_0(\kappa(x))))\leq b_{00}$$ for all $x\in U(t_0).$
 In particular $b_{00}!\cdot n_0\cdot \Tor(\CB_0(\kappa(x)))$ is the trivial group. Thus 
 by the structure of finite subgroups of abelian varieties, we have 
 $$
 \#\Tor(\CB_0(\kappa(x))))\leq (b_{00}!\cdot n_0)^{2\dim(\CB_{C_U}/C_U)}.
 $$
 and we choose $b_0:=(b_{00}!\cdot n_0)^{2\dim(\CB_{C_U}/C_U)}.$ 
 Finally we let as before $V$ be the open set $H_0\backslash g(H\backslash U)$. 
By construction, we have $g^{-1}(V)\subseteq U$. Thus 
every point of $V(t_0)$ lifts to a point of $U(t_0)$ (since $g$ is finite) and we see that 
$V$ has the required properties.
\endProof

\section{Ampleness of the Hodge bundle and inseparable points}

\label{appendix}

The terminology 
of this section is independent of the terminology of the rest of the article and its appendices.
In this appendix, we shall prove a mild extension of the main result of \cite{Rossler-On-the-group}. 

Let $k$ be a perfect field and let 
$S$ be a geometrically connected, smooth and proper curve over $k$. Let $K:=\kappa(S)$ be its function field. 
Suppose from now on that $k$ has characteristic $p>0$. 

Let $\pi:\CA\to S$ be a smooth commutative group scheme and let $A:=\CA_K$ be the generic fibre of $\CA$. 
Let $\epsilon_{\CA/S}:S\to\CA$ be the zero-section and let $\omega:=\epsilon^*_{\CA/S}(\Omega^1_{\CA/S})$ be the Hodge 
bundle of $\CA$ over $S$. 

%\begin{theor}
%Suppose that $A$ is an ordinary abelian variety. Then $\bar\mu_\min(\omega)\geqslant 0$. 
%\label{mprop2}
%\end{theor}
%In other words, if $A$ is ordinary then $\omega$ is a nef vector bundle.
%\beginProof The proof of \cite[Th. 1.2]{Rossler-On-the-group} goes through verbatim.\endProof

\begin{theor}
Suppose that  $\CA/S$ is semiabelian and that $A$ is an abelian variety.  Suppose that $\bar\mu_\min(\omega)>0$. Then there exists $\ell_0\in\mN$ such the natural injection $A(K^{p^{-\ell_0}})\hookrightarrow A(K^\perf)$ is surjective (and hence a bijection). 
\label{theorA1}
\end{theor}

\NB\, In \cite[Th. 1.1]{Rossler-On-the-group}, Theorem \ref{theorA1} was proven under the assumption that 
$A$ is principally polarised and that $k$ is algebraically closed. In can be shown that 
the condition $\bar\mu_\min(\omega)>0$ is equivalent to the requirement that 
$\omega$ is an ample bundle (see \cite[Introduction]{Rossler-On-the-group} for detailed references).

\beginProof Notice first that in our proof of Theorem \ref{theorA1}, we may replace $K$ by a finite extension field $K'$ without restriction of generality. 
We may thus suppose that $A$ is endowed with an $m$-level structure for some $m\geqslant 3$ with $(m,p)=1$. 

If $Z\to W$ is a $W$-scheme and $W$ is a scheme of characteristic $p$, then for any $n\geqslant 0$ we shall write 
$Z^{[n]}\to W$ for the $W$-scheme given by the composition of arrows $$Z\to W\stackrel{F_W^{n}}{\to}W.$$

Now fix $n\geqslant 1$ and suppose that $A(K^{p^{-n}})\backslash A(K^{p^{-n+1}})\not=\emptyset$. 

Fix $P\in A^{(p^n)}(K)\backslash A^{(p^{n-1})}(K)=A(K^{p^{-n}})\backslash A(K^{p^{-n+1}}).$
The point $P$ corresponds to a  commutative diagram of $k$-schemes
$$
\xymatrix{
& & A\ar[d]\\
\Spec\ K^{[n]}\ar[rru]^{P}\ar[rr]^{\,\,\,\,\,\,\,\,\,\,\,\,\,F_K^n} & &\Spec\ K
}
$$
such that the residue field extension $K|\kappa(P(\Spec\ K^{[n]}))$ is of degree $1$ (in other words $P$ is birational onto its image).  In particular, the  
map of $K$-vector spaces $P^*(\Omega^1_{A/k})\to\Omega^1_{K^{[n]}/k}$ arising from the diagram is non zero. 

Now recall that there is a canonical exact sequence
$$
0\to \pi^*_K(\Omega^1_{K/k})\to \Omega^1_{A/k}\to\Omega^1_{A/K}\to 0.
$$
Furthermore the map $F_K^{n,*}(\Omega^1_{K/k})\stackrel{F_K^{n,*}}{\to}\Omega^1_{K^{[n]}/k}$ vanishes. Also, 
we have a canonical identification $\Omega^1_{A/K}=\pi_K^*(\omega_K)$ (see \cite[chap. 4., Prop. 2]{Bosch-Raynaud-Neron}). Thus the natural surjection 
\mbox{$P^*(\Omega^1_{A/k})\to\Omega^1_{K^{[n]}/k}$} gives rise to a non-zero map 
$$\phi_n=\phi_{n,P}:F_K^{n,*}(\omega_K)\to \Omega^1_{K^{[n]}/k}.$$

The next crucial lemma examines the poles of the morphism $\phi_n$.  

We let $E$ be the reduced closed subset, which is the union of the points $s\in S$, such that the fibre $\CA_s$ is not complete. 

\begin{lemma}
The morphism $\phi_n$ extends to 
a morphism of vector bundles $$F_S^{n,*}(\omega)\to \Omega^1_{S^{[n]}/k}(E).$$
\label{flem}
\end{lemma}
\beginProof (of \ref{flem}).   
First notice 
that there is a natural identification \mbox{$\Omega^1_{S^{[n]}/k}(\log E)=\Omega^1_{{S^{[n]}}/k}(E)$,} because there is a 
sequence of coherent sheaves
$$
0\to\Omega_{{S^{[n]}}/k}\to \Omega^1_{{S^{[n]}}/k}(\log E)\to \CO_E\to 0
$$
where the morphism onto $\CO_E$ is the residue morphism. Here the sheaf $\Omega^1_{{S^{[n]}}/k}(\log E)$ is the 
sheaf of differentials on $S^{[n]}\backslash E$ with logarithmic singularities along $E$. See \cite[Intro.]{Illusie-Reduction} for this result and more details on these notions. 

We may also suppose without restriction of generality that $A$ is principally polarised. Indeed, consider 
the following reasoning.  By Zarhin's trick, the abelian variety $B:=(A\times_K A^\vee)^4$ is principally 
polarised. Also, $B$ can be endowed with an $m$-level structure 
compatible with the given $m$-level structure on $A$, since $A^\vee$ is isogenous to $A$. Let $\CB:=(\CA\times_K \CA^\vee)^4$, where (abusing language) we have 
written $\CA^\vee$ for the connected component of the zero-section of the N\'eron 
model of $A^\vee$. The group scheme  $\CA^\vee$ is also semiabelian, 
since $A^\vee$ is isogenous to $A$ over $K$. 
The morphism $P\times 0\times 0\times\dots\times 0\,(\textrm{seven times})$ gives a point in $B^{(p^n)}(K)$ and there is a commutative diagram
\begin{equation}
\xymatrix{
F_K^{n,*}(\omega_{\CB,K})\ar[rr]^{\phi_{n,P\times 0\times\dots}}\ar[d] & &\Omega^1_{K^{[n]}/k}\\
F_K^{n,*}(\omega_{\CA,K})\ar[rr]^{\phi_{n,P}}&& \Omega^1_{K^{[n]}/k}\ar[u]^{=}
}
\label{PPeq}
\end{equation}
where the vertical arrow on the left is the pull-back map induced by the closed immersion $\lambda\mapsto 
\lambda\times 0\times 0\times\dots\times 0\,(\textrm{seven times})$. Now since $B$ is principally polarised, we know that if 
Lemma \ref{flem} holds for principally polarised abelian varieties, the upper 
row of the diagram \refeq{PPeq} extends to a morphism $F_S^{n,*}(\omega_{\CB})\to \Omega^1_{S^{[n]}/k}(E)$ 
(note that the set of points, where $\CB$ is not complete coincides with the set of points, where $\CA$ is not 
complete). Since $F_S^{n,*}(\omega_\CA)$ is a direct summand of $F_S^{n,*}(\omega_\CB)$, we see that Lemma \ref{flem} holds 
for $A$ if it holds for $B$, thus completing the reduction of Lemma \ref{flem} to the 
principally polarised case.\endProof

The rest of the proof of Theorem \ref{theorA1} is identical word for word with the proof of 
Theorem 1.1 in \cite{Rossler-On-the-group} (from the beginning of the 
proof of Lemma 2.1).\endProof

\section{Specialisation of the Mordell-Weil group}

\label{secRMW}

The terminology 
of this section is independent of the terminology of the rest of the article and its appendices.

In this section, we shall prove a geometric analog of N\'eron's result on 
the specialisation of the generic Mordell-Weil group to a fibre in a family of abelian varieties over number fields (see \cite[chap. 9, Cor. 6.3]{Lang-Fund}). The following results are reminiscent of some 
results proven by Hrushovski in a mixed characteristic context (see \cite{Udi-Manin}) and they are probably already known to many people but we include complete proofs for lack of a reference.

Let $l_0$ be an algebraically closed field. Let $U$ be a smooth and connected 
quasi-projective variety over $l_0$. Let $\CB$ be an abelian scheme over 
$U$. Suppose given an immersion $\iota:U\hookrightarrow\mP^N$ for some 
$N\geq 0$. Let $K$ be the function field of $U$ and let $B:=\CB_K$. 

\begin{prop} Suppose that $\CB(U)$ is finitely generated. For almost all linear subspaces $L\subseteq\mP^N$ of codimension 
$\dim(U)-1$, the intersection $C:=L\cap U$ is smooth, connected, non empty,  
the specialisation map
$$
\CB(U)\to\CB_{C}(C)
$$
is injective and $\Tr_{\kappa(C)|l_0}(\CB_{\kappa(C)})=0$.
\label{redMW}
\end{prop}
Recall that the linear subspaces $L\subseteq\mP^N$ of codimension 
$\dim(U)-1$ are classified by the 
Grassmannian $\Gr(\dim(U)-1,N)$, which is smooth and projective over $l_0$. The words "almost all" stand for "for all the $l_0$-rational points of some dense Zariski open subset of $\Gr(\dim(U)-1,N)$".

Recall that by a theorem of Weil, the restriction map $\CB(U)\to B(K)$ is a bijection. Thus, by the Lang-N\'eron theorem, the condition that $\CB(U)=B(K)$ is finitely generated is equivalent 
to the condition $\Tr_{K|l_0}(B)=0$ . 

For the proof of Proposition \ref{redMW}, we shall need a few lemmata:
\begin{lemma}
Let $N$ be a finite \'etale group scheme over $U$. 
Let $t\in H^1_{\rm et}(U,N)$ and suppose that $t\not=0$. Then for almost all linear subspaces $L\subseteq\mP^N$ of codimension 
$\dim(U)-1$, the intersection $C:=L\cap U$ is smooth, connected, non empty and 
the restriction $t_C\in H^1_{\rm et}(C,N_C)$ of $t$ to $C$ does not vanish.
\label{redetlem}
\end{lemma}
\beginProof 
Let $T\to U$ be a torsor under $N$.  
Note that the torsor $T$ is non trivial iff for all the irreducible components $T'$ of $T$, the (automatically flat and finite) morphism $T'\to U$ has degree $>1$. 
The same remark applies to the restriction of $T$ to a smooth and connected closed subscheme of $U$. 

Let $(T_i)$ be the set of irreducible components of $T$. 

By Bertini's theorem in Jouanolou's presentation (see \cite[p. 89, Cor. 6.11]{Jouanolou-Bertini}), for almost all linear subspaces $L\subseteq\mP^N$ of codimension 
$\dim(U)-1$,

- the intersection $C:=L\cap U$ is smooth, connected and non empty;

and 

- all the $T_{i,C}$ are irreducible.

Let $C$ be in this class. Suppose that $T\to U$ is not trivial. 
By construction, the irreducible components of $T_C$ are the $T_{i,C}$. Since $T_{i,C}\to C$ is flat and finite of the same degree as $T_i\to U$, we see that 
the irreducible components of $T_C$ all have degree $>1$ over $C$. Hence 
the torsor $T_C$ is not trivial.\endProof

\begin{lemma}
Let $N$ be a finite \'etale group scheme over $U$. Suppose 
that $N(U)=0$.  Then for almost all linear subspaces $L\subseteq\mP^N$ of codimension 
$\dim(U)-1$, the intersection $C:=L\cap U$ is smooth, connected, non empty and 
$N_C(C)=0$.
\label{redtorlem}
\end{lemma}
\beginProof
Let $(N_i)$ be the set of irreducible components of $N$, excluding the component of the identity. The condition that $N(U)=0$ is equivalent to the condition that for all $i$, the morphism $N_i\to U$ has degree $>1$. 

As before, by Bertini's theorem, for almost all linear subspaces $L\subseteq\mP^N$ of codimension 
$\dim(U)-1$,

- the intersection $C:=L\cap U$ is smooth and connected;

and 

- all the $N_{i,C}$ are irreducible.

Let $C$ be in this class.  
By construction, the irreducible components of $N_C$ outside of the component of the identity are the $N_{i,C}$. Since $N_{i,C}\to C$ is flat and finite of the same degree as $N_i\to U$, we see that 
the irreducible components of $N_C$ outside of the component of the identity all have degree $>1$ over $C$. Hence $N_C(C)=0$.\endProof

\begin{lemma}
Let $G\subseteq\CB(U)$ be a finite group. 
For almost all linear subspaces $L\subseteq\mP^N$ of codimension 
$\dim(U)-1$, the intersection $C:=L\cap U$ is smooth and connected and 
the reduction map
$$
G\to\CB_{C}(C)
$$
is injective.
\label{injtorlem}
\end{lemma}
\beginProof Left to the reader.\endProof

Finally, we need an elementary but very insightful lemma, due to in essence to N\'eron. 
The following version is due to Hrushovski (see \cite[lemma 1]{Udi-Manin}):

\begin{lemma}[N\'eron-Hrushovski]
Let $r:G\to H$ be a map of abelian groups. Let $l$ be a prime number. 
Suppose that $\Tor_l(H)=0$ and that the induced map $G/lG\to H/lH$ is injective. 
Then $\ker\, r\subseteq \bigcap_{j\geq 0}l^jG$.
\label{HNlem} 
\end{lemma}
\beginProof Let $g\in\ker\, r$. Suppose for contradiction that $g\not\in\bigcap_{j\geq 0}l^jG.$ Let $m\geq 0$ be the smallest  natural number such that $g\not\in 
l^mG.$ Then there is $g'\in G$ such that $l^{m-1}g'=g$ and thus 
$r(g')\in \Tor_l(H)$ so that from the assumptions we have $r(g')=0$. Since the map $G/lG\to H/lH$ is injective, there is $g''\in G$ such that $lg''=g'$. Hence $g=l^m g''$, a contradiction.\endProof

\beginProof (of Proposition \ref{redMW}). Let $l$ be a prime number such that 
$\Tor_l(\CB(U))=0$ and such that $l$ is not the characteristic of $l_0$. Note that 
for any closed subscheme $C$ of $U$, we have an injection 
$\delta_C:\CB(C)/l\CB(C)\hookrightarrow H^1_{\rm et}(C,\ker\,[l]_{\CB,C})$ and this injection is functorial for restrictions to smaller closed subschemes $C_1\hookrightarrow C$. According to Lemmata \ref{redtorlem}, \ref{redetlem} and \ref{injtorlem}, for almost all linear subspaces $L\subseteq\mP^N$ of codimension 
$\dim(U)-1$,

- the intersection $C:=L\cap U$ is smooth and connected;

- the restriction map $H^1(U,\ker\,[l]_{\CB})\to H^1(C,\ker\,[l]_{\CB,C})$ 
is injective on the image of $\delta_U$;

- $(\ker\,[l]_{\CB,C})(C)=0$;

- the restriction map $\Tor(\CB(U))\to\CB(C)$ is injective.

Let $C$ be in this class. 
By construction, the map $\CB(U)/l\CB(U)\to \CB(C)/l\CB(C)$ is injective 
and $\Tor_l(\CB(C))=0$. Let $F$ be a free subgroup of $\CB(U)$, which is a direct summand of $\Tor(\CB(U))$. We have $F\cap\big(\cap_{j\geq 0}l^j\CB(U)\big)=0$ since $\CB(U)$ is finitely generated and $F$ is free. Applying Lemma \ref{HNlem} 
to $G=\CB(U)$ and $H=\CB(C)$, we see that the restriction map 
$F\to\CB(C)$ is injective. Since the restriction map $\Tor(\CB(U))\to\CB(C)$ is also injective, we thus see that the restriction map $\CB(U)\to\CB(C)$ is injective. Finally, we have $\Tr_{\kappa(C)|l_0}(\CB_{\kappa(C)})=0$, for otherwise, we would have $\Tor_l(\CB(C))\not=0$. \endProof

\end{document}